\documentclass[11pt,oneside]{amsbook}

\usepackage{amsmath}
\usepackage{amsthm}
\usepackage{amssymb}
\usepackage{latexsym}
\usepackage{graphicx}

\newtheorem{thm}{Theorem}[chapter]
\newtheorem{lem}[thm]{Lemma}
\newtheorem{prop}[thm]{Proposition}

\newtheorem{cor}[thm]{Corollary}
\theoremstyle{definition}
\newtheorem{df}[thm]{Definition}
\newtheorem{exa}[thm]{Example}
\theoremstyle{remark}
\newtheorem{rem}[thm]{Remark}

\newcommand{\mf}[1]{\mathfrak{#1}}

\newcommand{\N}{\mathbb{N}}
\newcommand{\Z}{\mathbb{Z}}
\newcommand{\Q}{\mathbb{Q}}
\newcommand{\R}{\mathbb{R}}
\newcommand{\CC}{\mathbb{C}}
\newcommand{\HH}{\mathbb{H}}
\newcommand{\Gal}{\mathop{\mathrm{Gal}}}
\newcommand{\diag}{\mathop{\mathrm{diag}}}
\newcommand{\Ker}{\mathop{\mathrm{Ker}}}
\newcommand{\im}{\mathop{\mathrm{Im}}}
\newcommand{\C}{\mathop{\mathbf{C}}}
\newcommand{\M}{M}
\newcommand{\Span}{\mathop{\mathrm{span}}}

\newcommand{\A}{{\mathbb{A}}}
\renewcommand{\P}{{\mathbb{P}}}
\newcommand{\PxP}{\mathop{\mathbb{P}^1\times\mathbb{P}^1}}
\newcommand{\slxsl}{\mathop{\mf{sl}_2\oplus\mf{sl}_2}}
\newcommand{\Pic}{\mathop{\mathrm{Pic}}}
\newcommand{\Aut}{\mathop{\mathrm{Aut}}}
\newcommand{\GL}{\mathop{\mathbf{GL}}}

\newcommand{\PGL}{\mathop{\mathbf{PGL}}}
\newcommand{\gl}{\mathop{\mathfrak{gl}}}
\renewcommand{\sl}{\mathop{\mathfrak{sl}}}
\newcommand{\Int}{\mathop{\mathrm{Int}}}
\newcommand{\Ad}{\mathop{\mathrm{Ad}}}
\newcommand{\ad}{\mathop{\mathrm{ad}}}
\newcommand{\Rad}{\mathop{\mathrm{Rad}}}
\newcommand{\Sym}{\mathop{\mathrm{Sym}}}

\newcommand{\Div}{\mathop{\mathrm{Div}}}
\newcommand{\Der}{\mathop{\mathrm{Der}}}
\newcommand{\g}{\mathfrak{g}}

\newcommand{\keyw}[1]{{\tt #1}}

\def\co{\colon\thinspace}

\input{diagrams.tex}

\begin{document}

\begin{titlepage}
  \vfil
  \ \vskip 35pt
  \begin{flushleft}
    \huge Parametrizing algebraic varieties\\ using Lie algebras
  \end{flushleft}
  \par
  \vskip 10pt
  \hrule height 3pt
  \par
  \vskip 10pt
  \begin{flushright}
    \LARGE Jana P\'\i lnikov\'a \par
  \end{flushright}
  \vskip 70pt
  \large
  \begin{flushright}
    {\Large A dissertation}\\
    submitted in partial fulfillment of the requirements \\
    for the degree of\\
    Doctor of Natural Sciences
    \vskip 10pt
    at
    \vskip 10pt
    Research Institute for Symbolic Computations\\
    Johannes Kepler University\\
    Linz, Austria
  \end{flushright}
  \vskip 50pt
  \begin{flushright}
    Committee in charge:
    \vskip 10pt
    Prof. Dr. Josef Schicho, advisor \\
    Dr. Willem A. de Graaf
  \end{flushright}
  \vskip 50pt
  \begin{flushright}
    Autumn 2006
  \end{flushright}
  \vfil
\end{titlepage}

\pagenumbering{roman}

\ \vskip2.3cm
\begin{center}
{\bf Eidesstattliche Erkl\"arung}
\end{center}
\vskip0.5cm

Ich erkl\"are an Eides statt, dass ich die vorliegende Dissertation 
selbstst\"andig und ohne fremde Hilfe verfasst, andere als die angegebenen 
Quellen und Hilfsmittel nicht benutzt bzw. die w\"ortlich oder 
sinngem\"a\ss\ entnommenen Stellen als solche kenntlich gemacht habe.

\vskip1.9cm
\noindent
Jana P\'\i lnikov\'a

\noindent
Linz, September 2006

\ \newpage
\ \vskip2.3cm
\begin{center}
{\bf Zusammenfassung}
\end{center}
\vskip0.5cm

In dieser Arbeit wird eine neue Methode f\"ur die Parametrisierung
von algebraischen Variet\"aten \"uber einen K\"orper der Charakteristik
Null vorgestellt.
Das Parametrisierungs\-problem wird auf ein Problem von Finden eines 
Isomorphismus von Algebren reduziert.

Wir f\"uhren die Lie-Algebra einer Variet\"at als eine Lie-Algebra,
die mit die Gruppe der Automorphismen der Variet\"at zusammenh\"angt, ein.
Wenn wir einen Isomorphismus von dieser Algebra und irgendeiner klassischen 
Lie-Algebra gestalten 
(zum Beispiel die Algebra von Matrizen mit die Spur Null),
dann k\"onnen wir den Isomorphismus f\"ur die Parametrisierung
der Variet\"at verwenden.
Das Problem des Findens eines Isomorphismus der Lie-Algebren
wird weiter auf die Trivializierung einer assoziativen Algebra reduziert,
d.h.~auf das Finden eines Isomorphismus der gegebenen Algebra und einer Algebra
von Matrizen.
Wenn der Grundk\"orper nicht algebraisch abgeschlossen ist, ist
das letztere ein klassisches Problem aus der Zahlentheorie.
Wir pr\"asentieren Algorithmen zur Trivializierung von Algebren von Grad
bis zu 4 \"uber den rationalen Zahlen und \"uber Zahlk\"orpern.

In der vorliegenden Arbeit wird diese Parametrisierungsmethode auf 
Del-Pezzo-Fl\"achen von Grad 8 und 9 angewandt.
Der Algorithmus ist f\"ur den Grundk\"orper der rationalen Zahlen
implementiert.

\ \newpage
\ \vskip2.3cm
\begin{center}
{\bf Abstract}
\end{center}
\vskip0.5cm

In the thesis we present a new method for parametrizing algebraic varieties 
over the field of characteristic zero.
The problem of parametrizing is reduced to a problem 
of finding an isomorphism of algebras.

We introduce the Lie algebra of a variety as a Lie algebra related 
to its group of automorphisms.
Constructing an isomorphism of this one
and some classical Lie algebra 
(for example the algebra of matrices of the zero trace)
then leads to parametrizing the variety.
The problem of finding an isomorphism of Lie algebras is further reduced to 
trivializing an associative algebra, which means 
finding an isomorphism of the algebra and a full matrix algebra.
The last is a classical problem in number theory, when
regarded over algebraically non-closed fields.
We give algorithms for trivializing algebras of degrees up to 4
over number fields.

In our work we used the method to parametrize Del Pezzo
surfaces of degrees~8 and 9 over number fields.
The algorithms are implemented for the case of the field of the rationals.

\ \newpage
\thispagestyle{empty}
\ \vfill\vfill
\begin{center}
  {\it To the memory of my grandmother.}
\end{center}
\vfill\vfill\vfill

\tableofcontents
 
\ \vskip2.3cm
\thispagestyle{empty}
\begin{center}
{\bf Acknowledgements}
\end{center}
\vskip0.5cm

I want to thank the people who contributed significantly
to the results presented in this thesis: 
Josef Schicho for explaining the problem of parametrization 
and giving a lot of geometrical background and
Willem de Graaf for his work and help with Lie algebras.
I thank also Mike Harrison for his incredible insights and hints
during the work.

A part of the research was carried out at University of Sydney
where I was visiting the Magma group in September 2004.
I am grateful to John Cannon for making the visit possible.

Most of all I would like to express my gratitude to my
advisor Josef Schicho for introducing me to several 
areas of mathematics, sharing interesting research problems 
with me and for all his interest, encouragements and patience.

\ \newpage

\clearpage\setcounter{page}{1}
\pagenumbering{arabic}

\chapter{Introduction}

Finding a parametrization of algebraic varieties is a classical problem in
algebraic geometry. Informally speaking, the task is to find 
for a variety $X\subset\P^n$ 
given as a zero set of a finite collection of polynomials an invertible 
map $\P^d\to X$ for some $d$. The map is in coordinates described by homogeneous
polynomials of the same degree.
As an example we take the unit circle in $\P^2$ given implicitly
as a zero set of the polynomial $x_1^2 + x_2^2 - x_0^2$.
If we find one point $p_0$ on the circle, we can use 
the {\em stereographic projection} and find a parametrization
by considering the points of intersection of $X$ with the lines through $p_0$.
\begin{center}
    \includegraphics{figure1.ps}
\end{center}
A point $p$ on the horizontal axis has coordinates $(s:t:0)$ and the
line through $p$ and $p_0$ is then given by $t x_0 - s x_1 - t x_2$.
The line and the circle meet in points $p_0$ and $\varphi(p)$. 
The coordinates of the second point give a parametrization 
$\P^1\to X\subset\P^2$ of the circle:
\begin{displaymath}
  (s:t)\mapsto(s^2+t^2 : 2st : t^2-s^2).
\end{displaymath}

If a variety given over the field of the rationals $\Q$
can be parametrized over the
algebraic closure $\overline{\Q}$,
one may pose the question whether it is
possible to find a parametrization also over the rationals.
By ``parametrization over $k$'' we mean that the polynomials
describing the map $\P^d\to X$ have their coefficients in $k$.
In our example we were able to parametrize the circle over $\Q$
because the first point $p_0$ had its coordinates in $\Q$.
If we know the method of stereographic projection, then for a nonsingular conic these
two problems become equivalent: we can find a $\Q$-parametrization
of the conic $X$ if and only if we can find a single $\Q$-rational point on $X$.
A similar pattern emerges for many other varieties, including those of higher dimensions.
The stereographic projection is replaced by different methods
(e.g.~using divisor arithmetics) and then the knowledge of a single $\Q$-rational
point again leads to a $\Q$-parametrization.
But finding such a point turns out to be a very hard problem.
Some recent achievements in this research area can be found for example 
in~\cite{lang, tschinkel}.

Let us mention that the {\em Hasse principle} holds for conics.
It means that a conic defined over $\Q$ has a rational point 
(a so-called {\em global solution})
if and only if it has a point over every completion of $\Q$
(i.e.~with respect to any valuation which is possible to define on $\Q$,
a so-called {\em local solution}).
The latter condition basically means that there exists a real solution
to the defining equation of the conic, and also that for every 
prime number $p\in\mathbb{N}$ there is a solution modulo $p$
(or modulo some small power of $p$ in exceptional cases).
In case of conics, the Hasse principle has also a constructive form.
Namely, for a given curve one finds a finite set $\mathcal{S}$ of 
``critical primes''. It means
that if for each prime $p\in\mathcal{S}$ one gets a solution
modulo $p$, then it is possible to combine these local solutions
into a global one, so a rational point is constructed
(cf. e.g.~\cite{cassels}).
Apart from the constructive version of the Hasse principle, 
there are also other efficient methods for finding 
a rational point on a conic (see e.g.~\cite{simon2}).

For all varieties considered in our work, the Hasse principle holds too.
Nevertheless in these cases it gives only a statement on the existence
of a rational point. Though there are methods 
for finding a finite set of ``critical primes'',
still the problem of combining the local
solutions into a global one remains.
We do not pursue this direction in our work.

Here we will deal with some of {\em Del Pezzo surfaces}. 
By definition they are rational (i.e.~parametrizable over 
an algebraically closed field) smooth surfaces 
having their anticanonical divisor ample. 
A basic tool for constructing Del Pezzo surfaces is {\em blowing up} 
the projective plane in one or more points.
To blowup the plane in a point, one removes the point 
and replaces it with the projective line.
A surface obtained by blowing up the projective plane 
in at most 8 points in a general position 
(no three points lie on a line, no six points are on a conic)
turns out to be Del Pezzo.
After blowing up the projective plane in more than 8 points,
the anticanonical divisor of the resulting surface is not ample,
hence the surface is not Del Pezzo.
On the other hand we have that Del Pezzo surfaces come in degrees 1 till 9.
Those of degree 9 are isomorphic to the projective plane
after passing to an algebraically closed field. 
A Del Pezzo surface of degree 8 is either $\P^1\times\P^1$ or 
a blowup of the projective plane in one point.
Finally, Del Pezzo surface of degree $d$, where $1\le d\le 7$,
is a blowup of the projective plane in $9-d$ points in a general position,
see~\cite{manin}.
In particular, every nonsingular cubic surface in $\P^3$ is a blowup
of the projective plane in $6$ points.

The parametrization algorithm which we introduce is demonstrated
on Del Pezzo surfaces of degrees $9$ and $8$.
It follows that all these surfaces are parametrizable over an
algebraically closed field.
Here we give a method which decides whether a given surface has 
a parametrization over $\Q$,
and finds one in the affirmative case.
For blowups in one point, the decision part of the algorithm is not relevant
since blowups are always parametrizable over the field
where they are implicitly defined.
But the constructive part of the method gives an efficient
algorithm for parametrizing these surfaces, 
which is also preferable to known methods, see e.g.~\cite{manin}.
\begin{figure}[!htb]
\begin{center}
\begin{tabular}{c@{\hskip 3cm}c}
  \includegraphics{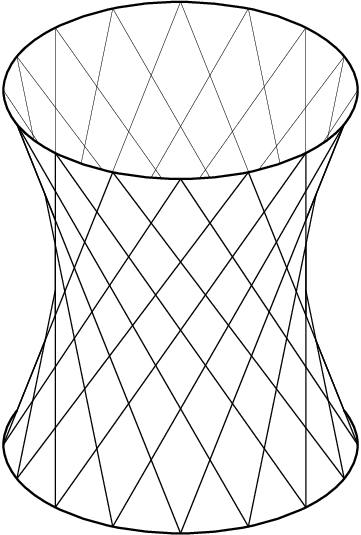}
  &
  \includegraphics{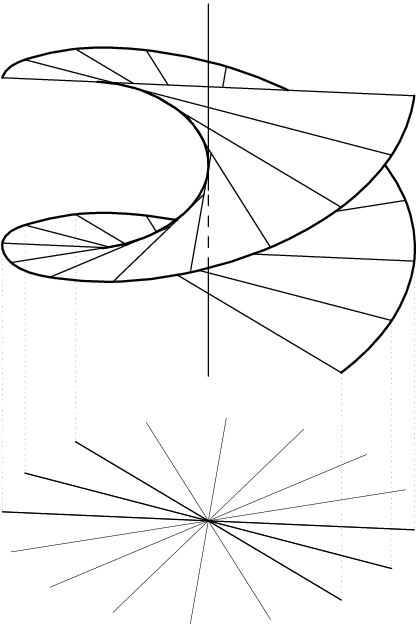}
\end{tabular}\\
  \small
  \begin{tabular}{c}
On the left: The surface $\P^1\times\P^1$ embedded into $\P^3$.\\
On the right: The blowup of $\P^2$ can be imagined as one turn of the helix.
\end{tabular}
\end{center}
\end{figure}

Parametrizing Del Pezzo surfaces of degree 9 and those 
of degree 8 which are embeddings of $\PxP$, 
is a part of the more general problem of 
parametrizing surfaces over the rationals.
In 1998, Schicho~\cite{josef2} gave an algorithm for parametrizing 
surfaces over algebraically closed fields.
In the algorithm one reduces to several base cases which,
except some trivial cases (e.g. $\mathbb{P}^2$), 
are Del Pezzo surfaces of degrees 5 till 9 
and conic fibrations (see also~\cite{Iskovskih:80}). 
Note that rational Del Pezzo surfaces of degrees
smaller than $5$ are mapped in this step to surfaces
of larger degree.
When parametrizing over the rationals, 
these particular cases have to be treated again so that
the properties of the field of the rationals are considered.
The case of conic fibrations was solved in~\cite{josef1}. 
Rational parametrization of Del Pezzo surfaces
of degree~5 is discussed in \cite{Shepperd-Barron:92}, 
degree~6 is dealt with in~\cite{hs}
and degree~7 can be found in e.g. \cite{manin}. 
So this work gives a $\Q$-rational parametrization
of the remaining cases.
The contained material is also presented in~\cite{ghps, gps}.
Thereby the whole problem of parametrizing surfaces over $\Q$ is solved.

Now we give a brief sketch of our method.
We reduce the problem of parametrization to a problem of
trivializing a central simple algebra, i.e.~finding an isomorphism
of an associative algebra and a full matrix algebra.

Let $X'$ denote the variety defined over $k$ which we want to parametrize.
We assume that by some previous analysis we know the isomorphism
class of $X'$. 
More precisely, we know that $X'$ is over the algebraic closure $\overline k$
isomorphic to a ``standard variety'' $Y$ and we know a $k$-parametrization of $Y$.
If we find a $k$-isomorphism between $X'$ and $Y$, we have also constructed
a desired parametrization of $X'$.

Assume that $Y\subset \P^n$ is anticanonically embedded.
The first step in finding a parametrization of $X'$ is to embed it
anticanonically, $\iota\co X'\hookrightarrow\P^n$. Let $X$ denote
the image $\iota(X')$. Then $X'$ and $Y$ are $k$-isomorphic exactly if
$X$ and $Y$ are projectively equivalent over $k$. 

Next we introduce the Lie algebra of a variety as the Lie algebra of the
algebraic group of the projective transformations fixing the variety.
In all examples considered here it is the group of all automorphisms of $X$.
The algebraic group considered in this construction is linear and so
is the Lie algebra as well. 
Hence we have its $(n+1)$-dimensional representation
(provided that $X\subseteq\P^n$).
If $X$ and $Y$ are projectively equivalent over $k$, then their
Lie algebras are two $k$-representations of the same Lie algebra
and the corresponding $(n+1)$-dimensional modules are isomorphic.
We show that under some assumptions also the other implication is true:
if the modules of Lie algebra representations corresponding to 
the varieties $X$ and $Y$ are $k$-isomorphic,
then so are $X$ and $Y$.
Furthermore, we prove that a module isomorphism coincides with
the projective equivalence of the varieties $X$ and $Y$
and hence yields the required parametrization of $X'$.

To find an isomorphism of modules, we first have to construct 
an isomorphism of Lie algebras of $X$ and $Y$. 
If the field $k$ is not algebraically closed, this is a very difficult task.
We transform this problem to finding an isomorphism of associative algebras. 
More precisely, we always reduce to an isomorphism of a given algebra and 
a full matrix algebra. 
Here we give algorithms for algebras up to degree~4.

The algorithms for parametrizing Del Pezzo surfaces described
in the thesis are implemented for $k=\Q$ and included into 
Magma~V2.13,~\cite{magma}.

\vskip0.9cm
\begin{center}
{\bf Notation and terminology}
\end{center}

Throughout the whole work, by $k$ we denote a field of characteristic $0$.
The main focus is to develop algorithms for $k=\Q$.
But except the characteristic being 0 we do not use any special
properties of the field of the rationals,
therefore we stick to the more
general notation ``$k$'' instead of ``$\Q$''.
The algebraic closure of $k$ is denoted by $\overline k$.

The dimension of a vector space $V$ over the field $k$ is denoted $[V:k]$.
The linear span of vectors $v_1,\dots,v_n$
(i.e.~the set $\{\sum_{i=1}^n c_iv_i \mid c_i\in k\}$) 
is denoted by $\Span\{v_1,\dots,v_n\}$ or $\Span_k\{v_1,\dots,v_n\}$.

If $G, H, G_i, \dots$ denote algebraic groups,
then, as is usually done, $\g, \mf{h}, \g_i, \dots$ denote 
their Lie algebras.

If $A$ is an algebra (associative or Lie) over a field $k$
and $k'\supset k$ a field extension, then the algebra $A\otimes_k k'$
obtained from $A$ by enlarging the field of coefficients is also
denoted by $A_{k'}$.

By a {\em twist} of a $k$-algebra $A$ we mean a $k$-algebra $A'$ 
such that there is a field extension $k'\supset k$ with the property
$A_{k'}\cong A'_{k'}$.
For algebraic varieties $X$ and $X'$ defined over $k$ we say analogously 
that $X'$ is a {\em twist} of $X$ if there is an isomorphism
$X\mapsto X'$ defined over an extension $k'$ of $k$.

Sometimes we are slightly sloppy in the language 
of representation theory in the following way.
Though an isomorphism of modules of algebraic groups is defined for
a single group (``two $G$-modules are isomorphic''), we say sometimes
for isomorphic groups $G_1$ and $G_2$ that 
{\em ``$G_1$- and $G_2$-module are isomorphic''}.
By this we mean that there is a fixed isomorphism 
$\varphi\co G_1\to G_2$ of groups
(and from the context it is always clear which $\varphi$ is taken)
and an invertible map $\mu\co V_1\to V_2$ of the modules 
of $G_1$ and $G_2$ respectively,
such that $\mu(gv) = \varphi(g)\mu(v)$.
This situation occurs only if $G_1$ and $G_2$ are 
subgroups of some $\GL_n(k)$ and $V_1$ and $V_2$ respectively are 
their natural modules.
Analogously we use this language for modules of linear Lie algebras.

\chapter{Preliminaries}\label{ch:preliminaries}

This chapter consists mainly of known definitions and facts.
An exception is the analysis of the algorithm {\tt EnvelopingAlgebra} 
in the second section. 

\section{Some notions from algebraic geometry}

By an {\em (algebraic) $k$-variety} $X$ we mean a quasiprojective variety
which is defined by polynomials over $k$.
By $\mathcal{V}(I)$ we mean the variety defined by ideal $I$
and $\mathcal{I}(X)$ denotes the vanishing ideal of $X$.

If $X$ is a $k$-variety (i.e.~defined over $k$) and $k'\supset k$ is a field
extension of $k$,
then $X(k')$ denotes the set of all $k'$-rational points on $X$.

We say that a variety $X$ is {\em $k'$-rational}, if there is a $k'$-rational map
$\P^d\to X(k')$ (called a {\em parametrization}) with a $k'$-rational inverse. 
In the thesis the input variety is always a $\overline k$-rational projective variety, 
almost always a surface.
For a given $\overline k$-rational surface we decide whether it is
also $k$-rational, and construct a parametrization over $k$ 
in the affirmative case.

\subsection{Tangent spaces}

Since the notion of the tangent space is local in its nature, here
we assume that $X\subseteq\A^n$ is an affine variety. 
Let $x\in X$.
The line $L\subset\A^n$ is {\em tangent} to $X$ at $x$,
if its intersection multiplicity with $X$ at $x$ is at least $2$.
The {\em tangent space} of $X$ at $x$ is the geometric locus of all points
on lines tangent to $X$ at $x$. We will denote this space by $\mathcal{T}_x(X)$.

For a point $x\in X$, by $\mf{m}_x$ we denote the maximal ideal of the
local ring $\mathcal{O}_x$ of all functions on $X$ which are regular at $x$.
Let $f\in k[X]$ be the restriction of $F\in k[\A^n]$, $f = F|_X$.
Then $d_x\co k[X]\to\mathcal{T}_x(X)^*$, 
such that $f\mapsto (F - F(x)\pmod{\mf{m}_x^2})|_{\mathcal{T}_x(X)}$, is well defined 
i.e. does not depend on the choice of $F$.
Further we have
\begin{displaymath}
  d_x(f+g) = d_x f + d_x g,\quad
  d_x(fg) = f(x)d_x g + d_x f g(x).
\end{displaymath}
The kernel of $d_x|_{\mf{m}_x}$ is exactly $\mf{m}_x^2$, therefore
$d_x$ is an isomorphism of $\mf{m}_x/\mf{m}_x^2$ and $\mathcal{T}_x(X)^*$.
(cf.~\cite{shafa1}, \S. II.1).
We have an alternative
\begin{df}
  The {\em tangent space} of $X$ at $x$ is defined by 
  $\mathcal{T}_x(X) = (\mf{m}_x/\mf{m}_x^2)^*$.
\end{df}

For more insight let us again consider the map
$d_x\co\mathcal{O}_x\to(\mathcal{T}_x(X))^*$, which is the unique
extension of $d_x$ on $k[X]$ as defined before. Its pullback $d_x^*$ maps
$(\mathcal{T}_x(X))^{**} = \mathcal{T}_x(X)\to\mathcal{O}_x^*$.
If $(a_1,\dots,a_n)\in\mathcal{T}_x(X)$, then by $d_x^*$ this
point is identified with the map
$(a_1,\dots,a_n)^*\circ d_x = \sum_i a_i\frac{\partial}{\partial t_i}|_{x}$
(differentiation followed by evaluation at $x$).
Therefore elements of the tangent space of $X$ at $x$ can be
viewed as {\em point derivations}, i.e.~linear maps $\delta\co\mathcal{O}_x\to k$
such that $\delta(fg) = \delta(f)g(x) + f(x)\delta(g)$.

Let now $X\subseteq\A^n$ and $Y\subseteq\A^m$ be affine varieties
and let $\varphi\co X\to Y$ be a morphism.
Let $x\in X$ and $y = \varphi(x)$.
The pullback $\varphi^*\co k[Y]\to k[X]$ then gives a well defined map
$\mf{m}_y/\mf{m}_y^2\to\mf{m}_x/\mf{m}_x^2$.
Its dual maps $\mathcal{T}_x(X)\to\mathcal{T}_y(Y)$ and is called
the {\em differential of $\varphi$} and denoted by $d_x \varphi$.
Here we give an explicit differentiation formula (see~\cite{humAG}, \S.5.4.)

Let $\varphi$ be given by $m$ coordinate functions $\varphi_i(t_1,\dots,t_n)$, 
$i=1,\dots,m$. Let us identify the tangent spaces at $x$ and $y=\varphi(x)$
with linear subspaces of $k^n$ and $k^m$.
The point $(a_1,\dots,a_n)\in k^n$ represents the point derivation 
$\mathcal{O}_x\to k$ given by $\sum_i a_i\frac{\partial}{\partial t_i}|_x$.
Then its image is
\begin{displaymath}
  d\varphi(a) = (b_1,\dots,b_m)\quad\textrm{where}\quad
  b_k = \sum_i a_i\frac{\partial\varphi_k}{\partial t_i}(x). 
\end{displaymath}

\subsection{Divisors and rational maps}

For an irreducible variety $X$, a collection of irreducible closed subvarieties
$C_1,\dots, C_r$ of codimension 1 in $X$ together with assigned multiplicities
$k_1,\dots,k_r\in\Z$ is called a {\em Weil divisor} on $X$ and denoted 
\begin{displaymath}
  D = k_1C_1 + \dots + k_rC_r.
\end{displaymath}
If all $k_i = 0$, we write $D=0$. 
If all $k_i\ge 0$ and some $k_i>0$, we write
$D>0$ and say that $D$ is {\em effective}.
All Weil divisors form a free abelian group generated by subvarieties
of codimension 1. 
We denote this group by $\Div X$.

A {\em Cartier divisor} (or a {\em locally principal divisor}) is 
a maximal family of pairs $\{(U_i, f_i)\}$ consisting of Zariski open 
subsets $U_i$ of $X$ and rational functions $f_i$ such that
\begin{enumerate}
  \item[(i)]  $\bigcup_i U_i = X$,
  \item[(ii)] for each $i,j$ the rational functions $f_i^{-1}f_j$ and $f_if_j^{-1}$
    are regular on $U_i\cap U_j$.
\end{enumerate}
For a pair $(U,f)$ in the Cartier divisor $D$ we also say that it
represents the divisor $D$ locally (on the open set $U$).
Cartier divisors also form a group: if a divisor $D$ is represented by $(U, f)$
and $D'$ by $(U,f')$, then $D+D'$ is represented by $(U, ff')$.

If the variety $X$ is nonsingular, then the groups of Weil and Cartier 
divisors are isomorphic. In the whole work we deal only with nonsingular
varieties therefore from now on we don't distinguish Weil and Cartier divisor.

Let $f\in k(X)$ be nonzero. Then $f$ defines a divisor
represented by $(U,f)$ for all open $U\subseteq X$ and denoted by $(f)$.
Such divisors are called {\em principal divisors}. 
For divisors $D$ and $D'$ we write $D\sim D'$ if $D-D'$ is principal.
All principal divisors form a subgroup of $\Div X$.
The factor group of all divisors modulo principal divisors
is called the {\em Picard group} of $X$ and denoted by $\Pic X$.

Let $\varphi\co X\to Y$ be a dominant map 
(i.e. the image $\varphi(X)$ is dense in $Y$). 
Since any divisor is principal in some affine subset of $X$, 
the pullback $\varphi^*$ defines a map $\Div Y\to\Div X$.
The pullback of a principal divisor is principal, therefore
we have also a map $\varphi^*\co\Pic Y\to\Pic X$.

If $C$ is a subvariety of $X$ of codimension 1 with the vanishing ideal
$\mathcal{I}(C)\subset k[X]$, 
we can define the valuation $v_C$ on $k[X]\setminus\{0\}$ 
by taking for $v_C(f)$ the largest $n\in\mathbb{N}$ such that 
$(\mathcal{I}(C))^n$ is in the ideal generated by $f$.
The valuation extends to $k(X)\setminus\{0\}$ 
by $v_C(f/g) = v_C(f)-v_C(g)$ for $f,g\in k[X]\setminus\{0\}$.

Let $D$ be a divisor on $X$.
The subset of $k(X)$
\begin{displaymath}
  \mathcal{L}(D) = \{f\in k(X) \mid D + (f) \ge 0\}
\end{displaymath}
together with $\{0\}$ is a linear space, 
since for each irreducible codimension 1 subvariety 
$C\subset X$ and all $f,g\in k(X)$ it holds 
$v_C(f+g) \ge\min\{v_C(f), v_C(g)\}$.
We call $\mathcal{L}(D)$ the {\em Riemann-Roch space of $D$}.
One can prove (see e.g.~\cite{shafa1}) that if $X$ is a projective
variety, then $\mathcal{L}(D)$ is finite-dimensional for all $D$.
In following we will always assume that $[\mathcal{L}(D):k]<\infty$.

If $f_0,\dots,f_n$ is a basis of $\mathcal{L}(D)$,
then we associate with $D$ a rational map
\begin{displaymath}
  \iota_D=(f_0:\dots:f_n)\co X\to\P^n.
\end{displaymath}
Choosing another basis of $\mathcal{L}(D)$ leads to an image of $X$ in $\P^n$
which is projectively equivalent to the former one.

\begin{lem}
  If $D\sim D'$ are divisors on a variety $X$, then
  $\iota_D(X)$ and $\iota_{D'}(X)$ are projectively equivalent.
\end{lem}

\begin{proof}
  Let $g\in k(X)$ be such that $D = D' + (g)$.
  If $f_0,\dots,f_n$ is a basis of $\mathcal{L}(D)$, 
  then $f_0g,\dots,f_ng$ is a basis of $\mathcal{L}(D')$.
  In projective space $(f_0:\dots:f_n)$ is the same map as
  $(f_0g:\dots:f_ng)$, hence the assertion of the lemma follows.
\end{proof}

\begin{lem}
  Let $\varphi\co X\to Y$ be an isomorphism of varieties.
  For a divisor $D'\in\Div Y$ and its pullback $D=\varphi^*(D')\in\Div X$
  then holds that 
  $\iota_D(X)$ and $\iota_{D'}(Y)$ are projectively equivalent.
\end{lem}

\begin{proof}
  If $D' + (g)\ge 0$, then also $\varphi^*(D' + (g)) = D + (\varphi^* g)\ge 0$.
  From $\varphi$ being an isomorphism then follows that
  if $g_0,\dots,g_n$ is a basis of $\mathcal{L}(D')$, then
  $\varphi^*(g_0),\dots,\varphi^*(g_n)$ is a basis of $\mathcal{L}(D)$.
  Then $\iota_D(X) = (\varphi^*(g_0):\dots:\varphi^*(g_n))(X) = 
  (g_0:\dots:g_n)(\varphi(X)) = \iota_{D'}(Y)$.
\end{proof}

We will describe now the very rough idea of our parametrization algorithm.
Let $X$ be a variety that we want to parametrize. Suppose that we know
(by some preceeding analysis) not only that $X$ is $\overline k$-rational
but also that it is isomorphic over $\overline k$ to a variety $X_0$
having a $k$-rational parametrization. Deciding the existence and finding
an isomorphism of $X_0$ and $X$ over $k$ would solve the problem of finding
a $k$-parametrization of $X$.

We say that the divisor $D$ on the variety $X$ 
is {\em very ample}, if $\iota_D$ is an embedding. 
Let $D$ be a very ample divisor such that the class of $D$ in $\Pic X$
is fixed under any automorphism of $X$. 
If $\varphi_1,\varphi_2\co X_0\to X$ are two isomorphisms,
then $\varphi_1^*(D)\sim\varphi_2^*(D)$.
By the two lemmata above it follows that $\iota_{D_0}(X_0)$ and $\iota_D(X)$
are projectively equivalent, where $D_0$ is the pullback of $D$ by any
isomorphism $X_0\to X$.

If $X_0$ and $X$ are $k$-varieties and the embeddings 
$\iota_{D_0}$ and $\iota_D$ are defined over $k$, 
then the problem of deciding existence and finding a $k$-isomorphism
$X_0\to X$ is therefore reduced to finding a projective equivalence over $k$.

In almost all examples considered, we take an anticanonical divisor
as the very ample divisor.
This divisor class is invariant under automorphisms of the variety.
Further, for a $k$-variety there is always an anticanonical divisor over $k$,
therefore the associated embedding is defined over $k$.
The only exception is Example~\ref{exa:sphere}, where 
for the sake of simplicity, another very ample divisor is taken.

\section{Associative and Lie algebras}\label{se:algebras}

\subsection{Structure of semisimple associative algebras}

Let $A$ be a finite-di\-men\-sional $k$-algebra and $V$ an $A$-module.
We say that the module $V$ is {\em semisimple}, if for each submodule 
$V_1\subset V$ there is a submodule $V_2\subset V$ such that
$V = V_1\oplus V_2$.
An algebra is called {\em semisimple}, if it is semisimple as a module over itself.
If there are no nontrivial (two-sided) ideals in $A$, then $A$ is said to be 
{\em simple}.

\begin{thm}[Wedderburn's structure theorem]
  Every semisimple $k$-algebra is a direct sum of finitely many simple algebras 
  and each summand is isomorphic to a matrix algebra over a unique division algebra.
\end{thm}

\begin{proof}
  \cite{jacBA2, pierce}.
\end{proof}

By Wedderburn's theorem the study of semisimple algebras can be reduced to 
study of simple ones. Further, since the center of a simple algebra is a field,
every simple algebra can be regarded as
a {\em central simple algebra $A$}, which by definition is
a finite-dimensional simple $k$-algebra such that the center $\mathbf{C}(A)$ 
is isomorphic to $k$.
For any such algebra $A$ there is $n\in\N$ and 
a unique central division algebra $\Delta$ over $k$ such that
$A$ isomorphic to $\M_n(\Delta)$, the algebra of all matrices over $\Delta$.

\begin{lem}
  If the field $k$ is algebraically closed, then the only division algebra 
  over $k$ is the field $k$ itself.
\end{lem}

\begin{proof}
  If $\Delta$ is a division algebra over $k$, then the minimal polynomial
  of any element in $\Delta$ is irreducible over $k$.
  Only linear polynomials over an algebraically closed field are irreducible,
  therefore each element in $\Delta$ is of the form $\gamma\cdot 1$ for $\gamma\in k$.
\end{proof}

If the algebra $A$ is central simple over $k$, then the algebra 
$A_{k'} = A\otimes_k k'$ obtained from $A$ by extending the coefficient field,
is central simple over $k'$ and its dimension 
$[A_{k'} : k'] = [A : k]$.
Moreover by the previous lemma, if $k'$ is algebraically closed, 
$A_{k'}\cong\M_n(k')$ for some $n\in\N$. 
Consequently, the dimension of any central 
simple algebra is always a square, therefore we can talk about the
{\em degree} of the algebra, which is the integer $\sqrt{[A:k]}$.
The field $k'$ such that $A_{k'} \cong\M_n(k')$ is called
a {\em splitting field of $A$}.

\subsection{Lie algebras}\label{se:LieAssoc}

A {\em Lie algebra} $\g$ is a vector space over $k$ together with 
a (non-associative) bilinear operation 
$[.,.]\co {\g}\times\g\to\g$ satisfying
\begin{enumerate}
  \item[(i)]  $[x,y] + [y,x] = 0$ (anticommutativity),
    \footnote{In case the characteristic of the field equals 2,
    the property (i) is replaced by $[x,x] = 0$ (in other characteristics
    these two are equivalent).}
  \item[(ii)] $[[x,y],z] + [[y,z],x] + [[z,x],y] = 0$ (Jacobi identity).
\end{enumerate}

If $A$ is an associative algebra over $k$, then it can be turned into a Lie algebra
if we replace the associative multiplication by the {\em Lie bracket}
$[x,y] = xy-yx$. We denote this Lie algebra by $A_{\mathrm{Lie}}$.

The Lie algebra is called {\em simple}, if it contains no nontrivial ideal.

Important examples are the Lie algebras $\mf{gl}_n(k)$ and $\mf{sl}_n(k)$, where
$\mf{gl}_n(k) = M_n(k)_{\mathrm{Lie}}$ and $\mf{sl}_n(k)$ is the subalgebra of 
$\mf{gl}_n(k)$ containing all $n\times n$ matrices of zero trace, so 
$[\mf{sl}_n(k) : k] = n^2-1$.
The algebra $\mf{gl}_n(k)$ has the ideal consisting of scalar matrices,
the algebra $\mf{sl}_n(k)$ is simple.

The following theorem about the structure of automorphisms of the Lie algebra 
$\mf{sl}_n(k)$ will be important.

\begin{thm}\label{thm:slAut}
  The group of automorphisms of the Lie algebra of $2\times 2$ matrices
  of trace $0$ is the set of mappings $x\mapsto a^{-1}xa$ for $a\in\GL_2(k)$.

  The group of automorphisms of the Lie algebra of $n\times n$ matrices, $n>2$
  of trace $0$ is the set of mappings $x\mapsto a^{-1}xa$ and $x\mapsto -a^{-1} x^t a$ 
  for $a\in\GL_n(k)$.
\end{thm}

\begin{proof}
  For algebraically closed fields this is Theorem~5 in \cite{jacLA}, Chapter~IX.
  If $k$ is not algebraically closed, let $\alpha\co x\mapsto a^{-1}xa$ for 
  $a\in\GL_n(\overline{k})$ be an automorphism of $\mf{sl}_n(k)$.
  Then $xa = a\alpha(x)$ for each $x\in\mf{sl}_n(k)$ gives
  a system of linear conditions for entries of the matrix $a$,
  where all these conditions are over $k$.
  Therefore $a$ can be taken to be from $\GL_n(k)$.
  The case $x\mapsto -a^{-1} x^t a$ is analogous.
\end{proof}

A Lie algebra homomorphism (i.e. a vector space homomorphism preserving
the multiplication) $\g\to\mf{gl}_n(k)$ is called 
a {\em representation of $\g$}.
The representation $\ad\co\g\to\mf{gl}_d(k)$ where $d = [\g:k]$ such that
$(\ad x) y = [x,y]$ is called the {\em adjoint representation of $\g$}.

Not for every Lie algebra $\g$ there exists an associative algebra $A$ 
such that $\g\cong A_\mathrm{Lie}$.
But to any Lie algebra we can assign an associative algebra as follows.
Let $\g\subset A_\mathrm{Lie}$ for an associative algebra $A$. Then 
the (associative) subalgebra of $A$ generated by $\g$ is called
an {\em enveloping algebra} of $\g$ (in $A$).
The {\em universal enveloping algebra $\mf{U}(\g)$ of $\g$} is the largest
enveloping algebra of $\mf{g}$ in a sense that any enveloping algebra of $\g$
is the universal enveloping algebra modulo an ideal in $\mf{U}(\g)$.

Here we describe and analyze the algorithm which for twists of $\mf{sl}_n(k)$,
(i.e.~such $\g$ that $\g\otimes\bar k\cong\mf{sl}_n(\bar k)$)
constructs an enveloping algebra. 
This is the case needed in our parametrization problem. 

\begin{lem}\label{le:slEnvelop}
  Let $\g$ be a twist of $\mf{sl}_n(k)$ such that there is 
  an $n^2$-dimensional enveloping algebra $A$ of $\g$.
  Then $A$ is central simple.
\end{lem}

\begin{proof}
  Let us fix en embedding $\rho\co\g\hookrightarrow A_\mathrm{Lie}$.
  Then $\rho(\g)$ has no center and hence does not contain the identity.
  Therefore $A_\mathrm{Lie} = \g\oplus\Span\{I\}$.
  Any two-sided ideal in $A$ is an ideal in $A_\mathrm{Lie}$.
  The only ideals of $A_\mathrm{Lie}$ are $0$, $k$ (the center of $A$), 
  $\rho(\g)$ and $A_\mathrm{Lie}$.
  But $k$ is not an ideal in $A$ and 
  $0$ and $A_\mathrm{Lie}$ are trivial ideals. 
  Since $\g$ is simple, the only nontrivial ideal in $A$
  has to coincide with $\rho(\g)$.
  This ideal cannot be nilpotent since otherwise
  (by Engel's theorem, cf.~\cite{humLA})
  also $\g$ would be nilpotent.
  Hence the radical of $A$ is zero. 
  So $A$ is the direct sum of simple ideals. 
  But the center of $A$, which is a direct sum of centers of direct summands, 
  has dimension $1$, therefore $A$ is central simple.  
\end{proof}

\begin{df}
  A {\em Cartan subalgebra} of a Lie algebra $\g$ is a nilpotent subalgebra 
  $\mf{h}\subseteq\g$ which is equal to its normalizer, 
  i.e.~if $[x,y]\in\mathfrak{h}$ for all $x\in\mathfrak{h}$, 
  then $y\in\mathfrak{h}$.
\end{df}

\vskip0.4cm
\noindent
\begin{tabular}{ll}
  {\sc Algorithm:} & {\tt EnvelopingAlgebra}\\
  {\sc Input:}     & $\g$ -- a Lie algebra which is a twist of $\mf{sl}_n(k)$.\\
  {\sc Output:}    & $A$ -- an associative algebra and an embedding 
                    $\g\hookrightarrow A_{\textrm{Lie}}$.
\end{tabular}
\begin{enumerate}
  \item $\mf{h}$ := a Cartan subalgebra of $\g$;
  \item $k'$ := a splitting field of $\mf{h}$, i.e.~a number field containing 
    the eigenvalues of $\ad h$ for $h\in\mf{h}$;
  \item {\tt if} $k' = k$ {\tt then} \\
      \phantom{xx} construct an isomorphism $\g\to\mf{sl}_n(k)$
      (e.g.~according to~\cite{gra6}, \S 5.11);\\
      \phantom{xx} {\tt return} $\M_n(k)$ and $\g\to(\M_n(k))_\textrm{Lie}$\\
    {\tt end if};
  \item construct an isomorphism 
    $\rho'\co\g\otimes k'\to\mf{sl}_n(k')\subset\mf{gl}(V')$
    using the split Cartan subalgebra $\mf{h}\otimes k'$;
  \item viewing the $k'$-space $V'$ of dimension $n$ as a vector space $V$ over $k$
    gives a representation 
    $\rho\co\g\to\mf{sl}_{dn}(k)\subset\mf{gl}(V)$ 
    with $[V:k] = dn$, $d$ being the degree of $k'$ over $k$; 
  \item $\rho(\g)^*$ := the algebra generated by $\rho(\g)$;
  \item {\tt return} $\rho(\g)^*$ and $\g\to(\rho(\g)^*)_\textrm{Lie}$.
\end{enumerate}
\vskip0.4cm

The following assertions show that in case of ``good'' twists 
of $\mf{sl}_n$ our algorithm indeed constructs the smallest possible 
enveloping algebra, and discuss its uniqueness.

\begin{lem}\label{le:LieAssoc}
  The representation $\rho\co\g\to (\rho(\g)^*)_{\mathrm{Lie}}$
  constructed by the algorithm {\tt EnvelopingAlgebra}
  is an injective homomorphism of the Lie algebras.
\end{lem}

\begin{proof}
  Follows from the fact that $\g$ as a twist of $\mf{sl}_n(k)$ is a simple algebra.
\end{proof}

\begin{lem}\label{le:findEnvelop}
  Let $A$ be a twist of $\M_n(k)$.
  Suppose $\g$ is a Lie algebra over $k$ such that $\g\oplus k\cong A_{\mathrm{Lie}}$.
  Let $\rho(\g)^*$ be the enveloping algebra constructed by the algorithm.
  Then $\rho(\g)^*$ is a twist of $\M_n(k)$.
  In particular, $[\rho(\g)^*:k] = n^2$.
\end{lem}

\begin{proof}
  Let $a_0,\dots,a_{n^2-1}$ be a $k$-basis of $A$.
  It can be chosen so that
  after viewing $a_i$'s as $n\times n$ matrices over some extension $k'$ of $k$
  we have that $a_0 = I_n$ is the identity matrix 
  and $a_i$ for $i\ge 1$ are matrices of zero trace.
  The structure constants in $A$ with respect to this basis are in $k$:
  \begin{displaymath}
    a_i a_j = \sum_k c_{ij}^k a_k,\quad c_{ij}^k\in k.
  \end{displaymath}
  From the assumption there is an isomorphism of Lie algebras 
  $\g\oplus k\to A_{\mathrm{Lie}}$. 
  Let $\varphi$ denote the restriction of this map to $\g$.
  After extending the field of coefficients to $k'$ we have
  $\varphi(\g\otimes k') = \mf{sl}_n(k')$, the linear span of 
  $a_1,\dots, a_{n^2-1}$.
  Let $b_1,\dots b_{n^2-1}$ be the basis of $\g$ such that 
  $a_i = \varphi(b_i)$.

  Let $\rho'\co(\g\otimes k')\to\mf{sl}_n(k')$ be the representation (over $k'$) 
  constructed in the step (4) of the algorithm.
  We have to show that the (associative) structure constants in $A$
  with respect to the basis $a_0,\rho'(b_1),\dots,\rho'(b_{n^2-1})$ are also in $k$.

  Denote $\alpha = \rho'\circ\varphi^{-1}$, 
  so $\alpha$ is an automorphism of $\mf{sl}_n(k')$.
  By Theorem~\ref{thm:slAut} then $\alpha$ either maps
  $x\mapsto a^{-1}xa$ or $x\mapsto -a^{-1}x^ta$.
  In the first case we have
  \begin{displaymath}
    \rho'(b_i)\rho'(b_j) = (\alpha\circ\varphi(b_i))(\alpha\circ\varphi(b_j)) = 
    \alpha(a_i)\alpha(a_j) = \dots = 
    c_{ij}^0 a_0 + \sum_{k=1}^{n^2-1} c_{ij}^k \rho'(b_i).
  \end{displaymath}
  Similarly in the second case we get 
  \begin{displaymath}
    \rho'(b_i)\rho'(b_j) = -c_{ji}^0 a_0 \sum_{k=1}^{n^2-1} (-c_{ji}^k) \rho'(b_i).
  \end{displaymath}
  Finally, the $k$-dimension of $\rho(\g)^*$ is equal 
  to the $k$-dimension of $\rho'(\g)^*$.
\end{proof}

Note that the assumption $\g\oplus k\cong A_{\mathrm{Lie}}$ is stronger then 
$\g$ being a twist of $\mf{sl}_n(k)$. 
Indeed, if $\g$ is the Lie algebra over $\Q$ of zero-trace skew Hermitian matrices 
of degree 3 (i.e. all matrices $x\in\M_3(\Q(i))$ such that $x^t + \overline x = 0$),
then $\g$ splits over $\Q(i)$, i.e. $\g\otimes\Q(i)\cong\mf{sl}_3(\Q(i))$.
On the other hand there is no twist of $\M_3(\Q)$ having splitting field of degree $2$.
In this case one finds out that $[\rho(\g)^*:k] = 18$, though $n=3$.

\begin{prop}\label{prop:uniqueEnvelop}
  Let $\g$ be a twist of $\mf{sl}_n$ such that there is an embedding 
  $\g\hookrightarrow A$, where $A$ is a twist of $\M_n(k)$.
  Let $A_1,A_2$ be two enveloping algebras of $\g$, both of dimension $n^2$.
  Then $A_1$ and $A_2$ are isomorphic or antiisomorphic,
  where the second possibility can occur only if $n > 2$.
\end{prop}

\begin{proof}
  Firstly, by Lemma~\ref{le:slEnvelop}, 
  any $n^2$-dimensional enveloping algebra of $\g$ is central simple.
  Therefore both $A_1$ and $A_2$ are twists of $M_n(k)$.

  Let now $\rho_i$ be a fixed embedding 
  $\g\hookrightarrow (A_i)_\mathrm{Lie}$, $i=1,2$.
  Since $A_i$ is a twist of $\M_n(k)$, for each $x\in\g$ the image 
  $\rho_i(x)$ can be understood as a matrix in $\M_n(\overline k)$,
  the associative multiplication in $A_i$ being then just a multiplication 
  of matrices.

  Let further $\varphi\co A_1\to A_2$ be a linear map such that
  $\varphi([a,b]) = [\varphi(a),\varphi(b)]$. 
  Such a map exists because the Lie algebra isomorphism
  $\rho_2\circ\rho_1^{-1}$ easily extends to $\rho_1(\g)\oplus I_n$.
  Clearly $\varphi(\rho_1(\g)) = \rho_2(\g)$, so $\varphi$ restricts
  to a Lie algebra isomorphism $\rho_1(\g)\to\rho_2(\g)$.
  Then $\varphi$ as a map from $A_1$ uniquely extents to a linear
  map $A_1\otimes\overline{k} = \M_n(\overline{k}) \to \M_n(\overline{k})$.
  Restriction to the matrices of trace 0 leads to a Lie algebra
  automorphism of $\mf{sl}_n(\overline{k})$.
  By Theorem~\ref{thm:slAut} there are two possibilities:
  (1) $\varphi$ maps $x$ to $a^{-1}x a$ for some regular $a$.
  Then $\varphi$ is an automorphism of $\M_n(\overline{k})$,
  i.e. an isomorphism $A_1\to A_2$.
  (2) $\varphi(ab) = -\varphi(b)\varphi(a)$, therefore
  $\psi\co A_1\to A_2$, $x\mapsto -\varphi(x)$ is an 
  antiisomorphism of algebras.
\end{proof}

\begin{cor}\label{cor:LieAsocIso}
  The Lie algebra $\g$ is isomorphic to $\mf{sl}_n(k)$ 
  if and only if $\rho(\g)^*$ is isomorphic to $\M_n(k)$.
\end{cor}

\begin{proof}
  There exists an embedding $\rho\co\mf{sl}_n(k)\to\M_n(k)_\mathrm{Lie}$.
  Since there is an antiisomorphism $\M_n(k)\to\M_n(k)$ (e.g.~$x\mapsto x^t$),
  the claim follows. 
\end{proof}

\subsection{Representations of semisimple Lie algebras}

For a Lie algebra $\g$ we define the sequence of ideals $\g^{(i)}$
(called the {\em derived series}) by 
$\g^{(0)} = \g$, $\g^{(i)} = [\g^{(i-1)}, \g^{(i-1)}]$.
The algebra $\g$ is called {\em solvable}, if $\g^{(n)} = 0$ fore some $n$.
The unique maximal solvable ideal in $\g$ is called the {\em radical of $\g$}
and denoted by $\Rad\g$.
If $\Rad\g = 0$, then $\g$ is called {\em semisimple}.
Clearly, a simple algebra is also semisimple.

In following we assume that $\g$ is a semisimple Lie algebra
with a fixed split Cartan subalgebra $\mf{h}$
(so all eigenvalues of $\ad h$ for all $h\in\mf{h}$ are in the field $k$).
It follows that $\ad\mf{h}$ is diagonalizable over $k$.
Therefore $\g$ is a direct sum of the subspaces
$\g_\alpha = \{x\in\g\mid[h,x] = \alpha(h)x\textrm{ for all }h\in\mf{h}\}$,
where $\alpha$ ranges over $\mf{h}^*$.
The space $\g_0$ is the centralizer of $\mf{h}$ in $\g$ and it is equal to $\mf{h}$.
The set of all nonzero $\alpha\in\mf{h}^*$ such that $\g_\alpha\ne 0$
is denoted by $\Phi$. 
The elements of $\Phi$ are called {\em roots of $\g$} relative to $\mf{h}$.
By $\Phi^+$ resp.~$\Phi^-$ we denote the set of positive resp.~negative roots.

If $V$ is any $\g$-module, then $\mf{h}$ acts diagonally on $V$:
$V = \coprod V_\lambda$, where $\lambda$ ranges over $\mf{h}^*$ and
$V_\lambda = \{v\in V\mid h\cdot v = \lambda(h)v\}$.
If $V_\lambda\ne 0$, we call it {\em weight space} and $\lambda$ is called
{\em weight of $V_\lambda$}.

For a $\g$-module $V$, each $\g_\alpha$ for $\alpha\in\Phi$ 
maps weight spaces to the weight spaces, 
more precisely $\g_\alpha$ maps $V_\lambda$ to $V_{\lambda+\alpha}$.
A nonzero vector $v\in V_\lambda$ for some $\lambda$ killed by all 
$\g_\alpha$, $\alpha\in\Phi^+$ is called {\em a highest weight vector} and denoted
by $v^+$. Then $\lambda$ is called a {\em highest weight}.

\begin{lem}
  Let $V$ be an irreducible finite dimensional $\g$-module, 
  with a highest weight vector $v^+$.
  Then $v^+$ is unique up to (nonzero) scalar multiples.
\end{lem}

\begin{proof}
  See~\cite{humLA}, Corollary to Theorem~20.2.
\end{proof}

\begin{thm}
  Let $V$, $W$ be finite dimensional $\g$-modules of the same highest weight $\lambda$.
  If $V$ and $W$ are irreducible, then they are isomorphic.
\end{thm}

\begin{proof}
  See~\cite{humLA}, Theorem~A in \S 20.3.
\end{proof}

If $v^+$ is the heighest weight vector of an irreducible finite dimensional 
module $V$, then $V = \mf{U}(\g)\cdot v^+$.
Actually, $V$ is spanned by vectors $x_{\alpha_1}\dots x_{\alpha_m} v^+$, 
where $\alpha_1,\dots,\alpha_m\in\Phi^-$. 
This gives a very fast algorithm for constructing an isomorphism of modules:
\begin{enumerate}
  \item Fix a split Cartan subalgebra $\mf{h}$ in $\g$.
  \item Check whether the two $\g$-modules $V$ and $V'$ have the same
    dimension and the same highest weight.
  \item Set $v^+$ be the highest weight vector in $V$ and
    $v'^+$ the highest weight vector in $V'$.
  \item Take $v_1,\dots,v_n$ to be a basis of $V$ such that 
    the action of $\mf{h}$ is diagonal with respect to this basis.
  \item For each $v_i$ take a chain $x_{\alpha_1}\dots x_{\alpha_m}\in\mf{U}(\g)$
    such that $v_i = x_{\alpha_1}\dots x_{\alpha_m} v^+$.
    Then $v_i$ is by the module isomorphism mapped to 
    $x_{\alpha_1}\dots x_{\alpha_m} v'^+$.
\end{enumerate}

The last notion to be defined is
a {\em Chevalley basis} of a semisimple algebra $\g$. 
It is a basis $\{x_\alpha, \alpha\in\Phi;\ h_i, i=1\dots m\}$ such that
$\{h_i, i=1\dots m\}$ is a basis of a split Cartan subalgebra of $\g$
and $x_\alpha$ are root vectors satisfying
\begin{itemize}
  \item[(i)]  $[x_\alpha,x_{-\alpha}]$ is a $\mathbb{Z}$-linear 
    combination of $h_1,\dots,h_l$,
  \item[(ii)] if $\alpha,\beta,\alpha+\beta\in\Phi$, 
    $[x_\alpha,x_\beta] = c_{\alpha,\beta}x_{\alpha+\beta}$, 
    then $c_{-\alpha,-\beta} = -c_{\alpha,\beta}$.
\end{itemize}
It holds that the structure constants of any algebra $\g$ 
with respect to its Chevalley basis are integral (cf. e.g.~\cite{humLA}).

\section{Algebraic groups and their Lie algebras}

\begin{df}
  An {\em algebraic group} is an algebraic variety $G$ together with
  \begin{enumerate}
    \item[(i)]   an element $e\in G$,
    \item[(ii)]  a morphism $\mu\co G\times G\to G$, $(g,h)\mapsto gh$,
    \item[(iii)] a morphism $\iota\co G\to G$, $g\mapsto g^{-1}$.
  \end{enumerate}
  with respect to which $G$ is a group.
  We call $G$ a {\em $k$-group}, if $G$ is a $k$-variety
  and morphisms $\mu$ and $\iota$ are defined over $k$.

  A {\em morphism} of algebraic groups is a morphism of algebraic varieties
  which is also a homomorphism of groups.
\end{df}

Since an algebraic group has a transitive group of automorphisms (namely
$g\mapsto gh$) and the singular locus is a proper closed subset of the variety, 
an algebraic group is a non-singular variety.

For a finite-dimensional space $V$ over $k$, the group $\GL(V)$ of all
automorphisms of $V$ is an algebraic $k$-group since it is a principal open
subset of $\M_d(k)$ (given by the $k$-polynomial equation $\det g\ne 0$) 
and both multiplication 
and inverse of a matrix are polynomial maps with integer coefficients.
A closed subgroup of $\GL(V)$ is called a {\em linear algebraic group}.
A morphism $G\to\GL(V)$ of algebraic groups is called 
a {\em rational (linear) representation} of $G$.

With an algebraic group $G$ we can associate a {\em Lie algebra} $L(G)$ 
which is the tangent space $\mathcal{T}_e(G)$ of $G$ at $e$ endowed 
with the multiplication as follows.

Let $A$ be a commutative $k$-algebra and $M$ be an $A$-module.
A {\em $k$-derivation} from $A$ to $M$ is a $k$-linear map $\delta\co A\to M$
such that $\delta(ab) = \delta(a)b + a\delta(b)$ for all $a,b\in A$.
The set of all $k$-derivations from $A$ to $M$ is denoted $\Der_k(A,M)$.

Then the {\em Lie algebra of $G$} is the space of all left invariant derivations 
in $\Der_k(k[G],k[G])$, i.e.
\begin{displaymath}
  L(G) = \{\delta\in{\Der}_k(k[G],k[G]) \mid 
    \delta\lambda_g = \lambda_g\delta\ \forall g\in G\}
\end{displaymath}
where $\lambda_gf(h) = f(g^{-1}h)$ for all $f\in k[G]$, $g,h\in G$
is a {\em left translation}.

We can relate the Lie algebra of an algebraic group $G$ with the tangent space
of the group at the identity $e$ 
(we recall that $\mathcal{T}_e(G) = (\mf{m}_e/\mf{m}_e^2)^*$):

\begin{thm}
  Let $\Theta\co L(G)\to\mathcal{T}_e(G)$ be defined by
  $(\Theta\delta)(f) = (\delta f)(e)$.
  Then $\Theta$ is a vector space isomorphism.
\end{thm}

\begin{proof}
  See~\cite{humAG}, Theorem~9.1.
\end{proof}

For an algebraic group $G$, we denote $\mf{g} = \mathcal{T}_e(G)$ together
with the Lie algebra structure carried over by $\Theta$.

The Lie algebra $L(\GL(V))$ with $[V:k]=n$ is $\mf{gl}_n(k)$, the set of all
$n\times n$ matrices, where $[x,y] = xy-yx$.
If $G\subseteq\GL(V)$ is a linear algebraic group, then $L(G)$ is a subalgebra
of $\mf{gl}_n(k)$. For more details see~\cite{borel, humAG, springer}.

\begin{prop}\label{prop:lieAlgIso}
  Let $G_1,G_2\subseteq\GL_m(k)$ be algebraic groups and 
  $\g_1,\g_2\subseteq\mf{gl}_m(k)$
  their Lie algebras.
  If $h\in\GL_m(k)$ is such that $g\mapsto hgh^{-1}$
  is an isomorphism from $G_1$ to $G_2$,
  then its differential is an isomorphism from $\g_1$ to $\g_2$ given by
  $x\mapsto hxh^{-1}$.
\end{prop}

\begin{proof}
  Let us denote $\varphi:G_1\mapsto G_2$, $g\mapsto hgh^{-1}$
  We have 
  \begin{displaymath}
    \varphi(g)_{ij} = (hgh^{-1})_{ij} = \textstyle\sum_{kl}h_{ik}g_{kl}(h^{-1})_{lj}.
  \end{displaymath}
  By differentiating $(i,j)$-th coordinate of $\varphi$ we get
  \begin{displaymath}
    d\varphi(x)_{ij} = d(\textstyle\sum_{kl}h_{ik}g_{kl}(h^{-1})_{lj})(x)
    = \textstyle\sum_{kl}h_{ik}x_{kl}(h^{-1})_{lj}.
  \end{displaymath}
  So $d\varphi(x) = hxh^{-1}$.
\end{proof}

\begin{cor}\label{cor:moduleIso}
  Let $\alpha_i\co G\to G_i\subseteq\GL_m(k)$, $i=1,2$ be two faithful 
  rational representations of the algebraic group $G$ such that the two
  $m$-dimensional modules $V_1$, $V_2$ are isomorphic 
  ($V_i$ is the natural $G_i$-module).
  Then the isomorphism of $G$-modules $V_1$ and $V_2$ is also
  an isomorphism of $V_1$ and $V_2$ as $\g$-modules.
\end{cor}

\begin{proof}
  Let $\mu$ denote the $G$-module isomorphism $V_1\to V_2$
  and let $h$ be the matrix of $\mu$, so $\mu\co v\mapsto hv$.

  Let $\varphi\co G_1\to G_2$ be an isomorphism of algebraic groups 
  such that $\mu(g v) = \varphi(g)\mu(v)$.
  Then for any $g\in G_1$ and $v\in V_1$ we have
  \begin{displaymath}
    h g v = \mu(gv) = \varphi(g)\mu(v) = \varphi(g) h v.
  \end{displaymath}
  Since $v\in V_1$ is arbitrary, it follows $\varphi(g) = hgh^{-1}$ 
  for all $g\in G_1$.
  Now by Proposition~\ref{prop:lieAlgIso}, the isomorphism $d\varphi\co\g_1\to\g_2$
  maps $x\mapsto hxh^{-1}$. 
  It is straightforward to check that then $v\mapsto hv$
  is a corresponding isomorphism of the $\g$-modules. 
\end{proof}

The $G\subset\GL_m(k)$ acts on itself by inner automorphisms
\begin{displaymath}
  \Int h\co g\mapsto hgh^{-1},\quad g,h\in G.
\end{displaymath}
The differential of $\Int h$ is denoted by $\Ad h$. 
By Proposition~\ref{prop:lieAlgIso} we have that 
\begin{displaymath}
  \Ad h(x) = hxh^{-1},\quad h\in G, x\in\mathfrak{g}.
\end{displaymath}

\begin{prop}\label{prop:LieConjug}
  Let $G_1\subseteq\GL_n(k),G_2\subseteq\GL_m(k)$ be algebraic groups and 
  $\mf{g}_1\subseteq\mf{gl}_n(k),\mf{g}_2\subseteq\mf{gl}_m(k)$ their Lie algebras.
  Let $\varphi\co G_1\to G_2$ be a morphism and 
  $d\varphi\co\mf{g}_1\to\mf{g}_2$ its differential.
  Then for every $h\in G_1(k)$ we have
  \begin{displaymath}
    d\varphi\circ\Ad h = \Ad\varphi(h)\circ d\varphi.
  \end{displaymath}
\end{prop}

\begin{proof}
  Write $\varphi = (\varphi_{ij})_{i,j = 1}^m$, so 
  $d\varphi = (d\varphi_{ij})_{i,j = 1}^m$, where
  \begin{displaymath}
    d\varphi_{ij}(x) = 
      \sum_{k,l} x_{kl}\frac{\partial}{\partial t_{kl}}\varphi_{ij}(g)\big|_{g=e}.
  \end{displaymath}
  Further we have
  \begin{displaymath}
    \varphi_{ij}(hgh^{-1}) = (\varphi(h)\varphi(g)\varphi(h^{-1}))_{ij} =
    \sum_{u,v}\varphi_{iu}(h)\varphi_{uv}(g)\varphi_{vj}(h^{-1}).
  \end{displaymath}
  Putting these two together gives
  \begin{displaymath}
    d\varphi_{ij}(hxh^{-1}) = (\varphi(h)d\varphi(x)\varphi(h^{-1}))_{ij}
  \end{displaymath}
  and the assertion of the proposition follows.
\end{proof}

\chapter{Automorphisms of a variety}\label{ch:aut}

This chapter gives a theoretical ``basement'' for the parametrization method.
We introduce the Lie algebra of a variety and give an algorithm for computing it.
We will see that $k$-isomorphic varieties have $k$-isomorphic Lie algebras.
We analyze also the other implication: 
If two $\overline k$-isomorphic varieties have $k$-isomorphic Lie algebras,
are then these varieties isomorphic also over $k$?
And how can an isomorphism of Lie algebras be used
for finding an isomorphism of the varieties?

\section{The Lie algebra of the variety}

Let $X\subseteq\P^n$ be an embedded projective variety 
so that $\mathcal{I}(X)\subset k[x_0,\dots,x_n]$. 
We define the following linear group:
\begin{displaymath}
  G(X,k) = \left\{g\in{\GL}_n(k) \mid f(gp) = 0
  \ \forall p\in X\ \forall f\in\mathcal{I}(X)\right\}.
\end{displaymath}
Then the projectivization $\P G(X,k) = G(X,k)/kI_n$,
where $I_n$ is the identity matrix, are exactly 
the automorphisms of $\P^n$ fixing $X$. 
Since $G(X,k)$ is described by polynomials, it is an algebraic group, 
so it has a Lie algebra which we denote $\g(X,k)$.

\begin{df}\label{df:embeddedLA}
  Let $X\subseteq\P^n$ and let $\g(X,k)$ be the Lie algebra of $G(X,k)$.
  We define $\g_0(X,k) = \g(X,k)/kI_n$, where $I_n\in\gl_n(k)$
  is the identity matrix.
  We call $\g_0(X,k)$ {\em the Lie algebra of the (embedded) variety} 
  $X(k)$. 
\end{df}

Of course $\g_0(X,k)$ by this definition depends on the embedding of $X$.
Moreover $\P G(X,k)$ can be a very small subgroup of all automorphisms of $X$:

\begin{exa}
  Let $C = \mathcal{V}(x^3+y^3+z^3 - 3mxyz) \subset\P^2$ with $m^3\ne 1$,
  be a curve over $k = \overline\Q$.
  Then $C$ is an elliptic curve, so 
  $\{\alpha_p\co C\to C$, $q\mapsto pq\mid\forall p\in C\}$ 
  is a continuous subgroup of its automorphisms.
  Nevertheless, there are only finitely many linear transformations of $\P^2$ 
  fixing $C$, as each such transformation has to fix the set of 9 flexes of $C$.
  Finally, since $G(C,k)$ is discrete, the Lie algebra of $C$ is trivial.
\end{exa}

By $\Aut(X)$ we denote the group of automorphisms of $X$.
To be able to study the structure of $G(X,k)$, 
we next introduce a subgroup of $\Aut(X)$, which is independent on
the embedding of $X$:
\begin{displaymath}
  {\Aut}_0(X) = \left\{\varphi\in\Aut(X) \mid
  \varphi^* \textrm{ acts trivially on }\Pic(X)\right\}.
\end{displaymath}

\begin{thm}
  Let $X$ be embedded into $\P^n$ by a map associated
  to a very ample divisor $D$ on $X$.
  Then $\Aut_0(X)$ is embedded into $\P G(\iota_D(X),k)$
  as an algebraic variety.
\end{thm}

\begin{proof}
  Let $D$ be a very ample divisor, $f_0,f_1,\dots,f_n$ a basis of the
  Riemann-Roch space $\mathcal{L}(D)$, so $(f_0:f_1:\dots:f_n)$ is the 
  associated embedding $\iota_D\co X\hookrightarrow\P^n$.

  Let $\alpha\in\Aut_0(X)$. Then the pullback $\alpha^*\co\Div X\to\Div X$
  maps $D$ into $D - (g)$ for some $g\in k(X)$. Further, the image of an effective
  divisor under $\alpha$ is again effective. So for $f\in\mathcal{L}(D)$ we have 
  $\alpha^*(D + (f)) = D + (g^{-1} f\circ\alpha) \ge 0$.
  It follows that $f\mapsto g^{-1} f\circ\alpha$ is a projective transformation
  of $\mathcal{L}(D)$, i.e. there are $c_{ij}\in k$ such that
  $g^{-1}f_i\circ\alpha = \sum_{j=0}^n c_{ij}f_j$.
  Now by assumption $(f_0:\dots:f_n)$ is the embedding of $X\hookrightarrow\P^n$.
  Then $\alpha(X)$ is embedded by $(f_0\circ\alpha:\dots:f_n\circ\alpha) = 
  (g\sum_j c_{0j}f_j:\dots:g\sum_j c_{nj}f_j) = 
  (\sum_j c_{0j}f_j:\dots:\sum_j c_{nj}f_j)$.

  Since $f_0,\dots,f_n$ are linearly independent, $\iota_D(X)\subseteq\P^n$ 
  is not contained in any hyperplane. 
  Therefore a transformation of $\P^n$ fixing $X$ pointwise is the identity on $\P^n$. 
  Hence $\Aut_0(X)$ is injectively mapped into $\PGL_{n+1}(k)$.

  Finally, the set of elements $(c_{ij})_{i,j=0}^n\in\PGL_{n+1}(k)$
  in the image of $\Aut_0(X)$ is the zero set of polynomials
  obtained from the ideal generated by
  $f\circ\alpha = g\sum_{j=0}^n x_{ij}f_j$ for $i=0,\dots,n$,
  after eliminating all unknowns except $x_{ij}$'s.
  So $\Aut_0(X)$ is indeed an algebraic variety.
\end{proof}

We have seen that the group $\Aut_0(X)$ is algebraic,
therefore we can define the Lie algebra of an abstract variety.

\begin{df}\label{df:unembeddedLA}
  {\em The Lie algebra of a variety $X$} such that 
  $X$ has a very ample divisor 
  is the Lie algebra of the algebraic group $\Aut_0(X)$.
\end{df}

Note that the Lie algebra of an abstract variety $X$ and
the Lie algebra of its
embedding using a very ample divisor $D$ might differ.
Indeed, the image of the group $\Aut_0(X)$ after the embedding
$\iota_D$ is only a subgroup of (not equal to) $\P G(\iota(X),k)$,
so the Lie algebra of $X$ can happen to be only a subalgebra of $\g_0(X,k)$.
Nevertheless, in all examples considered in our work 
these two algebras are isomorphic.

\begin{exa} \label{exa:autPn}
  Let $X = \P^n$, $n>0$, so $\Pic(X)\cong\Z$. 
  Every automorphism of $\P^n$ fixes the Picard group.
  It follows that $\Aut_0(X) = \Aut(X) = \PGL_{n+1}(k)$,
  therefore for the Lie algebras we have
  $\g(X,k)=\gl_{n+1}(k)$, and $\g_0(X,k)=\sl_{n+1}(k)$.
\end{exa}

\begin{exa}\label{exa:autPxP}
  Let $X = \P^1\times\P^1$, so $\Pic X\cong \Z^2$.
  For automorphisms we have 
  $\Aut(\P^1\times\P^1) = G_1 \times G_2 \rtimes \left< g_0 \right>$,
  where $G_1, G_2 \cong \PGL_2(k)$ and $g_0$ switches the $\P^1$'s:
  $g_0: (s_0:s_1; t_0:t_1)\mapsto(t_0:t_1; s_0:s_1)$.

  The Segre embedding $X\hookrightarrow\P^3$, 
  $(s_0:s_1; t_0:t_1) \mapsto (s_0t_0 : s_0t_1 : s_1t_0 : s_1t_1)$
  is the map associated to a divisor of bidegree $(1,1)$.
  Here also the automorphism $g_0$ is represented by a matrix in $\PGL_4(k)$.
  We conclude that 
  $\Aut(X) = \left\{g\in\PGL_4(k) \mid gp\in X\ \forall p\in X\right\}$.
  The same holds for the anticanonical embedding of $\P^1\times\P^1$, 
  which is the map associated to a divisor of bidegree $(2,2)$.

  Finally, $\Aut_0(\P^1\times\P^1) = G_1 \times G_2$ is a subgroup of 
  $\Aut(\P^1\times\P^1)$ of finite index, therefore the Lie algebras 
  of an abstract variety (as in Definition~\ref{df:unembeddedLA})
  and the embedded ones (as in Definition~\ref{df:embeddedLA}) coincide.
\end{exa}

\begin{exa}\label{exa:autBlowup}
  Let $X$ be the blowup of $\P^2$ in $(1:0:0)$. Then $\Pic X \cong \Z^2$.
  Every automorphism of $X$ fixes the exceptional line.
  More precisely, $\Aut(X)$ is a subgroup of $\Aut(\P^2)$ fixing 
  the point $(1:0:0)$. Therefore $\Pic X$ is fixed by $\Aut(X)$, 
  hence $\Aut(X) = \Aut_0(X)$.
\end{exa}

\section{Computing $\g_0(X,k)$}\label{se:computeL}

The reason why we define the Lie algebra for an embedded variety is,
that such definition immediately delivers also tools for computing 
$G(X,k)$ and $\g_0(X,k)$. 

Let $\mathcal{I}(X) \subset k[x_0,\dots,k_n]$ be the vanishing ideal 
of $X\subset\P^n$ and let $\mathcal{B} = \{f_1,\dots,f_k\}$ be 
a finite set of homogeneous equations generating $\mathcal{I}(X)$.
For $g=(g_{ij})_{i,j=1}^{n+1}\in\GL_{n+1}(k)$ we have that 
\begin{eqnarray*}
  g\in G(X,k) & \Leftrightarrow & 
    f(gp) = 0 \textrm{ for all } p\in X \textrm{ and } f\in\mathcal{I} \\
  & \Leftrightarrow & f\circ g\in\mathcal{I}(X) \textrm{ for each } f\in\mathcal{B}.
\end{eqnarray*}
The last expression is valid if $f = \sum a_if_i$ 
with $a_i\in k[x_0,\dots,x_n]$ homogeneous and 
$f_i\in\mathcal{B}$ such that $\deg f_i \le \deg f$. 
For each $f\in\mathcal{B}$ therefore we form the following set of polynomials:
$\mathcal{S} = \{x^u f_i \mid f_i\in\mathcal{B}, \deg f_i \le \deg f, 
  |u| = \deg f - \deg f_i \}$,
where $u$ is the vector $(u_0,\dots,u_n)\in\N^{n+1}$, 
$|u| = u_0+\dots+u_n$ and
$x^u$ is a shortcut for $x_0^{u_0}\dots x_n^{u_n}$.
Then $f\circ g\in\mathcal{I}(X)$ if 
$f\circ g = \sum_{h_j\in\mathcal{S}} c_jh_j$ for some $c_j$'s in $k$.
In this way for each $f\in\mathcal{B}$ we get a polynomial equation
in $g_{ij}$'s and $c_j$'s defining together an (affine) algebraic variety.
Then $G(X,k)$ is its projection discarding all $c_j$'s. 
Let $J$ denote the ideal obtained in this way, i.e.~$G(X,k) = \mathcal{V}(J)$.
To find the Lie algebra of $X$ we have further to compute the radical $\sqrt{J}$
and differentiate its generators.

In the examples we are interested in, the ideal of the variety $X$ is 
generated by quadratic equations. 
In this case, we have much more straightforward method for computing 
the Lie algebra of $X$.

\begin{lem}\label{le:linearLA}
  Let $\pi\co G\to\GL(V)$ be a rational representation,
  and let $N\subset M$ be vector subspaces of $V$. Put
  \begin{displaymath}
    H = \{g\in G \mid \pi(g)N = N, \pi(g)M = M, \pi(g)_{M/N} = e \}.
  \end{displaymath}
  Then
  \begin{displaymath}
    L(H) = \{ x\in\mathfrak{g} \mid d\pi(x)M\subset N\}.
  \end{displaymath}
\end{lem}

\begin{proof}
  See~\cite{borel}, Lemma~7.4.
\end{proof}

\begin{thm}\label{thm:quadricsLA}
  Let $X\subset\P^n$ be embedded such that
  the ideal of $X$ is generated by quadrics $z^t A_i z$, 
  $i=1,\dots,r$, where $A_i$'s are symmetric matrices of degree $n+1$,
  and $z = (z_0\ z_1\ \dots\ z_n)^t$.
  Let $\mathcal{A} \subset\M_{n+1}(k)$ denote the subspace generated by all $A_i$'s.
  Then 
  \begin{displaymath}
    \g(X,k) = 
      \{x\in{\gl}_{n+1}(k) \mid x^tA_i+A_ix\in\mathcal{A}\ \ \forall i\}.
  \end{displaymath}
\end{thm}

\begin{proof}
  First we describe $G(X,k)$. 
  Let $\lambda_1,\dots,\lambda_s$ be linear forms defining 
  the subspace $\mathcal{A}\subset\M_{n+1}(k)$.
  By definition of $G(X,k)$ for $g\in\GL_{n+1}(k)$ we then have
  \begin{eqnarray*}
    g\in G(X,k) & \Leftrightarrow & z^t(g^tA_kg)z \in\mathcal{I}(X)\quad \forall k \\
    & \Leftrightarrow & g^tA_kg \in \mathcal{A}\quad \forall k \\
    & \Leftrightarrow & \lambda_l(g^tA_kg) = 0\quad \forall k,\ \forall l.
  \end{eqnarray*}
  So $G(X,k)$ is described by $kl$ quadratic forms in $k[g_{ij} \mid i,j=0,\dots,n]$.

  Except acting on on $\P^n$, $\GL_{n+1}(k)$ acts also on the space 
  $\M_{n+1}(k)$ by $g(A) = g^tAg$. 
  Then $G(X,k)$ is exactly the subgroup of $\GL_{n+1}(k)$ fixing $\mathcal{A}$.

  Let $m = (n+1)^2 = [\M_{n+1}(k):k]$.
  Let $V$ be $m$-dimensional vector space representing matrices in $\M_{n+1}(k)$.
  We will write $v_A$ for the vector in $V$ representing the matrix $A$.
  Let $\pi$ be the representation of $\GL_{n+1}(k)\to\GL_m(k)$, so
  $\pi(g)v_A = v_{g^tAg}$. 
  If we set $M=N=\mathcal{A}$, then by Lemma~\ref{le:linearLA} the Lie algebra 
  of $G(X,k)$ is
  \begin{displaymath}
    L(G) = \{x\in{\gl}_{n+1}(k) \mid d\pi(x)\mathcal{A}\subset\mathcal{A}\}.
  \end{displaymath}
  It remains to compute the differential of $\pi$.

  For $A = (A_{ij})_{i,j=0}^{n}\in\M_{n+1}(k)$ let 
  $v_A$ be the vector $(A_{00}, A_{01}, \dots, A_{0n}, A_{10},\dots, A_{nn})^t$.
  The coordinates in $v_A$ are naturally indexed by vectors:
  $(v_A)_{(i,j)} = A_{ij}$.
  Then for $g = (g_{ij})_{i,j=0}^n \in\GL_{n+1}(k)$ we have
  $(\pi(g))_{(i,j)(k,l)} = g_{ki}g_{lj}$.
  By differentiating in this coordinate we get 
  \begin{align*}
    (d\pi(x))_{u,v} 
     &= (x_{ki}\frac{\partial}{\partial g_{ki}} + 
         x_{lj}\frac{\partial}{\partial g_{lj}}) g_{ki}g_{lj} \mid_{g=e} 
     \ =\ (x_{ki} g_{lj} + x_{lj} g_{ki})\mid_{g=e} \\
     &= x_{ki}g_{lj}+x_{lj}g_{ki} \mid_{g=e} = \left\{
    \begin{array}{ll}
      x_{ii} + x_{jj} & \textrm{ if } i=k, j=l \\
      x_{lj}          & \textrm{ if } i=k, j\ne l \\
      x_{ki}          & \textrm{ if } i\ne k, j=l \\
      0               & \textrm{ if } i\ne k, j\ne l.
    \end{array} \right.
  \end{align*}
  The last tedious computation shows that 
  $(d\pi(x)v_A)_{u} = (v_{(xA^t + Ax)})_u$ for all $u=1,\dots,m$ 
  and the assertion of the theorem follows.
\end{proof}

\begin{exa}
  Let $X\subset\P^2$ be the zero set of $z_0^2 + z_1^2 + z_2^2$
  and let $k=\Q$.
  Though $X(\Q) = \emptyset$, there are still automorphisms of $X$ 
  defined over $\Q$. 
  We find the Lie algebra $\g(X,\Q)$.

  The group $G(X,\Q)$ consists of $g\in\GL_3(\Q)$ such that 
  $g^tI_3g = \lambda I_3$ for some $\lambda\in\Q$
  where $I_3$ is the identity matrix. 
  Therefore 
  \begin{displaymath}
    \g(X,\Q) = \{x\in{\gl}_3(\Q) \mid x^t + x =\lambda I_3 
    \textrm{ for some } \lambda\in\Q\},
  \end{displaymath}
  so it is generated by all antisymmetric matrices together with $I_3$.
\end{exa}

\section{Varieties and modules}

Let $X_1,X_2\subseteq\P^n$ be two varieties which are 
projectively equivalent over $k$, 
i.e.~there is $\mu\co k^{n+1}\to k^{n+1}$ (and consequently $\P^n\to\P^n$) 
given by a matrix $M\in{\GL}_{n+1}(k)$ such that 
\begin{equation}\label{eq:M}
  \mu\co p\mapsto Mp\ \textrm{ for all }\ p\in\P^n
  \ \textrm{ and }\ \mu(X_1) = X_2.
\end{equation}
Then $\mu$ induces a map $\tilde\mu\co G(X_1,k)\to G(X_2,k)$ 
such that the following diagram is commutative:
\begin{displaymath}
  \begin{diagram}
    X_1 & \rTo^{\mu} & X_2 \\
    \dTo^{g} &    &   \dTo_{\tilde\mu(g)} \\
    X_1 & \rTo^{\mu} & X_2 
  \end{diagram}
\end{displaymath}
i.e.~$\tilde\mu(g) = \mu\circ g\circ\mu^{-1}$.
It is straightforward to see that $\tilde\mu$ is an isomorphism of groups,
mapping $g$ to $MgM^{-1}$.

\begin{prop}\label{prop:lieAlgVarIso}
  Suppose that $X_1$ and $X_2\subseteq\P^n$ are projectively equivalent over $k$ 
  and let $M$ be the matrix as in~(\ref{eq:M}). Then
  \begin{enumerate}
    \item[(i)] 
      the map $d\tilde\mu\co x\mapsto MxM^{-1}$ is a Lie algebra
      isomorphism from $\g(X_1,k)$ to $\g(X_2,k)$,
    \item[(ii)] $\mu\co V_1\to V_2$, $v\mapsto Mp$ is an isomorphism
      of $(n+1)$-dimensional natural modules of $\g_0(X_1,k)$ and $\g_0(X_2,k)$.
  \end{enumerate}
\end{prop}

\begin{proof}
  The first assertion follows from Proposition~\ref{prop:lieAlgIso} and the fact 
  that $\tilde\mu\co g\mapsto MgM^{-1}$ is an isomorphism from 
  $G(X_1,k)$ to $G(X_2,k)$.

  The second assertion follows from Corollary~\ref{cor:moduleIso} and the fact
  that the restriction of $\varphi_L$ to $\g_0(X_1,k)$ leads to
  an isomorphism of $\g_0(X_1,k)$ and $\g_0(X_2,k)$.
\end{proof}

For our purposes, a converse of the last proposition is more useful. 
Namely, given two algebraic varieties in $\P^n$, if we can find
an isomorphism of their Lie algebras and corresponding modules, 
can we ``lift'' it to an
isomorphism of the varieties, provided it exists?
The following proposition gives a partial answer.

\begin{prop}\label{AlgVarEq}
  Suppose that $X_1$ and $X_2$ in $\P^n$ are projectively equivalent over $k$.
  Assume further that any automorphism of $\g_0(X_1,k)$ is a conjugation by 
  a matrix from $G(X_1,k)$. 
  Then any isomorphism $\psi\co \g_0(X_1,k)\to \g_0(X_2,k)$ leads to 
  an isomorphism $\nu$ of their modules.
  If $N\in\M_{n+1}(k)$ is the matrix such that $\nu\co v\mapsto Nv$
  then the matrix $N$ moreover defines
  a projective equivalence $X_1\to X_2$, $p\mapsto Np$.
\end{prop}

\begin{proof}
  By assumption there is a matrix $M\in\GL_{n+1}(k)$ describing a projective 
  equivalence of $X_1$ and $X_2$, as in~(\ref{eq:M}).
  By Proposition~\ref{prop:lieAlgVarIso} it leads to an isomorphism $d\tilde\mu$ 
  of Lie algebras $\g_0(X_1,k)\to \g_0(X_2,k)$ mapping $x$ to $MxM^{-1}$. 
  Then $d\tilde\mu^{-1}\circ\psi$ is an automorphism of $\g_0(X_1,k)$.
  By assumption there exists $G\in G(X_1,k)$ such that 
  \begin{displaymath}
    d\tilde\mu^{-1}\circ\psi(x) = GxG^{-1}. 
  \end{displaymath}
  It follows that
  \begin{displaymath}
    \psi(x) = d\tilde\mu(GxG^{-1}) = (MG)x(MG)^{-1}.
  \end{displaymath}
  Hence $\psi$ also leads to an isomorphism of modules, namely $v\mapsto MGv$,
  so $N = MG$.
  Finally, the matrix $N$ describes the composition of an automorphism of $X_1$
  and the projective equivalence of $X_1$ and $X_2$, hence the last assertion follows.
\end{proof}

The last proposition gives an algorithm for deciding existence and finding 
a linear projective isomorphism of two algebraic varieties provided 
their Lie algebras satisfy its assumptions:

\vskip0.4cm
\noindent
\begin{tabular}{lp{11.9cm}}
  {\sc Algorithm:} & {\tt FindProjectiveEquivalence}\\
  {\sc Input:}     & $X_1, X_2\subseteq\P^n$ -- projective varieties defined over $k$ 
                     such that each automorphism of $\g_0(X_1,k)$ is a conjugation by 
                     a matrix from $G(X_1,k)$,\\
  {\sc Output:}    & $M\in \GL_{n+1}(k)$ such that 
                     $Mp\in X_2$ for all $p\in X_1$ if such $M$ exists,\\
                   & \keyw{failed} otherwise.
\end{tabular}
\begin{enumerate}
  \item Compute $\g_0(X_1,k)$ and $\g_0(X_2,k)$ as described in Section~\ref{se:computeL}.
  \item Find a Lie algebra isomorphism $\psi\co \g_0(X_1,k)\to \g_0(X_2,k)$.
    If such isomorphism does not exist, return \keyw{failed}.
  \item Construct an isomorphism $\mu$ between the modules $V_1$ and $V_2$ 
    of algebras $\g_0(X_1,k)$ and $\g_0(X_2,k)$, 
    $\mu\co v\mapsto Mv$ for some $M\in\GL_{n+1}(k)$.
    If such isomorphism does not exist, return \keyw{failed}.
  \item If $\mu$ transforms $X_1$ to $X_2$, return $M$.
    Otherwise return \keyw{failed}.
\end{enumerate}
\vskip0.4cm

We will use this algorithm for finding a parametrization of a given
variety. More precisely, as the variety $X_1$ we take a variety with 
known parametrization. Then we find an isomorphism of $X_1$ and an
implicitly given variety $X_2$, obtaining so a parametrization of the latter.
Next we illustrate the parametrization algorithm in detail on 
a simple example of conics.

\section{Illustrating example: $\P^1$ and $\sl_2$}\label{se:P1}

Of course there is a lot of ways how to parametrize a given conic.
But for the sake of simplicity we use this easy example to 
illustrate our method.

We will find a parametrization of a conic $C\subset\P^2$ given 
by a quadratic form in $k[z_0,z_1,z_2]$. 
We will do so by finding an isomorphism of $C$ and $C_0: z_0z_2 - z_1^2$, 
where $C_0$ has a parametrization $(s:t)\mapsto(s^2:st:t^2)$. 

The automorphism group of $\P^1(k)$ is $\PGL_2(k)$, therefore we have
$G(\P^1,k)\cong\GL_2(k)$ and the Lie algebra 
$\g(\P^1,k) = L(\GL_2(k)) = \gl_2(k)$.
If we factor out the one dimensional subalgebra spanned by the identity matrix,
we get 
\begin{displaymath}
  \g_0(\P^1,k)\cong{\sl}_2(k).
\end{displaymath}
If $C$ has a $k$-rational parametrization (or equivalently $C$ is over $k$
projectively equivalent to $C_0$), then also the Lie algebra $\g_0(C,k)$ 
is isomorphic to $\sl_2(k)$. 

Next we need to check the assumption of Proposition~\ref{AlgVarEq}.
Firstly, Theorem~\ref{thm:slAut} asserts that any automorphism of $\mf{sl}_2(k)$
is a conjugation by a regular $2\times 2$ matrix.

We will now examine the conic $C_0$.
Let $v_0,v_1$ be the standard basis of $V=k^2$. Let
$W=\Sym^2(V)$ with the basis
$v_0^2$, $2v_0v_1$, $v_1^2$.
Let $\varphi'_0\co V\to W$ be defined by
$\varphi'_0(v) = v^2$. We write the coordinates of an element in
$W$ with respect to the basis above. Then the image of the induced
map $\varphi_0\co\P(V)\to\P(W)$ in these coordinates is exactly $C_0$ 
(see also~\cite{fh}, \S 11.3).
The group $G(\P^1,k) = \GL_2(k)$ acts on $W$ by
$g\cdot(v_1v_2) = (gv_1)(gv_2)$.
It follows that an automorphism of $\P^1$ 
represented by the matrix $g = (g_{ij})_{i,j=0,1}$,
\begin{displaymath}
  g p = \left( \begin{array}{@{\ }cc@{\ }}
    g_{00} & g_{01} \\
    g_{10} & g_{11}
  \end{array} \right)
  \left( \begin{array}{@{\ }c@{\ }}
    p_0 \\ p_1
  \end{array} \right),
\end{displaymath}
is by the parametrization $\varphi_0\co (s:t)\mapsto(s^2:st:t^2)$
transformed to a linear automorphism of $\P^2$ fixing $C_0$:
\begin{equation*}
  g\cdot\varphi_0(p) = \tilde\varphi_0(g) \varphi_0(p) = \left( 
  \begin{array}{@{\ }ccc@{\ }}
    g_{00}^2 & 2g_{00}g_{01} & g_{01}^2 \\
    g_{00}g_{10} & g_{00}g_{11} + g_{01}g_{10} & g_{01}g_{11} \\
    g_{10}^2 & 2g_{10}g_{11} & g_{11}^2 
  \end{array} \right)
  \left( \begin{array}{@{\ }c@{\ }}
    p_0^2 \\ p_0p_1 \\ p_2^2
  \end{array} \right),
\end{equation*}
so $\tilde\varphi_0$ is a group homomorphism $G(\P^1, k) = \GL_2(k)\to\GL_3(k)$.
The differential of $\tilde\varphi_0$ maps the Lie algebra of $G(\P^1,k)$
isomorphically to the Lie algebra of $G(C_0,k)$. 
It describes the action of the Lie algebra of $G(\P^1,k)$ on $W$, 
$x.(v_1v_2) = d\tilde\varphi_0(x)(v_1v_2) = (xv_1)v_2 + v_1(xv_2)$.
Explicitly it is given by
\begin{equation}\label{eq:d_varphi_sl2}
  d\tilde\varphi_0:\left( \begin{array}{@{\ }cc@{\ }}
    a_{00} & a_{01} \\
    a_{10} & a_{11}
  \end{array} \right)\mapsto
  \left( \begin{array}{@{\ }ccc@{\ }}
    2a_{00} & 2a_{01}         & 0 \\
    a_{10}  & a_{00} + a_{11} & a_{01} \\
    0       & 2a_{10}         & 2a_{11} 
  \end{array} \right)
\end{equation}
and maps isomorphically $\sl_2(k)\to \g_0(C_0,k)$.

From Theorem~\ref{thm:slAut} and Proposition~\ref{prop:LieConjug} now follows
that the algorithm {\tt Find\-Pro\-jective\-Equivalence}
can be used to parametrize a given conic $C$.

In the second step of the algorithm we have to find an isomorphism $\psi$
of Lie algebras of $C_0$ and $C$. We have already constructed an
isomorphism $d\tilde\varphi_0\co\sl_2\to \g_0(C_0,k)$, (see~(\ref{eq:d_varphi_sl2})).
Therefore it is enough to find an isomorphism $d\tilde\varphi\co\sl_2\to \g_0(C,k)$
\footnote{The prefix ``$d$'' in the name of the map $d\tilde\varphi$
suggests that $d\tilde\varphi$
is a differential of a map between algebraic groups. 
This indeed is true,
nevertheless we don't prove and don't use this fact. 
We look for the map $d\tilde\varphi$ 
independently on that.}
and compose $\psi = d\tilde\varphi\circ(d\tilde\varphi_0)^{-1}$.
Finding an isomorphism between Lie algebras 
is usually the most difficult part of the algorithm,
and is discussed in Chapter~\ref{ch:csa}.
Nevertheless, assume for now that we can find $d\tilde\varphi$. 
Then the last construction to be made
is finding an isomorphism between $3$-dimensional 
natural modules of $\g_0(C_0,k)$ and $\g_0(C,k)$.

\begin{lem}
  The $3$-dimensional $\sl_2$-module afforded by the representation
  $d\tilde\varphi_0\co\sl_2\to \g_0(C_0,k)$ is irreducible with
  highest weight $2$.
\end{lem}

\begin{proof}
  We have already seen that the $3$-dimensional $\g_0(C_0,k)$-module 
  is isomorphic to $\Sym^2(V)$, and therefore it is irreducible.
  The highest weight of this module is $2$.
\end{proof}

Let us fix a Chevalley basis of $\g_0(C_0,k)$ to be the image of the standard
Chevalley basis in $\sl_2(k)$ under $d\tilde\varphi_0$:
\begin{equation}\label{eq:sl2Chev}
  x_0 
  = \left( \begin{array}{@{\ }ccc@{\ }}
    0 & 2 & 0 \\
    0 & 0 & 1 \\
    0 & 0 & 0 
  \end{array} \right),\quad
  y_0 = \left( \begin{array}{@{\ }cccr@{\ }}
    0 & 0 & 0 \\
    1 & 0 & 0 \\
    0 & 2 & 0
  \end{array} \right),\quad
  h_0 = \left( \begin{array}{@{\ }cc@{\ }r@{\ }}
    2 & 0 & 0 \\
    0 & 0 & 0 \\
    0 & 0 & -2
  \end{array} \right).
\end{equation}
Then the highest weight vector in the $\g_0(C_0,k)$-module is
\begin{equation}\label{eq:C0eigenspaces}
  v_0^+ = e_1,
\end{equation}
where by $e_i$ we denote the 3-dimensional column vector with $1$ at the $i$-th
place and $0$ elsewhere. 
For the rest of the module we get 
\begin{equation}\label{eq:C0eigenspaces2}
  y_0\cdot v_0^+ = e_2 \quad\textrm{ and }\quad y_0\cdot(y_0\cdot v_0^+) = 2e_3.
\end{equation}

\begin{rem}
  Surely, the computation just made also shows that the natural 
  $\g_0(C_0,k)$-module is irreducible.
  Nevertheless the conceptual argumentation about the structure
  of the module is preferable
  since it easily generalizes to varieties of higher dimension,
  where the explicit computation would be too lengthy.
\end{rem}

The conic $C$ is given implicitly, therefore we have to find its Lie algebra
$\g_0(C,k)\subset\sl_3(\Q)$ by solving the linear system according 
to Theorem~\ref{thm:quadricsLA}, and intersecting with $\sl_3(\Q)$.
 
\begin{prop}\label{prop:conic-sl2}
  The conic $C\subset\P^2$ defined over $k$ has a $k$-rational parametrization
  if and only if $\g_0(C,k)\cong\sl_2(k)$.
\end{prop}

\begin{proof}
  If $C$ has a parametrization over $k$, then it is projectively equivalent 
  to $C_0$ ($C_0$ as above) and hence by Proposition~\ref{prop:lieAlgVarIso}, 
  $\g_0(C,k)\cong\sl_2(k)$.

  Let on the other hand $\g_0(C,k)\cong\sl_2(k)$. Let $k'\supseteq k$ be
  a field extension such that $C$ has a $k'$-rational parametrization.
  Then the natural modules of $\g_0(C,k')$ and $\g_0(C_0,k')$ are 
  two isomorphic $\mf{sl}_2(k')$-modules.
  Hence the natural $\g_0(C,k')$-module is irreducible with the highest weight $2$.
  But $\g_0(C,k)$ has a split Cartan subalgebra, therefore the corresponding 
  highest weight vector in the natural $\g_0(C,k')$-module 
  is defined over $k$ and so is the
  isomorphism of modules and of $C_0$ and $C$.
\end{proof}

We can assume that the Lie algebra $\g_0(C,k)\cong\sl_2(k)$
since otherwise $C$ does not have a $k$-rational parametrization.
From the proposition it follows that if there is an isomorphism $\sl_2\to \g_0(C,k)$,
then also the natural $\g_0(C_0,k)$- and $\g_0(C,k)$-modules are $k$-isomorphic.
Suppose we can find a Chevalley basis $x,y,h$ in $\g_0(C,k)$, 
which is equivalent to finding an isomorphism with $\sl_2(k)$. 
We find then a module isomorphism by finding the highest weight vector $v^+$ 
in the natural $\g_0(C,k)$-module and mapping
$v_0^+ \mapsto v^+$, 
$y_0\cdot v_0^+ \mapsto y\cdot v^+ $ and
$y_0\cdot(y_0\cdot v_0^+) \mapsto y\cdot(y\cdot v^+)$.

\begin{exa}
  Let us parametrize the unit circle $C\subset\P^2$
  given by the polynomial $z_1^2 + z_2^2 - z_0^2$ or, equivalently by
  the diagonal matrix $A\in M_3(\Q)$ with entries $(-1,1,1)$.

  We find $\g_0(C,\Q)\subset\sl_3(\Q)$ by solving the linear system
  $x^t A + Ax\in \Span_\Q\{A\}$ for $x\in\sl_3(\Q)$. 
  A Chevalley basis of $\g_0(C,\Q)$ is for example
  \begin{displaymath}
    x = \left( \begin{array}{@{\ }r@{\ }rr@{\ }}
      0 & 0 & 1 \\
      0 & 0 & 1 \\
      1 & -1 & 0 
    \end{array} \right),\quad
    y = \left( \begin{array}{@{\ }cc@{\ }r@{\ }}
      0 & 0 & 1 \\
      0 & 0 & -1 \\
      1 & 1 & 0
    \end{array} \right),\quad
    h = \left( \begin{array}{@{\ }ccc@{\ }}
      0 & 2 & 0 \\
      2 & 0 & 0 \\
      0 & 0 & 0
    \end{array} \right).
  \end{displaymath}

  The maximal weight vector and the derived vectors in the natural 
  $\g_0(C,k)$-module are
  \begin{displaymath}
    v^+ = e_1 + e_2,\quad y\cdot v^+ = 2e_3 \quad\textrm{ and }\quad 
    y\cdot(y\cdot v^+) = 2(e_1 - e_2).
  \end{displaymath} 
  Comparing with~(\ref{eq:C0eigenspaces}) and (\ref{eq:C0eigenspaces2})
  yields that the matrix of an isomorphism of modules is
  \begin{displaymath}
    M = \left(\begin{array}{@{\ }ccc@{\ }}
      1 & 0 & 1  \\
      1 & 0 & -1 \\
      0 & 2 & 0
    \end{array}\right),
  \end{displaymath}
  which gives a parametrization of the circle
  \begin{displaymath}
    \left(\begin{array}{@{\ }c@{\ }}
      z_0 \\ z_1 \\ z_2
    \end{array}\right)
    =
    \left(\begin{array}{@{\ }ccc@{\ }}
      1 & 0 & 1  \\
      1 & 0 & -1 \\
      0 & 2 & 0
    \end{array}\right)
    \left(\begin{array}{@{\ }c@{\ }}
      s^2 \\ st \\ t^2
    \end{array}\right)
    =
    \left(\begin{array}{@{\ }c@{\ }}
      s^2 + t^2 \\
      s^2 - t^2 \\
      2st 
    \end{array}\right).
  \end{displaymath}
\end{exa}

\section{Conics and their Lie algebras}\label{se:conics}

Here we explain how to use the known algorithms for parametrizing
conics when we want to find an isomorphism of a Lie algebra and $\mf{sl}_2$.

In the previous example the unexplained step in the algorithm was finding
an isomorphism $\mf{sl}_2(\Q)\to \g_0(C,\Q)$. 
This is as difficult as finding a rational point on a conic. 
There are already efficient algorithms for solving the latter 
(cf.~\cite{cremona, simon2}), some even implemented in Magma.
Therefore in the case of conics
we usually do a reduction in the other direction:
If for a given Lie algebra $\g$ over $\Q$ we have to find
an isomorphism $\mf{sl}_2(\Q)\to \g$, then we first find a conic
$C$ such that $\g$ is isomorphic to its Lie algebra, 
and afterwards we use the known fast
algorithms for finding a rational point on $C$ and consequently
a parametrization of $C$, which will lead to constructing an
isomorphism of $\g$ and $\mf{sl}_2(\Q)$.
We have to show that for each $\g$ which is 
a twist of $\mf{sl}_2$ it is possible to construct 
a corresponding conic.

Here again $C_0\subset\P^2$ is the conic given by the form $z_0z_2 - z_1^2$ 
and parametrized by monomials $(s : t) \mapsto (s^2 : st : t^2)$.

\begin{lem}\label{le:basis}
  Let $\g$ be a twist of $\mf{sl}_2(\Q)$.
  There is a basis $h,x,y$ of $\g$ such that
  the multiplication in $\g$ is given by 
  $[h,x] = y$, $[h,y]=ax$, $[x,y]=bh$ for some $a,b\in\Q^*$.
\end{lem}

\begin{proof}
  Since $\g$ is a twist of $\mf{sl}_2(\Q)$, it has a 1-dimensional
  Cartan subalgebra. Let $h$ spans a Cartan subalgebra of $\g$.
  Then we do a Fitting decomposition, $\g$ = $\Span_\Q\{h\} + [h,\g]$.
  The subalgebra $[h,\g]$ is generated by two
  vectors $x,y$, where $x$ is chosen not to be an eigenvector of $\ad h$, 
  and $y = [h,x]$.
  Then we also have $[x,y]=bh$ and $[h,y]=ax$ with both $a,b\in\Q$ nonzero. 
  The first multiplication rule can be seen after passing to a splitting field of $\g$, 
  the second follows from Jacobi identity.
\end{proof}

If the Cartan subalgebra, which we find while looking for a basis as 
in the last lemma, happens to be split, we are done, 
because we can easily find an isomorphism with $\mf{sl}_2(\Q)$,
for example by using the method of \cite{gra6}, \S 5.11.
Otherwise we continue with the computation.

The adjoint representation of $\g$ with respect to this basis is 
a representation by the following 3 by 3 matrices over the rationals:
\begin{displaymath}
  \ad x = \left( \begin{array}{@{\ }cc@{\ }r@{\ }}
    0 & 0 & 0 \\
    0 & 0 & -1 \\
    0 & b & 0 
  \end{array} \right),\quad
  \ad y = \left( \begin{array}{@{}rc@{\ }r@{\ }}
    0 & 0 & -a \\
    0 & 0 & 0 \\
   -b & 0 & 0
  \end{array} \right),\quad
  \ad h = \left( \begin{array}{@{\ }ccc@{\ }}
    0 & a & 0 \\
    1 & 0 & 0 \\
    0 & 0 & 0
  \end{array} \right).
\end{displaymath}
It is a straightforward computation to check that this is
the Lie algebra of the conic $C$ given by
\begin{equation}\label{eq:conic}
  -bz_0^2 + abz_1^2 + az_2^2 = 0,\quad a,b\in\Q^* 
\end{equation}
We use known algorithms to find a rational point on $C$. If there is
none, by Proposition~\ref{prop:conic-sl2} the given Lie algebra $\g$
was a proper twist of $\mf{sl}_2$. Otherwise we can construct
a Lie algebra isomorphism.

\begin{lem}\label{le:conicPar}
  Let $(p_0 : p_1 : p_2)$ be a point on the conic
  \begin{equation}\label{eq:conicGen}
    Az_0^2 + Bz_1^2 + Cz_2^2 = 0,
  \end{equation}
  such that $p_2$ is nonzero. 
  Then we have the following parametrization:
  \begin{displaymath}
    \left( \begin{array}{@{\ }c@{\ }} z_0 \\ z_1 \\ z_2 \end{array} \right) = 
    \left( \begin{array}{@{}rcr@{\ }} 
       A p_0  &  2B p_1  &  -B p_0 \\ 
      -A p_1  &  2A p_0  &   B p_1 \\ 
       A p_2  &  0       &   B p_2
    \end{array} \right)
    \left( \begin{array}{@{\ }c@{\ }} s^2 \\ st \\ t^2 \end{array} \right)    
  \end{displaymath}
\end{lem}

\begin{proof}
  Direct verification proves the lemma.
  When looking for the matrix, we followed the method 
  described in~\cite{mordell}: 
  One looks for the parametrization in the shape 
  $(z_0 : z_1 : z_2) = (s + p_0 w : t + p_1 w : p_2 w)$
  and eliminates $w$.
\end{proof}

By plugging in $-b, ab, a$ for $A,B,C$ in the lemma we have constructed 
an isomorphism of $C$ and the conic $C_0$ given by $z_0z_2 - z_1^2$.
Let us denote the matrix describing the isomorphism by $M$,
i.e.~for every $p\in C_0$ we have $Mp\in C$.
Let $x_0, y_0, h_0$ be as in (\ref{eq:sl2Chev}).
Then by Proposition~\ref{prop:lieAlgVarIso}(i), 
$Mx_0M^{-1}$, $My_0M^{-1}$ and $Mh_0M^{-1}$ is a Chevalley basis of $\g$,
so we have constructed an isomorphism $\mf{sl}_2(\Q)\to \g$.

\begin{rem}
  Taking the adjoint representation with respect to the basis 
  as in Lemma~\ref{le:basis} leads to constructing a conic 
  in the diagonal shape.
  As pointed out to us by Michael Stoll, this usually introduces
  additional difficulties when looking for a rational point on the conic
  (one has to factorize too large integers).
  Therefore in our algorithms we do not look for the special
  basis of the Lie algebra but immediately construct the adjoint 
  representation, and then the corresponding conic by solving the 
  linear system in the Theorem~\ref{thm:quadricsLA} for an unknown 
  matrix $A$. The results in this section assure that there is
  always a conic corresponding to the adjoint representation 
  of the Lie algebra.
\end{rem}

\section{Variety as an orbit of an algebraic group}

The core theorem in the example of the conic (see Section~\ref{se:P1}) 
is the first part of 
Theorem~\ref{thm:slAut} saying that each automorphism of $\mf{sl}_2$
(represented as $2\times 2$ matrices of trace 0) is a conjugation by 
a regular matrix. Unfortunately this is not true for all algebras 
appearing in examples we would consider. Therefore
in this section we develop another machinery to prove an assertion similar to
one in Proposition~\ref{AlgVarEq}.

Let $G$ be an algebraic group and
let the variety $X$ be a {\em $G$-set}, 
i.e.~there is an action $G\times X\to X$,
$(g,x)\mapsto g\cdot x$, such that
\begin{displaymath}
  e\cdot x = x\quad\textrm{and}\quad g\cdot(h\cdot x) = (gh)\cdot x
\end{displaymath}
for all $x\in X$ and all $g,h\in G$.
Let $x\mapsto g\cdot x$ be an automorphism of $X$ for each $g\in G$.
Then we refer to this situation by saying that 
{\em $G$ acts morphically on the variety $X$} or also that
$X$ is a {\em $G$-variety}.
If both $G$ and $X$ are given over $k$, we say also that 
$G$ acts {\em $k$-morphically} on $X$.

As $X$ we will usually take a projective
variety $X\subseteq\P^n$ and as the group $G$ a linear subgroup
of $\GL_{n+1}(k)$ representing linear transformations of $\P^n$ fixing $X$,
hence the action of $G$ on $X$ is just a multiplication $(g,x)\mapsto gx$.
Of course, $G$ acts on the whole $(n+1)$-dimensional vector space
(which we refer to as a {\em $G$-module}) and $X$ is a subset fixed by $G$.

\begin{prop}[Closed orbit lemma]
  Let $G$ be an algebraic group acting morphically on a non-empty variety $X$. 
  Then each orbit is a smooth variety which is open in its closure in $X$.
  Its boundary is a union of orbits of strictly lower dimension.
  In particular, the orbits of minimal dimension are closed.
\end{prop}

\begin{proof}
  See~\cite{borel}, Lemma~1.8.
\end{proof}

The {\em derived series} $\mathcal{D}^n G$ of a group $G$ is defined
inductively by:
\begin{displaymath}
  \mathcal{D}^0G = G,\quad
  \mathcal{D}^{n+1}G = (\mathcal{D}^nG, \mathcal{D}^nG),\quad
  n\ge 0,
\end{displaymath}
where by $(.,.)$ we denote the commutator subgroup.

A group $G$ is {\em solvable}, if $\mathcal{D}^nG = \{e\}$ 
for some $n$.

A connected solvable group $G$ is {\em $k$-split}, 
if it has a composition series 
$G = G_0\supset G_1\supset\dots\supset G_s=\{e\}$
consisting of connected $k$-subgroups such that $G_i/G_{i+1}$
is $k$-isomorphic to $\mathbf{G}_a$ (the additive group of the field)
or $\GL_1$ ($0\le i<s$).

\begin{thm}[Borel's fixed point theorem]\label{thm:fixedPt}
  Let $G$ be connected, solvable and $k$-split group acting $k$-morphically
  on a projective variety $X$. If $X(k)\ne\emptyset$, 
  then there is an $x\in X$ with $gx = x$ for all $g\in G$.
\end{thm}

\begin{proof}
  See~\cite{borel}, Proposition~15.2.
\end{proof}

Let $G$ be an algebraic group. A maximal connected solvable subgroup
of $G$ is called a {\em Borel subgroup}.

\begin{cor}\label{cor:1orbit}
  Let $G\subseteq\GL(V)$ be an algebraic group and 
  $B\subseteq G$ a $k$-split Borel subgroup.
  If $B$ fixes a unique point in $\P(V)$, 
  then $G$ has a unique closed orbit in $\P(V)$
  (i.e. there is a unique projective variety $X\subseteq\P(V)$
  on which $G$ acts transitively).
\end{cor}

\begin{proof}
  Let $X^\prime = Gx$, where $x\in\P(V)$ is the unique point fixed by $B$.
  Let $X\subseteq\P(V)$ be a closed orbit of $G$. 
  By the Closed orbit lemma $X$ exists and is an algebraic variety.
  By Borel's fixed point theorem, $x\in X$, so $X^\prime\subseteq X$.
  Since $X$ itself is an orbit under $G$, we have that $X^\prime = X$.
\end{proof}

\begin{exa}
  Let us again take the conic $C_0$ which is the zero set of $z_0z_2 - z_1^2$,
  see Section~\ref{se:P1}.
  It is the image of $\P^1$ under the map $\varphi_0\co\P^1\to\P^2$, 
  $(s:t)\mapsto(s^2:st:t^2)$.
  Let $V=k^2$ and let $\varphi_0'$ be the map $V\to\Sym^2(V)$, $v\mapsto v^2$.
  Let $v_1,v_2$ be a basis of $V$ and let us take $v_1^2, 2v_1v_2$ and $v_2^2$
  as the basis of $\Sym^2(V)$. Then for the coordinates of the image
  $\varphi_0'(V)$ with respect to the fixed basis of $\Sym^2(V)$ we have
  $\P(\varphi_0'(V)) = C_0$.
  The group $\GL_2(k)$ acting naturally on $V$ acts on $\Sym^2(V)$ by
  $g\cdot(v_i v_j) = (g v_i)(g v_j)$ so that $\varphi_0'$ is $\GL_2(k)$-equivariant.
  Therefore $C_0$ is fixed by $\GL_2(k)$ and hence we have a representation
  $\tilde\varphi_0\co\GL_2(k)\to G(C_0,k)\subset\GL_3(k)$.

  The group $\GL_2(k)$ has a Borel subgroup
  \begin{displaymath}
    B = \left\{\left(\begin{array}{cc}
      \ast & \ast \\ 0 & \ast \end{array}\right)\in{\GL}_2(k)\right\},
  \end{displaymath}
  the group of all upper triangular matrices.
  Since $\GL_2(k)$ acts transitively on $C_0$, the conic is a single orbit of $G(C_0,k)$.
  On the other hand, via the action $\tilde\varphi_0$, $B$ has a unique fixed point 
  in $\P(\Sym^2(V))$, namely $(1:0:0)$.
  Therefore $C_0$ is the unique closed orbit of $G(C_0,k)$ in $\P(\Sym^2(V))$.
\end{exa}

\begin{prop}\label{prop:ProjEqBorel}
  Let $X_1$ and $X_2$ in $\P^n$ be projectively equivalent over $k$.
  Suppose that $G(X_i,k)$ $(i=1,2)$ acts transitively on $X_i$ and
  that a $k$-split Borel subgroup of $G(X_1,k)$ fixes a unique $x\in\P^n$.
  Suppose further 
  that the natural $\g_0(X_i,k)$-module $V_i$ $(i=1,2)$ is absolutely irreducible.
  If there is an isomorphism $\varphi$ of Lie algebras $\g_0(X_1,k)\to \g_0(X_2,k)$
  and an invertible map $\nu\co V_1\to V_2, v\mapsto Nv$ such that
  $\nu(xv) = \varphi(x)\nu(v)$ for all $x\in \g_0(X_1,k)$ and $v\in V_1$,
  then $p\mapsto Np$ is also a projective equivalence $X_1\to X_2$.
\end{prop}

\begin{proof}
  We claim that there is an isomorphism $\psi$ of algebraic groups
  $G(X_1,k)\to G(X_2,k)$ such that $\varphi = d\psi$ and the natural
  modules of $G(X_1,k)$ and $G(X_2,k)$ are isomorphic, i.e.~there is
  $\nu'\co V_1\to V_2$ such that $\nu'(gv) = \psi(g)\nu'(v)$ 
  for all $g\in G(X_1,k)$ and $v\in V_1$.
  Then by Corollary~\ref{cor:moduleIso}, 
  $\nu'$ is also an isomorphism of natural modules of $\g_0(X_1,k)$ and $\g_0(X_2,k)$
  with respect to $\varphi$, so~$\nu'(xv) = \varphi(x)\nu'(v)$.
  Since the modules are absolutely irreducible,
  by the Schur's lemma (see e.g.~\cite{humLA}, \S~6.1) 
  we have $\nu'=\nu$ up to multiplication by scalars.

  Now by Corollary~\ref{cor:1orbit}, there is a unique closed orbit in the natural 
  $G(X_1,k)$-module which is by $\nu$ mapped to a unique closed orbit 
  in the natural $G(X_2,k)$-module.
  Since $G(X_i,k)$ fixes $X_i$ ($i=1,2$) and $X_i$ is a closed orbit of $G$ 
  with respect to the action described by $G(X_i,k)$, we have that
  $\nu$ is also a projective equivalence $X_1\to X_2$.

  It remains to prove the claim in the beginning of the proof.
  By assumption there is a projective equivalence $\mu\co v\mapsto Mv$, 
  $\mu(X_2) = X_1$. Then $\psi_\mu\co g\mapsto MgM^{-1}$ is an isomorphism 
  $G(X_2,k)\to G(X_1,k)$ and $\mu$ is the associated module isomorphism.
  Then also the modules of their Lie algebras are isomorphic, 
  $\mu(xv) = d\psi_\mu(x)\mu(v)$ for all $x\in \g_0(X_2,k)$ 
  and $v$ in the natural $\g_0(X_2,k)$-module.
  (see Corollary~\ref{cor:moduleIso} or Proposition~\ref{prop:lieAlgVarIso}). 
  For the composition $\mu\circ\nu$ we then get
  $\mu\circ\nu(xv) = (d\psi_\mu\circ\varphi(x))(\mu\circ\nu(v))$.
  Let $A$ be a regular matrix such that $\mu\circ\nu(v) = Av$,
  So $Axv = (d\psi_\mu\circ\varphi(x))Av$. 
  Since the last equality holds for all $v$ in the module,
  we can conclude that the automorphism $d\psi_\mu\circ\varphi$ of $\g_0(X_1,k)$
  maps $x$ to $AxA^{-1}$.
  Then $\psi'\co g\mapsto AgA^{-1}$ is an injective homomorphism of the group $G(X_1,k)$.
  The image $\psi'(G(X_1,k))$ is an algebraic group with the same Lie algebra as the one 
  of $G(X_1,k)$, therefore $\psi'(G(X_1,k)) = G(X_1,k)$.
  If we define $\psi\co G(X_1,k)\to G(X_2,k)$ as the composition
  $\psi_\mu^{-1}\circ\psi'$, $g\mapsto M^{-1}A g (M^{-1}A)^{-1}$,
  then $d\psi = \varphi$.
\end{proof}

Proposition~\ref{prop:ProjEqBorel} leads to an algorithm
for finding a projective equivalence of two varieties
under more general assumptions:

\vskip0.4cm
\noindent
\begin{tabular}{lp{11.9cm}}
  {\sc Algorithm:} & {\tt FindProjectiveEquivalence}\\
  {\sc Input:}     & $X_1, X_2\subseteq\P^n$ -- projective varieties defined over $k$ 
                     such that at least one of the following conditions is satisfied:\\
                   & \hangafter 1 \hangindent 19pt $\quad\bullet$ 
                     each automorphism of $\g_0(X_1,k)$ is a conjugation 
                     by a matrix from $G(X_1,k)$,\\
                   & \hangafter 1 \hangindent 19pt $\quad\bullet$ 
                     $G(X_1,k)$ acts transitively on $X_1$,
                     a Borel subgroup of $G(X_1,k)$ fixes a unique point in $\P^n$,
                     and the natural $\g_0(X_1,k)$-module is absolutely irreducible.\\
  {\sc Output:}    & $M\in \GL_{n+1}(k)$ such that 
                     $Mp\in X_2$ for all $p\in X_1$ if such $M$ exists,\\
                   & \keyw{failed} otherwise.
\end{tabular}

\begin{enumerate}
  \item Compute $\g_0(X_1,k)$ and $\g_0(X_2,k)$ 
    as described in Section~\ref{se:computeL}.
  \item Find a Lie algebra isomorphism $\psi\co \g_0(X_1,k)\to \g_0(X_2,k)$
    such that the natural modules of $\g_0(X_1,k)$ and $\g_0(X_2,k)$ are isomorphic.
    If such isomorphism does not exist, return \keyw{failed}.
  \item Construct an isomorphism $\mu$ between the natural modules 
    of algebras $\g_0(X_1,k)$ and $\g_0(X_2,k)$, 
    $\mu\co v\mapsto Mv$ for some $M\in\GL_{n+1}(k)$.
  \item If $\mu$ transforms $X_1$ to $X_2$, return $M$.
    Otherwise return \keyw{failed}.
\end{enumerate}

\chapter{Del Pezzo surfaces of degree 8}

We illustrate how to use Lie algebras for parametrizing 
Del Pezzo surfaces of degree~8. 
As already mentioned, there are two kinds of these surfaces:
$\PxP$ and the blowup of the projective plane in a point.
Here we give a parametrization algorithm over $k$ for both of them.

\section{Embeddings of $\PxP$}

The automorphism group of $\PxP$ consists of automorphisms of particular
projective lines in the product and automorphisms switching the two
factors, the former fixing the Picard group of $\PxP$, the latter not.
So we have $\Aut(\PxP) = \PGL_2(k)\times \PGL_2(k)\rtimes \Z/2\Z$ and
$\Aut_0(\PxP) = \PGL_2(k)\times \PGL_2(k)$. 
It follows that $\g_0(\PxP,k) = \slxsl$.

The automorphisms of $\PxP$ can be represented by a linear group
$G\subset\GL_4(k)$ in the following way.
Let $G = G_0\rtimes \left<g\right>$, 
where $G_0\cong\GL_2\times\GL_2(k)$ consists of block diagonal 
matrices $\diag(g_1,g_2)$ with $g_1,g_2\in\GL_2(k)$,
and 
\begin{displaymath}
  g = \left( \begin{array}{@{\ }cc@{\ }}
    0   & I_2 \\
    I_2 & 0   
    \end{array} \right),\quad\textrm{where }
  I_2 = \left(\begin{array}{cc}
    1 & 0 \\
    0 & 1
    \end{array}\right).
\end{displaymath}
The subgroup $G_0$ acts on $\PxP$ by $\diag(g_1,g_2)\cdot(u,v) = (g_1 u, g_2 v)$
and $g$ acts by $g\cdot(u,v) = (v,u)$.
This readily gives a surjective group homomorphism $G\to\Aut(\PxP)$.
The kernel of this homomorphism is a $2$-dimensional subgroup
consisting of matrices $\diag(c_1I_2, c_2I_2)$ with $c_1,c_2\in k^*$.

Let $\g_0\subseteq\mf{sl}_4(k)$ be the Lie algebra of $G$ with
the ideal generated by the identity matrix factored out.
It consists of block diagonal matrices, each block being 
a $2\times 2$ zero trace matrix.
The algebra $\g_0(\PxP,k)$ is then the image of $\g_0$
under the differential of the group homomorphism $G\to\Aut(\PxP)$.
The automorphism group of $\g_0$ is generated by the automorphisms 
of each of the direct summands of $\g_0$, 
along with the map switching them. 
By Theorem~\ref{thm:slAut} the former are conjugations
by a matrix from $G_0$, and the latter is the conjugation by $g$.

\subsection{The basic algorithm}

The anticanonical embedding $\varphi_0\co\PxP\to X_0\subset\P^8$ is given by
\begin{equation}\label{eq:P1xP1antican}
  (s_0{:}s_1;\ t_0{:}t_1)\mapsto
    (s_0^2t_0^2:s_0^2t_0t_1:s_0^2t_1^2:
    s_0s_1t_0^2:s_0s_1t_0t_1:s_0s_1t_1^2:
    s_1^2t_0^2:s_1^2t_0t_1:s_1^2t_1^2).
\end{equation}
We want to decide whether a given surface in $\P^8$ projectively
equivalent to $X_0$ over $\overline k$, is isomorphic to $X_0$
also over $k$ and if so, to find an isomorphism.

We now analyze the natural $\g_0(X_0,k)$-module.
Let $V_1$ and $V_2$ be two natural $2$-dimensional $\GL_2$-modules.
The action of $G$ on $V_1\times V_2$ is analogous to the one
described on $\PxP$ in the beginning of the chapter.
This action extends to the module $W=\Sym^2(V_1)\otimes\Sym^2(V_2)$
such that the map $\varphi_0'\co V_1\times V_2\to W$,
$(v_1,v_2)\mapsto v_1^2\otimes v_2^2$
is $G$-equivariant. 
If $v_{11}, v_{12}$ is a basis of $V_1$, $v_{21}, v_{22}$ a basis of $V_2$,
we take the following basis of $W$:
\begin{align*}
  & v_{11}^2\otimes v_{21}^2, 
    \ 2v_{11}^2\otimes v_{21}v_{22}, 
    \ v_{11}^2\otimes v_{22}^2,
  \ 2v_{11}v_{12}\otimes v_{21}^2, 
    \ 4v_{11}v_{12}\otimes v_{21}v_{22},
    \ 2v_{11}v_{12}\otimes v_{22}^2,&\\
  & v_{12}^2\otimes v_{21}^2,
    \ 2v_{12}^2\otimes v_{21}v_{22},
    \ v_{12}^2\otimes v_{22}^2.&
\end{align*}
Then after comparing with~(\ref{eq:P1xP1antican}) we see that 
the projectivization of the image of $V_1\otimes V_2$ under $\varphi_0'$ 
in coordinates relative to this basis is exactly $X_0$.
Since the image is fixed by $G$, we have a map
$\tilde{\varphi}_0\co G\to G(X_0,k)\subset\GL_9(k)$.
The differential $d\tilde\varphi_0$ takes the Lie algebra $\g_0$ to $\g_0(X_0,k)$. 
By Proposition~\ref{prop:LieConjug} every automorphism of $\g_0(X_0,k)$
is a conjugation by a matrix from $G(X_0,k)$. 
Therefore by Proposition~\ref{AlgVarEq} it follows that the algorithm 
{\tt FindProjectiveEquivalence} can be applied
for finding a parametrization of a surface $X\subset\P^8$ projectively 
equivalent to $X_0$ over $k$.

In the algorithm we have to construct two isomorphisms: an isomorphism
of $\g_0(X,k)$ and $\slxsl$ and isomorphism of two $\slxsl$-modules
afforded by the representations $\g_0(X_0,k)$ and $\g_0(X,k)$.
For deciding and finding a Lie algebra isomorphism we do the following:

\begin{enumerate}
  \item Check whether $\g_0(X,k)$ is semisimple.
  \item Decompose $\g_0(X,k)$ as a sum of two simple subalgebras,
    $\g_0(X,k) = \g_1\oplus \g_2$. 
    (An algorithm for the decomposition can be found in ~\cite{gra6}.)
  \item Find isomorphisms $\g_i\to\mf{sl}_2(k)$, $i=1,2$ 
    as described in Chapter~\ref{ch:csa}.
\end{enumerate}

If any step of the algorithm cannot be carried out, the variety $X$
is not isomorphic to $\PxP$ over $k$.
Note, that if we work over $\Q$, the third step in the algorithm 
reduces to finding rational points on two conics, 
see Chapter~\ref{ch:aut}, Section~\ref{se:conics}.

Now we want to find an isomorphism of the natural 
$\g_0(X_0,k)$- and $\g_0(X,k)$-modules.
In the previous discussion we have seen that the natural $\g_0(X_0,k)$-module
is the $\slxsl$-module $\Sym^2(V_1)\otimes\Sym^2(V_2)$,
where both $V_1$ and $V_2$ are $2$-dimensional natural $\mf{sl}_2$-modules.
It follows that the $\g_0(X_0,k)$-module is irreducible 
and has the highest weight $(2,2)$.

Hence we can decide module equivalence
by checking irreducibility and computing the highest weight.
In the affirmative case, we can again construct a module isomorphism
using highest weight vectors. 

\subsection*{Timings}
The algorithm is implemented in Magma.
Some statistics about timings for surfaces over $\Q$ 
can be found in Table~\ref{tab:PxP}.
For testing, the examples were constructed as follows.
We start with the standard embedding of $\PxP$ into $\P^8$ 
(see~(\ref{eq:P1xP1antican})) which is implicitly given by 20 binomials.
Then a $9\times 9$ matrix is generated, containing
random integer numbers with absolute values up to a given maximal
number (this is written in the first column of Table~\ref{tab:PxP}).
This matrix is then used as the matrix of a linear transformation of 
projective space obtaining so a different system of implicit equations.
\begin{table}[!htb]
  \begin{center}
  \begin{tabular}{|r|r|r|r|r@{.}l|r@{.}l|r@{.}l|}
    \hline
    perturb & eqns max & LA size & prm size & \multicolumn{2}{|c|}{time} & 
    \multicolumn{2}{|l|}{LA time} & \multicolumn{2}{|l|}{conic time}\\
    \hline
    1       & 4        & 11      & 18       & 4&56   & 4&49   & 0&00 \\
    5       & 73       & 47      & 70       & 21&93  & 21&66  & 0&03 \\
    10      & 255      & 55      & 84       & 28&46  & 28&11  & 0&09 \\
    50      & 5026     & 84      & 130      & 48&75  & 48&15  & 0&22 \\
    100     & 25304    & 111     & 166      & 61&02  & 60&15  & 0&34 \\
    300     & 225440   & 134     & 200      & 75&86  & 73&00  & 2&14 \\
            & 208199   & 136     & 204      & 89&52  & 73&15  & 15&77 \\
    400     & 335499   & 143     & 213      & 77&99  & 76&31  & 0&93 \\
            & 418185   & 141     & 210      & 152&21 & 77&91  & 73&56 \\
            & 545728   & 140     & 208      & 482&69 & 74&50  & 407&53 \\
    500     & 720193   & 147     & 222      & 91&11  & 82&24  & 8&10 \\
            & 525179   & 145     & 216      & 80&95  & 78&91  & 1&29 \\
            & 546787   & 143     & 218      & 176&13 & 78&51  & 96&96 \\
    \hline
  \end{tabular}
  \end{center}
  \small
  \begin{tabular}{r@{ -- }p{11.3cm}}
    perturb & maximal entry allowed in perturbation matrix,\\
    eqns max & maximal absolute value of the coefficients 
               in the implicit equations,\\
    LA size & maximal length of the numerator/denominator 
              of the structure constants of the Lie algebra,\\
    prm size & maximal length of the numerator/denominator 
               of the coefficients in the parametrization,\\
    time & total time (in seconds) needed for parametrizing,\\
    LA time & time (in seconds) needed for finding the Lie 
              algebra (is a part of ``time'' in the previous column).\\
    conic time & time (in seconds) needed for finding rational points
                 on two conics constructed to identify two summands 
                 $\mf{sl}_2(\Q)$ (is a part of ``time''). \\
  \end{tabular}
  \caption{Parametrizing $\PxP$.}\label{tab:PxP}
\end{table}

For a ``small'' perturbation, almost the whole computational time is spent 
for finding the Lie algebra of the surface. As the coefficients of 
the linear transformation grow, finding a rational point on the conic 
starts to play the main role in the time complexity.

\subsection{Proper twists of $\PxP$}

If a surface is isomorphic to $\PxP$ over an extension of $k$ but not 
over $k$ itself, it still can have a proper parametrization over $k$. 
Here we investigate this case.

\begin{thm}\label{thm:CxC}
  Let $X\cong C_1\times C_2$ be defined over $k$, 
  where $C_1$ and $C_2$ are twists of $\P^1$.
  Then $X$ has a $k$-parametrization if and only if 
  $C_1\cong\P^1$ and $C_2\cong\P^1$ over $k$.
\end{thm}

\begin{proof}
  Assume that $X$ has a parametrization. Then in particular there is
  a $k$-rational point $p\in X$.
  The two projections $\pi_i\co X\to C_i$ ($i=1,2$) are defined over $k$,
  therefore give $k$-rational points $\pi_1(p)\in C_1$ and 
  $\pi_2(p)\in C_2$. 
  It follows that $C_1\cong\P^1$ and $C_2\cong\P^1$ over $k$.
  The other direction is trivial.
\end{proof}

Theorem~\ref{thm:CxC} implies that if we are looking
for the proper twists of $\PxP$ with a $k$-rational parametrization,
we have to investigate varieties which do not decompose as a product
over $k$. 
Let us therefore study, how is the fact, that the variety
is a product over $k$, reflected by its Lie algebra.

\begin{thm}
  Let $X$ be a twist of $\PxP$ defined over $k$. 
  Then the following are equivalent:
  \begin{enumerate}
    \item[(i)]  $X$ is a product of two twists of $\P^1$,
    \item[(ii)] the Lie algebra $\g_0(X,k)$ is a direct sum
      of two twists of $\mf{sl}_2$.
  \end{enumerate}
\end{thm}

\begin{proof}
  (i)$\Rightarrow$(ii): 
  If $X\cong C_1\times C_2$, then $\Aut_0(X)$ is the direct
  product of the two normal subgroups $\Aut(C_1)$ and $\Aut(C_2)$. It follows
  that $\g_0(X,k) = \g_0(C_1,k)\oplus \g_0(C_2,k)$.

  (ii)$\Rightarrow$(i): 
  If $X$ is not a product, we set $k'$ to be 
  a (finite) Galois extension of $k$ such that $X(k')\cong \PxP$. 
  Then $\Pic X(k')\cong\Z^2$. 
  For any $\sigma\in G$ the image $\sigma(D)$ of a divisor $D$ is again
  a divisor and if $D_1\sim D_2$ then $\sigma(D_1)\sim\sigma(D_2)$,
  so $G$ acts on $\Pic X(k')$.
  Now we use the following a bit technical
  \begin{lem}\label{le:interchange}
    Let $X$ be a twist of $\PxP$ defined over $k$ which is not a product
    and let $k'$ be a finite Galois extension of $k$ such that $X(k')\cong\PxP$.
    Then there is $\sigma$ in $\Gal(k'|k)$ interchanging the divisor classes 
    $(1,0)$ and $(0,1)$.
  \end{lem}
  Since $\Gal(k'|k)$ interchanges the two classes defining the two projections,
  it also interchanges the two normal subgroups of $\Aut_0(X(k'))$ of dimension~3,
  and hence it also interchanges the two ideals of $\g_0(X,k')$.
  It follows that these ideals are not defined over $k$, hence $\g_0(X,k)$
  is simple.
\end{proof}

\begin{proof}[Proof of Lemma~\ref{le:interchange}]
  The action of $\Gal(k'|k)$ preserves the intersection numbers and maps
  effective divisors again to effective divisors.
  Therefore the divisor classes $(1,0)$ and $(0,1)$ can either be
  fixed or interchanged.
  We claim, that there is $\sigma$ in $\Gal(k'|k)$ interchanging these two classes.
  For, suppose by contradiction that every $\sigma\in\Gal(k'|k)$ fixes $(1,0)$.
  Take a divisor $D\in\Div X(k')$ such that its class $[D]=(1,0)$.
  Then the divisor $\sum_{\sigma\in\Gal(k'|k)}\sigma(D)$ is fixed by $\Gal(k'|k)$
  and its class is $(|\Gal(k'|k)|,0)$.
  Therefore the associated map is a projection over $k$ to a twist of $\P^1$,
  contradicting our assumption that $X$ is not a product.
\end{proof}

The {\em centroid} of a Lie algebra $\g$ is the centralizer
of $\ad \g$ in $\mf{gl}(\g)$. The centroid of a Lie algebra $\g$
is denoted by $\Gamma(\g)$.

\begin{lem}\label{le:centroid}
  Let $X$ be a twist of $\PxP$ defined over $k$
  and let $\Gamma_0$ be the centroid of $\g_0(X,k)$.
  \begin{enumerate}
    \item[(i)]  If $X$ is a product, then $\Gamma_0 = k^2$.
    \item[(ii)] If $X$ is not a product, then $\Gamma_0$ is
      a quadratic field extension of $k$ 
      and $X(\Gamma_0)$ is a product.
  \end{enumerate}
\end{lem}

\begin{proof}
(i)
  By~\cite{jacLA}, \S X, Theorem~1, the centroid of a simple Lie algebra is a field.
  If $\g$ is a twist of $\mf{sl}_2(k)$, we conclude then that $\Gamma(\g) \cong k$.
  Write $\g_0(X,k) = \g_1\oplus \g_2$.
  Let $T\co\g_0(X,k)\to \g_0(X,k)$ be a linear map commuting with all $\ad x$ 
  for $x\in \g_0(X,k)$. Then $T(\g_i)\subseteq \g_i$, $i=1,2$.
  For, let $x\in \g_1$ and write $T(x) = y_1+y_2$, $y_i\in \g_i$.
  If $y_2\ne 0$ then take $z\in \g_2$ such that $[z,y_2]\ne 0$.
  Then $\ad z\circ T(x) = [z,y_2] \ne 0$.
  But $T\circ \ad z(x) = T([z,x]) = 0$, a contradiction.
  If follows, that $\Gamma(\g_1\oplus \g_2) = \Gamma(\g_1)\oplus\Gamma(\g_2) = k^2$,
  since both $\g_1$ and $\g_2$ are twists of $\mf{sl}_2(k)$.
  
  (ii)
  Let $k'$ be a field extension of $k$ such that $\g_0(X,k')$ decomposes
  as the sum of two ideals.
  Because computing the centralizer commutes with base field extension,
  $\Gamma_0\otimes_k k'$ is equal to the centroid of $\g_0(X,k')$, which is
  $(k')^2$. Hence $\Gamma_0$ is a 2-dimensional (associative) $k$-algebra without
  nilpotent elements. This implies that either $\Gamma_0$ is a quadratic
  field extension of $k$ or $\Gamma_0\cong k^2$ (with componentwise multiplication).
  Since $\g_0(X,k)$ is simple, the first case is true.
  To prove the last statement we use again the commuting property of computing 
  the centroid and extending the field to $\Gamma_0$:
  $\Gamma(\g_0(X,\Gamma_0)) = \Gamma(\g_0(X,k))\otimes_k\Gamma_0 \cong \Gamma_0^2$.
  The two principal idempotents $(1,0)$ and $(0,1)$ in $\Gamma(\g_0(X,\Gamma_0))$ 
  lead to the desired decomposition of $\g_0(X,\Gamma_0)$.
\end{proof}

From Theorem~\ref{thm:CxC} and Lemma~\ref{le:centroid} it follows,
that a proper twist of $\PxP$ can have a $k$-rational parametrization
only if $X$ is not a product over $k$ and it decomposes as a product
of two projective lines over the centroid of its Lie algebra.
In following we will see that this is not only necessary but also
a sufficient condition.

\begin{prop}
  Let $X$ be a twist of $\PxP$ defined over $k$ which is not a product
  such that for the centroid $k'=\Gamma(\g_0(X,k))$ of $\g_0(X,k)$ we have 
  $X(k')\cong\PxP$.
  Then $X$ is $k$-rational.
\end{prop}

\begin{proof}
  By Lemma~\ref{le:centroid}, $k'$ is a quadratic extension of $k$ and
  $\g_0(X,k') \cong \slxsl$ is the decomposition over $k'$.
  Denote by $\sigma$ the generator of $\Gal(k'|k)$.
  By Lemma~\ref{le:interchange} then $\sigma$ interchanges 
  the divisor classes $(1,0)$ and $(0,1)$. 
  Let us denote the projection associated to $(1,0)$ by $\pi$.
  The intersection number of the two classes is 1,
  therefore for a point $p\in\P^1(k')$ we have that 
  $\pi^{-1}(p) \cap \sigma\circ\pi^{-1}(p)$ is a single point, let us 
  denote it by $q$.
  Further we have that $q$ is fixed by $\Gal(k'|k)$, hence $q\in X(k)$.
  Let $\alpha\in k'$ be such that $k' = k(\alpha)$.
  The map mapping $\P^1(k')\to X(k)$ just described is defined by polynomials.
  Its composition with $\P^2(k)\to\P^1(k')$ mapping $(s:t:u)\mapsto(s+\alpha t:u)$
  is a rational map over $k$ and has a $k$-rational inverse.  
\end{proof}

The proposition finishes our analysis of proper twists of $\PxP$.
Together with Theorem~\ref{thm:CxC} we have that a proper twist $X$ of $\PxP$ 
is $k$-rational if and only if the Lie algebra $\g_0(X,k)$ is simple 
and splits completely over its centroid.
The proof of the proposition gives also the last step in the process of parametrizing
rational proper twists of $\PxP$. To sum up, we have the following algorithm:
\begin{enumerate}
  \item Check whether $\g_0(X,k)$ is simple and compute its centroid 
    $k' = \Gamma(\g_0(X,k))$. 
    By Lemma~\ref{le:centroid} (ii), $k'$ is a quadratic extension of $k$; 
    let $\sigma$ denote the generator of the Galois group $\Gal(k'|k)$.
  \item Decompose $\g_0(X,k')$ as a sum of two simple subalgebras,
    $\g_0(X,k') = \g_1\oplus \g_2$. 
  \item Find isomorphisms $\g_i\to\mf{sl}_2(k')$, $i=1,2$ 
    by reducing to a norm equation (see Chapter~\ref{ch:csa}).
  \item Construct an isomorphism of modules of $\g_0(X_0,k')$ and $\g_0(X,k')$ 
    obtaining so a projective equivalence of $X_0$ and $X$ over $k'$
    and also an $k'$-embedding $\varphi\co\PxP\to X$.
  \item For a fixed $(\alpha{:}\beta)\in\P^1(k')$, let $l_{(\alpha{:}\beta)}$ 
    be the line $\varphi(\{(\alpha{:}\beta;\ s{:}t)\mid (s{:}t)\in\P^1(k')\})$.
    Then $\P^1(k')\to X(k)$, 
    $(\alpha:\beta)\mapsto l_{(\alpha{:}\beta)}\cap\sigma(l_{(\alpha:\beta)})$ 
    together with $\P^2(k)\to\P^1(k')$, $(s:t:u)\mapsto(s+\alpha t:u)$
    leads to a parametrization $\P^2(k)\to X(k)$.
\end{enumerate}

\begin{exa}\label{exa:sphere}
  To illustrate the method we find a parametrization of the unit sphere,
  which is a twist of Segre embedding of $\PxP$. 
  For simplicity we omit here the anticanonical embedding.

  The canonical surface $S_0$ is parametrized by 
  $(s_0{:}s_1;\ t_0{:}t_1)\mapsto(s_0t_0: s_0t_1 : s_1t_0 : s_1t_1)$
  and implicitly given by the polynomial $z_0z_3 - z_1z_2$.
  The surface $S$ is given by $z_1^2 + z_2^2 + z_3^2 - z_0^2$. 
  The Lie algebra $\g_0(S,\Q)$ is simple, 
  the centroid $\Gamma(\g_0(S,\Q))$ is the extension $\Q(i)$, $i^2=-1$. 
  Over this field the algebra splits as the sum 
  $\g_0(S,\Q(i)) = \mf{sl}_2(\Q(i))\oplus\mf{sl}_2(\Q(i))$.
  We find an isomorphism of the Lie algebras and afterwards the parametrization
  of $S$ over $\Q(i)$:
  \begin{equation}\label{eq:sphere}
    \left(\begin{array}{@{\ }c@{\ }}
      z_0 \\ z_1 \\ z_2 \\ z_3
    \end{array}\right)
    =
    \left(\begin{array}{@{\ }cccc@{\ }}
      0 & -1 &  1 &  0 \\
      1 &  0 &  0 & -1 \\
      i &  0 &  0 &  i \\
      0 &  1 &  1 &  0
    \end{array}\right)
    \left(\begin{array}{@{\ }c@{\ }}
      s_0t_0 \\ s_0t_1 \\ s_1t_0 \\ s_1t_1
    \end{array}\right)
    =
    \left(\begin{array}{@{\ }c@{\ }}
      -s_0t_1 + s_1t_0 \\
      s_0t_0 - s_1t_1 \\
      is_0t_0 + is_1t_1 \\
      s_0t_1 + s_1t_0
    \end{array}\right).
  \end{equation}
  The line $l_{(\alpha:\beta)}$ is parametrized by fixing
  $(s_0:s_1)$ to $(\alpha:\beta)\in\P^1(\Q(i))$ in~(\ref{eq:sphere}):
  $l_{(\alpha:\beta)} = (\beta t_0 - \alpha t_1 : \alpha t_0 - \beta t_1 :
  i\alpha t_0 + i\beta t_1 : \beta t_0 + \alpha t_1)$.

  Every point $(\alpha:\beta)\in\P^1(\Q(i))$ can be written as
  $(a+ib:c)$ such that $a,b,c\in\Q$. 
  A rational parametrization of the sphere $S$ is then obtained by
  mapping every point $(a:b:c)\in\P^2(\Q)$ into the intersection of
  $l_{(a+ib:c)}$ and its conjugate:
  \begin{eqnarray*}
    l_{(a+ib:c)} &=& (c t_0 - (a+ib) t_1 : (a+ib) t_0 - c t_1 :
    i(a+ib) t_0 + ic t_1 : c t_0 + (a+ib) t_1),\\
    \overline{l}_{(a+ib:c)} &=& (c s_0 - (a-ib) s_1 : (a-ib) s_0 - c s_1 :
    -i(a-ib) s_0 - ic s_1 : c s_0 + (a-ib) s_1),
  \end{eqnarray*}
  obtaining such
  \begin{displaymath}
    (a:b:c) \mapsto l_{(a+ib:c)} \cap \overline{l}_{(a+ib:c)} =
    (c^2 + a^2 + b^2 : 2ac : -2bc : c^2 - a^2 - b^2).
  \end{displaymath}
\end{exa}

\subsection*{Timings}
For testing the algorithm we constructed examples as follows.
We have chosen $d\in\Z$ such that $d\not\in\Q^2$
(given in the first column of Table~\ref{tab:sphere}).
Then the quadratic surface in $\P^3$ given by $z_0^2 - z_1^2 = z_2^2 - dz_3^2$
is isomorphic to $\PxP$ over $\Q(\sqrt{d})$ but not over $\Q$.
We anticanonically embedded the surface over $\Q$ into $\P^8$ 
obtaining such a surface described by 14 binomials and 6 polynomials
with 4 terms. Afterwards we made a linear transformation similar 
to the previous case, just here the generated matrix is 
sparser, to obtain examples solvable in practice. 
Since we have to identify two $\mf{sl}_2$'s over $\Q(\sqrt{d})$,
we have to solve two relative norm equations.
This is very time consuming, therefore we were able to parametrize
only ``small'' examples.
\begin{table}[!htb]
  \begin{center}
  \begin{tabular}{|r|r|r|r|r|r@{.}l|r@{.}l|r@{.}l|}
    \hline
     & perturb & eqns & LA & prm & 
    \multicolumn{2}{|l|}{} & \multicolumn{2}{|c|}{LA} & 
    \multicolumn{2}{|c|}{normeq} \\
    discr & (sparse) & max & size & size & 
    \multicolumn{2}{|c|}{time} & \multicolumn{2}{|c|}{time} & 
    \multicolumn{2}{|c|}{time} \\
    \hline
    -1  & 1  & 3   & 3   & 9    & 2&460   & 0&670   & 1&010 \\
    3   & 1  & 5   & 3   & 23   & 3&620   & 1&030   & 1&700 \\
    8   & 1  & 15  & 5   & 1135 & 211&340 & 1&270   & 123&230 \\
    -1  & 2  & 10  & 5   & 92   & 41&690  & 1&250   & 38&670 \\
    \hline
  \end{tabular}
  \end{center}
  \small
  \begin{tabular}{r@{ -- }p{11.3cm}}
    discr & square of the primitive element used for the construction,\\
    normeq time & time (in seconds) needed for solving two relative 
                  norm equations (is a part of ``time''). \\
  \end{tabular}
  Description of the other columns: as in Table~\ref{tab:PxP}.
  \caption{Parametrizing the proper twists of $\PxP$.}\label{tab:sphere}
\end{table}

\section{Blow-up of the projective plane}

In this section the surface $X$ is the {\em blowup} of the projective plane
in a single point, i.e. an embedding of
\begin{displaymath}
  \Pi=\left\{(x_0:x_1:x_2;\ y_1:y_2)\in\P^2\times\P^1 \mid x_1y_2 = x_2y_1\right\}.
\end{displaymath}
The line $l\subset\Pi$, $l=\{(1:0:0;\ s:t)\}$ is the {\em exceptional line} of $\Pi$
and it is fixed by every automorphism of $\Pi$.
Therefore the group of automorphisms of $\Pi$ is isomorphic
to the subgroup of automorphisms of $\P^2$
fixing the point $(1:0:0)$.
We have the group of matrices representing these automorphisms:

\begin{equation*}
  H = \left\{ g\in{\GL}_3(k)\ \Big|\ g = \left( \begin{array}{@{\ }ccc@{\ }}
      g_{00} & g_{01} & g_{02} \\
      0      & g_{11} & g_{12} \\
      0      & g_{21} & g_{22} 
    \end{array} \right) \right\}.
\end{equation*}
Therefore the Lie algebra of $\Pi$ (after discarding the identity) is isomorphic to
\begin{equation}\label{eq:BlowupLA}
  \mf{h}_0 = \left\{ a\in\mf{sl}_3(k)\ \Big|\ a = \left( \begin{array}{@{\ }ccc@{\ }}
    2a & b_1   & b_2 \\
    0  & c_1-a & c_2 \\
    0  & c_3   & -c_1-a 
  \end{array} \right) \right\}.
\end{equation}

The semisimple part of the Lie algebra $\mf{h}_0$ is isomorphic to $\mf{sl}_2$.
Let $\mf{k}$ denote the Levi subalgebra of $\mf{h}_0$
consisting of those matrices in~(\ref{eq:BlowupLA}) where
$a, b_1$ and $b_2$ vanish.
The Lie algebra $\mf{h}_0$ further contains 
a three-dimensional solvable radical $\Rad\mf{h}_0$
containing a two-dimensional nilradical $\mf{N}$.
If $e_{ij}$ denotes the matrix having a $1$ at
$(i,j)$-th position and zero elsewhere, then 
$\mf{N}$ is spanned by $e_{12}$ and $e_{13}$
and $\Rad\mf{h}_0$ in addition also by $e=2e_{11}-e_{22}-e_{33}$.

Now we want to prove a proposition implying that also for embeddings
of $\Pi$ associated to a very ample divisor, we can use
the algorithm {\tt FindProjectiveEquivalence}.
Before doing that we need some well-known technical lemmata.

\begin{lem}\label{le:LA_nilpElts}
  Let $\g\subseteq\gl_m(k)$ be a semisimple Lie algebra.
  If $x\in\g$ is nilpotent, then also $\ad x$ is nilpotent. 
\end{lem}

\begin{proof}
  \cite{humLA}, Lemma 3.2.
\end{proof}

\begin{lem}\label{le:LA_exp}
  Let $G\subseteq\GL_m(k)$ be an algebraic group and
  let $\g\subseteq\gl_m(k)$ be its Lie algebra
  If $x\in\g$ is nilpotent, then $\exp(x)\in G$.
\end{lem}

\begin{proof}
  \cite{borel}, \S 7.3.
\end{proof}

\begin{lem}\label{le:LA_Inn}
  Let $\g\subset \gl_m(F)$ be a Lie algebra. 
  Let $x\in \g$ be such that both $\ad x$ and $x$ are nilpotent.
  Then $\exp( \ad x)(x') = (\exp x) x'( \exp -x)$.
\end{lem}

\begin{proof}
  \cite{jacLA}, p. 282.
\end{proof}

\begin{prop}
  Every automorphism of $\mf{h}_0$ is a conjugation by a matrix from $H$.
\end{prop}

\begin{proof}
  As before, let $\mf{k}$ be the Levi subalgebra spanned by 
  $h = e_{22}-e_{33}, e_{23}$ and $e_{32}$. We first show that
  every automorphism $\xi$ of $\mf{h}_0$ such that its restriction on $\mf{k}$
  is the identity, acts on $\Rad\mf{h}_0$ by
  $e\mapsto e$, $e_{12}\mapsto\lambda e_{12}$ and $e_{13}\mapsto\lambda e_{13}$
  for some $\lambda\in k$.

  The automorphism $\xi$ leaves both $\mf{N}$ and $\Rad\mf{h}_0$ fixed.
  Therefore $\xi(e) = \alpha_1 e + \alpha_2 e_{12} + \alpha_3 e_{13}$,
  $\xi(e_{12}) = \beta_2 e_{12} + \beta_3 e_{13}$ and
  $\xi(e_{13}) = \gamma_2 e_{12} + \gamma_3 e_{13}$.
  From $\xi([e,h]) = [\xi(e),\xi(h)]$ we get $\alpha_2 = \alpha_3 = 0$.
  Then from $\xi([e,e_{12}]) = [\xi(e),\xi(e_{12})]$ it follows that $\alpha_1 = 1$.
  To finish the proof of the claim we use 
  $\xi([e_{32},e_{12}]) = [\xi(e_{32}),\xi(e_{12})]$,
  $\xi([e_{23},e_{13}]) = [\xi(e_{23}),\xi(e_{13})]$ and
  $\xi([e_{32},e_{13}]) = [\xi(e_{32}),\xi(e_{13})]$.

  Let $A$ be the subgroup of inner automorphisms generated by
  $\exp(\ad x)$ for $x\in\mf{N}$. Since for $x$ in the nilradical
  $\ad x$ is nilpotent (see Lemma~\ref{le:LA_nilpElts}),
  it follows from Lemmata~\ref{le:LA_exp} and \ref{le:LA_Inn} 
  that automorphisms in $A$ are conjugations by a matrix from $H$.
  Let now $\sigma$ be any automorphism of $\mf{h}_0$.
  Then $\sigma(\mf{k})$ is a Levi subalgebra and by Theorem of
  Malcev-Harish-Chandra (cf.~\cite{jacLA},\S III.)
  there is an automorphism $\tau$ of $\mf{h}_0$ from $A$ such that
  $\tau(\mf{k}) = \sigma(\mf{k})$.
  Together with Theorem~\ref{thm:slAut} we obtain the claim
  of the proposition.
\end{proof}

For finding an isomorphism of $\mf{h}_0$ and the Lie algebra $\g_0(X,k)$
of the given surface isomorphic to $\Pi$ it is useful to note
that the nilradical $\mf{N}$ is a 2-dimensional irreducible module
of the Levi subalgebra $\mf{k}$. Thus we have the following algorithm:

\begin{enumerate}
  \item Set $\mf{k}$ to be a Levi subalgebra of $\g_0(X,k)$.
  \item Set $\mf{N}$ to be the nilradical of $\g_0(X,k)$.
  \item Using the action of $\mf{k}$ on $\mf{N}$ construct an isomorphism
    $\mf{sl}_2\to\mf{k}$.
  \item Set $E_{12}\in\mf{N}$ to be a vector generating one-dimensional eigenspace
    of $\ad h$ with eigenvalue $-1$, where $h\in\mf{k}$ is a vector generating 
    a split Cartan subalgebra of $\mf{k}$.
  \item Similarly one can find the images $E_{13}$ and $E$
    of elements $e_{13}$ and $e$ in the Lie algebra $\mf{h}_0$.
    (Note that this step is not necessary for finding a module isomorphism
    as will become clear soon.)
\end{enumerate}

Next we analyze the 9-dimensional module afforded by the 
anticanonical embedding of $\Pi$.
This  embedding is given by 
\begin{displaymath}
  (x_0{:}x_1{:}x_2;\ y_1{:}y_2)\mapsto
    (x_0^2y_1 : x_0^2y_2 : 
     x_0x_1y_1 : x_0x_1y_2 : x_0x_2y_2 :
     x_1^2y_1 : x_1^2y_2 : x_1x_2y_2 : x_2^2y_2)
\end{displaymath}
where $x_1y_2 = x_2y_1$,
and $\Pi$ embedded in this way is parametrized by
\begin{equation}\label{eq:blowup}
  (s:t:u) \mapsto (s^2t:s^2u : st^2:stu:su^2 : t^3:t^2u:tu^2:u^3),
\end{equation}
The image of this embedding will be denoted as usual by $X_0$.

Let $V$ be a $3$-dimensional vector space with basis $v_0,v_1,v_2$. 
Consider the symmetric power $\Sym^3(V)$ with the basis
$v_0^3$, $3v_0^2v_1$, $3v_0^2v_2$, 
$3v_0v_1^2$, $6v_0v_1v_2$, $3v_0v_2^2$, 
$v_1^3$, $3v_1^2v_2$, $3v_1v_2^2$, $v_2^3$.
Let $\varphi_0^\prime:V\to \Sym^3(V)$ be given by 
$\varphi_0^\prime(v) = v^3$.

Let $G=\GL_3(\Q)$ act naturally on $V$.
Then $G$ acts also on $\Sym^3(V)$, making $\varphi_0'$ $G$-equivariant.

Let $U$ be the subspace of $\Sym^3(V)$ spanned by $v_0^3$.
Let $\pi:\Sym^3(V)\to\Sym^3(V)/U = W$ be the projection discarding 
the coordinate at $v_0^3$.
For $\varphi_0 = \pi\circ\varphi_0^\prime$ we have 
that $X_0$ is the projectivization of $\varphi_0(V)$.
Note that $H$ is exactly the stabilizer of $U$ in $G$
therefore we have a well-defined action of $H$ 
on the whole module $W$.

\begin{lem}\label{le:BlowupModule}
  As a $\mf{k}$-module, $W$ decomposes as a direct sum
  $W = W_2 \oplus W_3 \oplus W_4$,
  where $W_i$ is an $i$-dimensional irreducible $\mf{k}$-module.
  As $\mf{h}_0$-module, $W$ is irreducible.
\end{lem}

\begin{proof}
  When restricting to the Levi subalgebra $\mf{k}$, 
  the module $\Sym^3(V)$ (see the discussion before the lemma) 
  becomes an $\mf{sl}_2$-module and as such
  decomposes as a sum of four irreducible modules:
  $W_1 = U$,
  $W_2$ is the module spanned by $3v_0^2v_1, 3v_0^2v_2$ 
  and is a 2-dimensional $\mf{sl}_2$-module $k^2$,
  $W_3$ is spanned by $3v_0v_1^2, 6v_0v_1v_2, 3v_0v_2^2$
  and isomorphic to $\Sym^2(k^2)$, and lastly
  $W_4$ is spanned by $v_1^3, 3v_1^2v_2, 3v_1v_2^2, v_2^3$
  and isomorphic to $\Sym^3(k^2)$.
  It follows that $W$ as $\mf{sl}_2$-module decomposes 
  into the sum $W_2 \oplus W_3 \oplus W_4$.

  To prove the last assertion of the lemma,
  let us take any $b\in\mf{N}$, $b=b_1e_{12} + b_2e_{13}$.
  So if $w\in W_4$ is a basis vector, $w=v_1^iv_2^{3-i}$, then 
  $b\cdot w \in \Span\{v_0v_1^{i-1}v_2^{3-i}, v_0v_1^iv_2^{2-i}\}\subset W_3$.
  Similarly for $w\in W_3$ one gets $b\cdot w\in W_2$.
\end{proof}

Let $\psi\co W\to W$ be an isomorphism of $\g_0(X_0,k)$-modules. 
Then $\psi$ restricted to $W_i$ is multiplication by a scalar $\lambda_i$. 
Let $b=e_{12}\in\mf{N}$, and $w_4=v_1^3\in W_4$. 
Since $b\cdot v_1=v_0$, it follows $b\cdot w_4 = 3v_0v_1^2\in W_3$. 
Hence $\psi(b\cdot w_4)=\lambda_3b\cdot w_4$. 
On the other hand, $\psi(b\cdot w_4) = b\cdot \psi(w_4)=\lambda_4 b\cdot w_4$. 
We deduce that $\lambda_4=\lambda_3$. 
In the same way we find out that $\lambda_3=\lambda_2$,
so that $\psi$ is multiplication by a nonzero scalar. 

Now we find an isomorphism of $\g_0(X_0,k)$- and $\g_0(X,k)$-module
for a given surface $X$ as follows.
\begin{enumerate}
  \item Let $\mf{k}\subset\g_0(X,k)$ be a Levi subalgebra, so $\mf{k}\cong\sl_2(k)$.
  \item Using weight vectors find the decomposition 
    $W' = W'_2\oplus W'_3\oplus W'_4$ of the natural $\mf{k}$-module,
    where $[W_i':k] = i$.
  \item Find isomorphisms $\psi_i\co W_i\to W'_i$, $(i=2,3,4)$ of $\mf{sl}_2$-modules.
  \item Let $E_{12}$ be the image of the matrix $e_{12}$ in $\g_0(X,k)$,
    so $E_{12}$ is in the nilradical of $\g_0(X,k)$.
  \item $E_{12}$ maps the highest weight vector of $W_i$ (resp.~$W'_i$) 
    to the highest weight vector of $W_{i-1}$ (resp.~$W'_{i-1}$), $(i=3,4)$. 
    Use this to finish the construction of an isomorphism 
    of the natural $g_0(X_0,k)$- and $g_0(X,k)$-modules.
\end{enumerate}

\subsection*{Timings}
We tried our algorithm on examples which we constructed from the 
canonical surface (given by the binomial ideal with 20 generators) 
by a linear transformation of the projective space. 
The randomly generated matrix of the transformation has integral
entries with the given maximal absolute value
(the first column in Table~\ref{tab:blowup}).
We see that almost the whole time is spent for finding the Lie algebra
of the surface.
\begin{table}[htb]
  \begin{center}
  \begin{tabular}{|r|r|r|r|r@{.}l|r@{.}l|}
    \hline
    perturb & eqns max & LA size & prm size & 
    \multicolumn{2}{|c|}{time} & \multicolumn{2}{|l|}{LA time} \\
    \hline
    1       & 4        & 10      & 46       & 4&43   & 4&23 \\
    5       & 85       & 47      & 211      & 21&25  & 20&76 \\
    10      & 280      & 59      & 266      & 28&21  & 27&58 \\
    50      & 6372     & 93      & 424      & 51&66  & 50&43 \\
    100     & 26625    & 103     & 475      & 58&08  & 56&84 \\
    500     & 599186   & 145     & 666      & 82&81  & 80&89 \\
    1000    & 1926906  & 159     & 724      & 91&26  & 89&11 \\
    5000    & 60259495 & 207     & 957      & 118&99 & 115&94 \\
    10000   & 246171712& 219     & 1008     & 129&49 & 126&24 \\
    \hline
  \end{tabular}
  \end{center}
  \small
  Description of the columns: as in Table~\ref{tab:PxP}.
  \caption{Parametrizing blow-ups of $\P^2$.}\label{tab:blowup}
\end{table}

\chapter{Del Pezzo surfaces of degree 9}
Another surfaces, where we used Lie algebras for parametrizing,
are Del Pezzo surfaces of degree~9.
We not only give a parametrization algorithm
but also analyze in more detail the algebra associated to the surface
constructed during the computation.

\section{Severi-Brauer varieties}

{\em Severi-Brauer varieties} may be defined as twists of a projective space, 
i.e.~they are varieties $V$ such that $V\otimes k'\cong\P^n(k')$ for
some extension $k'$ of $k$.

Though geometrically are these varieties very simple, arithmetically they
are rather interesting objects. An alternative definition, giving also more insight,
relates the variety to a central simple algebra.
Namely, an $n$-dimensional Severi-Brauer variety is the set of $(n+1)$-dimensional 
left ideals of a central simple algebra of degree $n+1$ .
The variety associated to the algebra $A$ will be denoted by $\mathcal{V}_A$.
For a survey on the topic see~\cite{jahnel}

If $A\otimes_k k'\cong\M_{n+1}(k')$, then we say that $k'$ is a 
{\em splitting field} of $A$.
In such case the variety $\mathcal{V}_A$ is over the field $k'$ isomorphic 
to the $n$-dimensional projective space.
There is a $k'$-rational point on $\mathcal{V}_A$ if and only if $k'$ is a splitting
field of $A$.
For proofs of these properties see~\cite{jacFDDA}, \S~3.5.
Note that a Severi-Brauer variety is $k$-isomorphic to a projective space
if and only if it contains a $k$-rational point.

\begin{exa}\label{exa:quaternions}
  We will construct the Severi-Brauer variety corresponding to the
  quaternion algebra $\HH$ over $\R$.

  Since $\HH$ is a division algebra, it has no one-sided nontrivial 
  ideals, therefore there are no real points on $V_\HH$.
  The algebra splits over $\CC$ as follows:
  \begin{displaymath}
    \mathbf{1} \mapsto \left( \begin{array}{@{\ }cc@{\ }}
      1 & 0 \\
      0 & 1 
    \end{array} \right),\quad
    \mathbf{i} \mapsto \left( \begin{array}{@{\ }rr@{\ }}
      i & 0 \\
      0 & -i 
    \end{array} \right),\quad
    \mathbf{j} \mapsto \left( \begin{array}{@{\ }rr@{\ }}
       0 & 1 \\
      -1 & 0 
    \end{array} \right),\quad
    \mathbf{k} \mapsto \left( \begin{array}{@{\ }cc@{\ }}
      0 & i \\
      i & 0 
    \end{array} \right).
  \end{displaymath}
  Two-dimensional left ideals in $\HH\otimes\CC$ are parametrized by
  points of the projective line $\P^1(\CC)$. These ideals are of the form
  \begin{displaymath}
    \mathcal{L}_{\alpha:\beta}\ =\ \left\{\left( \begin{array}{cc}
        \alpha x & \beta x \\
        \alpha y & \beta y 
      \end{array} \right)\Big|\ x,y\in\CC\right\}\ =\ {\Span}_\CC\left\{
      \left(\begin{array}{@{\ }cc@{\ }}
        \alpha & \beta \\
        0      & 0
      \end{array}\right),
      \left(\begin{array}{@{\ }cc@{\ }}
        0      & 0 \\
        \alpha & \beta
      \end{array}\right)
    \right\}.
  \end{displaymath}
  In the basis $\mathbf{1}, \mathbf{i}, \mathbf{j}, \mathbf{k}$ of $\HH$ we can write
  \begin{displaymath}
    \mathcal{L}_{\alpha:\beta} = {\Span}_\CC\left\{
      \alpha(\mathbf{1}-i\mathbf{i}) + \beta(\mathbf{j}-i\mathbf{k}),
      \ -\alpha(\mathbf{j}+i\mathbf{k}) + \beta(\mathbf{1}+i\mathbf{i})
    \right\}.
  \end{displaymath}
  After embedding this set into Grassmannian $G(2,4)$ and afterwards 
  into $\P^5$ we have
  \begin{displaymath}\renewcommand{\arraystretch}{0.9}
    \left(\begin{array}{@{\ }cc@{\ }}
      \alpha   & \beta \\
      -i\alpha & i\beta \\
      \beta    & -\alpha \\
      -i\beta  & -i\alpha
    \end{array}\right)
    \mapsto
    \left(\begin{array}{@{\ }c@{\ }}
      p_{12} \\
      p_{13} \\
      p_{14} \\
      p_{23} \\
      p_{24} \\
      p_{34}
    \end{array}\right)
    = 
    \left(\begin{array}{@{\ }c@{\ }}
      2i\alpha\beta \\
      -\alpha^2 - \beta^2 \\
      -i\alpha^2 + i\beta^2 \\
      i\alpha^2 - i\beta^2 \\
      -\alpha^2 - \beta^2 \\
      -2i\alpha\beta
    \end{array}\right)
    \mapsto
    \left(\begin{array}{@{\ }c@{\ }}
      p_{12} \\
      p_{13} \\
      p_{14}
    \end{array}\right)
    = 
    \left(\begin{array}{@{\ }c@{\ }}
      2i\alpha\beta \\
      -\alpha^2 - \beta^2 \\
      -i\alpha^2 + i\beta^2
    \end{array}\right),
  \end{displaymath}
  where the last map is a projection onto linearly independent coordinates.
  The implicit equation of the Severi-Brauer curve corresponding to the 
  quaternion algebra therefore is
  \begin{displaymath}
    p_{12}^2 + p_{13}^2 + p_{14}^2 = 0.
  \end{displaymath}
\end{exa}

\section{Parametrizing Del Pezzo surfaces of degree 9}

Del Pezzo surfaces of degree 9
are anticanonically embedded Severi-Brauer surfaces,
so over an algebraic closure they are isomorphic to the projective plane.

The anticanonical embedding $\varphi_0$ of the projective plane into $\P^9$
is given by
\begin{equation}\label{eq:P2antican}
  (s:t:u)\mapsto
    (s^3:t^3:u^3:
     s^2t:t^2u:u^2s:
     st^2:tu^2:us^2:
     stu).
\end{equation}
The image $\varphi_0(\P^2)$ we denote $X_0$. Implicitly this surface
is given by 27 quadratic forms over the rationals.

For the projective plane we have that $\Aut(\P^2) = \Aut_0(\P^2) = \PGL_3(k)$.
Therefore the Lie algebra $\g_0(\P^2,k) = \mf{sl}_3(k)$.
In this example we cannot use Proposition~\ref{AlgVarEq} since the outer
automorphism of $\mf{sl}_3(k)$ mapping the matrix $x$ into $-x^t$ is
not a conjugation by any $g\in\GL_3(k)$ 
(here $\mf{sl}_3(k)$ is represented by $3\times 3$ matrices of trace~0).

\begin{lem}\label{le:Sym3}
  The 10-dimensional $\g_0(X_0,k)$-module 
  is an irreducible $\mf{sl}_3$-module.
\end{lem}

\begin{proof}
  Let $V=k^3$ and let $\varphi_0'$ be the map $V\to W=\Sym^3(V)$, $v\mapsto v^3$.
  The group $\GL_3(k)$ acts on $V$ naturally, and on $W$ by 
  $g\cdot (v_1 v_2 v_3) = (g v_1)(g v_2)(g v_3)$.
  The map $\varphi_0'$ is $\GL_3$-equivariant, so the image
  of $V$ under $\varphi_0'$ is left by $\GL_3(k)$ invariant.
  If we fix bases in $V$ and $W$, we will obtain a representation
  of $\GL_3(k)$ in $\GL_{10}(k)$.

  Let $v_0,v_1,v_2$ be the standard basis of $V=k^3$. 
  In $W=\Sym^3(V)$ we take the basis
  \begin{equation}\label{eq:Sym3basis}
    v_0^3, v_1^3, v_2^3,
    3v_0^2v_1, 3v_1^2v_2, 3v_2^2v_0,
    3v_0v_1^2, 3v_1v_2^2, 3v_2v_0^2,
    6v_0v_1v_2,
  \end{equation} 
  and by $\tilde\varphi_0$ we denote the corresponding faithful 
  representation $\GL_3(k)\to\GL_{10}(k)$.
  Now if we write the image of $V$ under $\varphi_0'$ in coordinates
  relative to the basis~(\ref{eq:Sym3basis}), we see that
  $\P(\varphi_0'(V)) = X_0$.
  So $\tilde\varphi_0$ maps $\GL_3(k)$ to $G(X_0,k)$.
  Since both groups are algebraic and isomorphic,
  $\tilde\varphi_0$ is an isomorphism.
  Then $d\tilde\varphi_0\co\mf{sl}_3(k)\to \g_0(X_0,k)$ is
  an isomorphism of Lie algebras.
  It is a representation of $\mf{sl}_3(k)$ with underlying module
  $\Sym^3(V)$.
  Therefore the 10-dimensional natural $\g_0(X_0,k)$-module is 
  irreducible with the highest weight $(3,0)$.
\end{proof}

The group $\Aut(\P^2)$ is the projectivization of the linear group $\GL_3(k)$.
We take the following Borel subgroup of $\GL_3(k)$:
\begin{displaymath}
  B = \left\{\left(\begin{array}{ccc}
      \ast & \ast & \ast \\
         0 & \ast & \ast \\
         0 &    0 & \ast 
  \end{array}\right)\in{\GL}_3(k)\right\},
\end{displaymath}
i.e.~the group of upper triangular matrices.
As in the case of conics, $\GL_3(k)$ acts transitively on $\P^2$ and hence
via $\varphi_0$ also on $X_0$. Therefore $X_0$ is a single orbit of $\GL_3(k)$
in $\P^9$.

\begin{lem}\label{le:P2Borel}
  The Borel subgroup $B$ of the group $\GL_3(k)$ has via the action 
  $\varphi_0$~(\ref{eq:P2antican}) a unique fixed point in $\P^9$.
\end{lem}

\begin{proof}
  As in the proof of Lemma~\ref{le:Sym3}, we have that 
  $\P^9 = \P(\Sym^3(V))$, $V = k^3$.
  The Borel subgroup $B$ fixes the unique line in $\Sym^3(V)$,
  namely the line spanned by the highest weight vector.
  Hence in $\P^9$ there is a unique point fixed by $B$.
\end{proof}

Let $H\subset\mf{sl}_3(k)$ be a fixed Cartan subalgebra with basis $h_1,h_2$ 
which are part of a Chevalley basis of $\mf{sl}_3(k)$. Let $\tau$ 
be a fixed automorphism of $\mf{sl}_3(k)$, such that $\tau(h_1)=h_2$
and $\tau(h_2)=h_1$ (such an automorphism exists by \cite{gra6}, \S 5.11).

\begin{lem}
  Let $d\tilde\varphi_0\co\mf{sl}_3(k)\to \g_0(X_0,k)$ be the 
  isomorphism of Lie algebra as above.
  Let $X\subset\P^9$ be projectively equivalent to $X_0$ over an extension of $k$.
  If there is a $k$-isomorphism $\rho\co\mf{sl}_3(k)\to \g_0(X,k)$,
  then either $d\tilde\varphi_0$ and $\rho$ or $d\tilde\varphi_0$ and $\rho\circ\tau$
  are isomorphic representations.
  Further, if $\mu$ is the module isomorphism, 
  then $\mu$ modulo scalar multiplication
  is a projective equivalence of $X_0$ and $X$ defined over $k$.
\end{lem}

\begin{proof}
  Since modules of $\g_0(X_0,\overline{k})$ and $\g_0(X,\overline{k})$ are 
  isomorphic (Proposition~\ref{prop:lieAlgVarIso}), $\g_0(X,k)$-module is also
  irreducible.
  There are exactly two irreducible $\mf{sl}_3$-modules of dimension $10$. 
  Let $h_1,h_2$ be as before the lemma.
  We represent a weight $\lambda\in H^*$ by the tuple 
  $(\lambda(h_1), \lambda(h_2))$. 
  Then the two irreducible $\mf{sl}_3$-modules of dimension $10$ 
  have highest weights $(3,0)$ and $(0,3)$ respectively.
  By composing $d\tilde\varphi_0$ with $\tau$ we change the highest weight
  of the corresponding module (from $(3,0)$ to $(0,3)$ or vice versa).
  Therefore, after possibly composing $\rho$ with $\tau$ we have that 
  the two representations have the same highest weight, and hence 
  are isomorphic.
  Now the last assertion follows from Lemma~\ref{le:P2Borel} and 
  Proposition~\ref{prop:ProjEqBorel}.
\end{proof}

Now we have an algorithm which for a surface $X$ in $\P^9$ given by 
quadratic forms over $k$ decides whether it is a Del Pezzo surface
of degree~9 containing a $k$-rational point, and if so, 
finds a $k$-parametrization.

\begin{enumerate}
  \item Find the Lie algebra $\g_0(X,k)$ as described in 
    Chapter~\ref{ch:aut}, Section~\ref{se:computeL}.
  \item Find an associative algebra $A$ such that 
    $A_\mathrm{Lie}\cong \g_0(X,k)\oplus k$ using the algorithm 
    {\tt EnvelopingAlgebra} 
    in Chapter~\ref{ch:preliminaries}, Section~\ref{se:algebras}.
  \item\label{step:csa} Construct an isomorphism $\rho'\co\M_3(k)\to A$ over $k$
    (see Chapter~\ref{ch:csa}),
    inducing so an isomorphism $\rho$ of Lie algebras $\mf{sl}_3(k)\to \g_0(X,k)$.
    If no such isomorphism exists, then $X$ is not isomorphic to $X_0$ over $k$.
  \item Check whether the modules $\rho$ and $d\tilde\varphi_0$ or 
    $\rho\circ\tau$ and $d\tilde\varphi_0$ are isomorphic. 
    If none of these isomorphisms exists, then $X$ was not 
    Del Pezzo surface of degree 9. 
  \item In the affirmative case in the previous step let 
    $M\in\M_{10}(k)$ be the matrix describing the module isomorphism.
    If $M$ does not map $X_0$ to $X$ then 
    $X$ was not Del Pezzo surface of degree 9.
\end{enumerate}

The construction in this section generalizes to finding an isomorphism 
of Severi-Brauer varieties of arbitrary dimension $n$.
The anticanonical class in $\P^n$ has degree $n+1$, so the associated
linear system is the space of all forms of degree $n+1$.
Therefore the involved module is the module of $(n+1)$-th symmetric powers 
of $k^{n+1}$, which is an irreducible $\mf{sl}_{n+1}$-module.
As a Borel subgroup we can take the group of upper triangular matrices
in $\GL_{n+1}(k)=G(\P^n,k)$. 
It fixes a unique point in the irreducible representation,
namely the highest weight vector.
The hardest problem is finding an isomorphism 
of an associative algebra and $\M_{n+1}(k)$. 
In Chapter~\ref{ch:csa} we give algorithms for algebras up to degree~4.

\subsection*{Testing the algorithm and timings}
In the step~(\ref{step:csa}) of the algorithm, the finding isomorphism
of an algebra and the full matrix algebra is reduced to solving 
a norm equation over a cubic field extension, see Chapter~\ref{ch:csa}.
When constructing examples by perturbing coefficients in the
implicit equations of the standard embedding~(\ref{eq:P2antican})
as in the case of Del Pezzo surfaces of degree~8, 
we were able to construct the norm equation but it 
was too hard to solve.
Therefore for generating Severi-Brauer surfaces
we used the {\em Galois descent}, a very efficient method based
on Galois cohomology. 
For expositions of the Galois descent see for example \cite{jacFDDA, krashen},
an explicit method can be found in~\cite{kang}. 
In this case the associative algebra corresponding to the surface is found 
in a very neat shape: we can very easily construct the norm equation
and moreover this equation is easy to solve,
even when the coefficients in the implicit equations of the surface 
are very large as is shown in Table~\ref{tab:P2}.

When using the Galois descent, one has first to take a cubic field
extension of the base field (in our case $\Q$) which is Galois.
This is equivalent to adjoining a root of a cubic irreducible polynomial
such that the discriminant of the polynomial is a square in the base field.
When generating such polynomial, we bound the absolute value of the
coefficients in the polynomial. 
The bound is given in the first column of the table.
\setcounter{table}{3}
\begin{table}[!htb]
  \begin{center}
  \begin{tabular}{|r|r|r|r|r@{.}l|r@{.}l|r@{.}l|}
    \hline
    field cfs & eqns size & LA size & prm size & \multicolumn{2}{|c|}{time} & 
    \multicolumn{2}{|l|}{LA time} & \multicolumn{2}{|l|}{normeq time}\\
    \hline
    100     & 4        & 3       & 18       & 1&94   & 1&48   & 0&16 \\
    500     & 6        & 7       & 23       & 2&56   & 1&85   & 0&21 \\
    1000    & 7        & 8       & 47       & 3&26   & 2&05   & 0&58 \\
    5000    & 11       & 13      & 42       & 3&63   & 2&62   & 0&27 \\
    10000   & 9        & 9       & 49       & 3&17   & 2&26   & 0&25 \\
    50000   & 10       & 10      & 48       & 3&14   & 2&29   & 0&27 \\
    100000  & 12       & 14      & 59       & 3&76   & 2&70   & 0&32 \\
    500000  & 13       & 15      & 82       & 4&33   & 2&60   & 0&75 \\
    \hline
  \end{tabular}
  \end{center}
  \small
  \begin{tabular}{r@{ -- }p{11.3cm}}
    field cfs & maximal entry allowed in the minimum polynomial
                of an algebraic element generating the extension,\\
    eqns size & maximal length of the numerator/denominator 
               of the coefficients in the implicit equations,\\
  \end{tabular}
  Description of the other columns: as in Tables~\ref{tab:PxP} and \ref{tab:sphere}.
  \caption{Parametrizing $\PxP$.}\label{tab:P2}
\end{table}

\section{More on the algebra of a Severi-Brauer variety}

In the previous section we have seen that a given Severi-Brauer variety 
was a split one if and only if the associative algebra constructed 
in our algorithm is isomorphic to a matrix algebra.
So in the split case we indeed found the algebra $A$ corresponding to 
the variety $\mathcal{V}_A$. Here we investigate the associative
algebra also in the case of non-split variety.

Using the construction of a variety as the set of left
ideals of a central simple algebra, we will show that for curves, 
the associative algebra constructed by our method is always the algebra
corresponding to the given Severi-Brauer curve.
For Severi-Brauer varieties of higher dimensions we are able to prove
a similar but slightly weaker assertion.

\begin{prop}\label{prop:autIso}
  Let $A$ be a central simple algebra over $k$ and let
  $\mathcal{V}_A$ be the associated Severi-Brauer variety.
  Then the groups of $k$-automorphisms of $A$ and $\mathcal{V}_A$
  are isomorphic. 
\end{prop}

\begin{proof}
  By $d$ we will denote the degree of $A$
  and let $k'$ be a splitting field of $A$, 
  so $A_{k'}\cong\M_d(k')$.

  Let $\alpha$ be an automorphism of the algebra $A$.
  By the Noether-Skolem Theorem, every automorphism of $A$ is inner
  (see~\cite{pierce}, \S~12.6), therefore 
  there is $c\in A$ such that $\alpha\co x\mapsto c^{-1}xc$.
  Since the image of a left ideal in $A_{k'}$ under $\alpha$ is again 
  a left ideal, 
  $\alpha$ permutes the set $\{\mathcal{L}\}$ 
  of all $d$-dimensional left ideals in $A_{k'}$.
  Now the situation is illustrated by the following diagram:
  \begin{displaymath}
    \begin{diagram}
      \{\mathcal{L}\} & \rTo^{\alpha}           & \{\mathcal{L}\} \\
      \dTo^{P}        &                            & \dTo_{P} \\
      \mathcal{V}_A   & \rDashto^{\tilde\alpha} & \mathcal{V}_A
    \end{diagram}
  \end{displaymath}
  where $P$ stands for Pl\"ucker embedding and $\tilde\alpha$ is a
  map of $\mathcal{V}_A$ to itself such that the diagram is commutative.
  From the commutativity then follows that 
  $\alpha\mapsto\tilde\alpha$ is a homomorphism of groups.
  Since $P$ is a bijection between $\{\mathcal{L}\}$ and $\mathcal{V}_A$,
  the map $\alpha\mapsto\tilde\alpha$ is injective.
  We claim that $\tilde\alpha$ is the restriction of a linear
  transformation of the projective space to $\mathcal{V}_A$
  and that this transformation is defined over the field $k$.

  Let us fix a $k$-basis in $A$.
  Let $b_1,\dots, b_d$ be vectors spanning $\mathcal{L}$ such that 
  their coordinates are with respect to the fixed $k$-basis of $A$.
  The Pl\"ucker embedding 
  is taking all $d\times d$ minors of the $d^2\times d$ matrix $b$
  containing the coordinates of the vectors $b_i$ as columns
  (compare with Example~\ref{exa:quaternions}).

  The automorphism $\alpha$ maps the left ideal $\mathcal{L}$
  to $\mathcal{L}c$, since for the invertible $c\in A$ we have that
  $c^{-1}\mathcal{L} = \mathcal{L}$.
  Let $C$ be the $d^2\times d^2$ matrix of the linear transformation
  $x\mapsto xc$ for all $x\in A$, i.e. if $v_x$ is the vector of
  coordinates of $x$ then $C v_x$ is the vector of coordinates of $xc$.
  Clearly, $C$ is a matrix over $k$.
  The Pl\"ucker coordinates of $\alpha(\mathcal{L})$ are all
  $d\times d$ minors of the matrix $Cb$.

  In following for an $m\times n$ matrix $M$ and subsets 
  $I\subset\{1,\dots,m\}$, $J\subset\{1,\dots,n\}$ such that
  $|I| = |J| = d$ we by $|M|_{I,J}$ mean the minor of the matrix $M$
  which is the determinant of the $d\times d$ matrix obtained from $M$
  by omitting all rows not in $I$ and all columns not in $J$.

  The Pl\"ucker coordinates of $\mathcal{L}$ are $p_I = |b|_{I,\{1,\dots,d\}}$ 
  ($I\subset\{1,\dots,d^2\}$, $|I|=d$).
  For the Pl\"ucker coordinates $p'_I$ of $\alpha(\mathcal{L})$ we then have
  \begin{displaymath}
    p'_I = |Cb|_{I,\{1,\dots,n\}} = 
      \sum_{I'} |C|_{I,I'} |b|_{I',\{1,\dots,n\}} =
      \sum_{I'} |C|_{I,I'} p'_{I'},
  \end{displaymath}
  where the sum is taken through all subsets $I'$ of $\{1,\dots,d^2\}$ of $d$ elements.
  Since $C$ is a matrix over $k$, also all its $d\times d$ minors are in $k$
  and the claim is proven.

  We have shown that for a $k$-automorphism $\alpha$ 
  is $\tilde\alpha$ a $k$-automorphism of $\mathcal{V}_A$.
  Moreover, the group homomorphism 
  $\Aut_k A\to\Aut_k\mathcal{V}_A$, $\alpha\mapsto\tilde\alpha$ 
  is an algebraic map between twists of $\PGL_d(k)$.
  We have seen that it is injective and therefore from being algebraic 
  it follows that it is also surjective.
  For, $\PGL_d(k)$ is an irreducible variety, so if the
  morphism of groups were not surjective, the image of $\Aut A$
  would be a subvariety in $\Aut\mathcal{V}_A$ of positive codimension,
  hence $\alpha\mapsto\tilde\alpha$ would not be injective, a contradiction.  
\end{proof}

\begin{lem}\label{le:csaAutLie}
  If $A$ is a central simple algebra, then for the Lie algebra 
  $L(\Aut A)$ of the group $\Aut A$ holds
  $L(\Aut A)\oplus k\cong A_{Lie}$.
\end{lem}

\begin{proof}
  Since each automorphism of $A$ is inner, $\Aut A$ is a projectivization
  of $A^*$, the group of all invertible elements in $A$.
  An element $a\in A$ is not invertible, if the matrix of the 
  regular representation of $a$ is singular,
  therefore $A^*$ is a principal open subset of the linear variety $A$.
  The Lie algebra $L(A^*)$ is the tangent space to $A^*$ at $e$, therefore
  it consists of elements of $A$ and the multiplication in $L(A^*)$
  is $[a,b] = ab-ba$, so $L(A^*) = A_\mathrm{Lie}$.
  The claim of the Lemma now follows from $\Aut A = A^*/\mathbf{C}(A^*)$.
\end{proof}

Summing up, we have the following:

\begin{prop}
  Let $A$ be a twist of $\M_d(k)$ and $\mathcal{V}_A$ the corresponding
  Severi-Brauer variety. 
  Further let $A'$ be the associative algebra constructed in our algorithm, 
  (i.e.~such that $\g_0(\mathcal{V}_A,k)$ is embedded into $A'_\mathrm{Lie}$).
  \begin{enumerate}
    \item[(i)]  If $d=2$ then $A'\cong A$.
    \item[(ii)] If $d>2$ then either $A'\cong A$ or $A'\cong A^{op}$, 
      the opposite algebra to $A$.
  \end{enumerate}
\end{prop}

\begin{proof}
  First, by Proposition~\ref{prop:autIso} there is a $k$-isomorphism
  of algebraic groups $\Aut(A)\to\Aut(\mathcal{V}_A)$.
  Its differential gives a $k$-isomorphism of their Lie algebras
  \begin{displaymath}
    L(\Aut(A))\to L(\Aut(\mathcal{V}_A)) = {\g}_0(\mathcal{V}_A,k).
  \end{displaymath}
  Suppose first that $d=2$.
  Then by Proposition~\ref{prop:uniqueEnvelop} there is a unique
  enveloping algebra $A'$ of $\g_0(\mathcal{V}_A,k)$ such that $[A':k] = 4$.
  (i.e.~an algebra $A'$ with an embedding 
  $\g_0(\mathcal{V}_A,k)\hookrightarrow (A')_\mathrm{Lie}$)
  and our algorithm finds it (see Lemma~\ref{le:findEnvelop}).
  On the other hand $A$ is by Lemma~\ref{le:csaAutLie} 
  an enveloping algebra of $L(\Aut(A))$. So again by uniqueness
  in degree $2$ case we conclude that $A'\cong A$.

  Now if $d>2$ then again by Lemma~\ref{le:findEnvelop} the enveloping
  algebra found by the algorithm has dimension $d^2$.
  But Proposition~\ref{prop:uniqueEnvelop} together with Lemma~\ref{le:csaAutLie} 
  in this case give that $A$ and $A'$ are either isomorphic or antiisomorphic.
\end{proof}

\chapter{Trivializing central simple algebras}\label{ch:csa}

Assume that $\g$ is a simple Lie algebra which is
a twist of $\mf{sl}_n(k)$, i.e.~$\g\otimes\bar k\cong\mf{sl}_n(\bar k)$. 
We want to decide whether $\g$ is isomorphic to $\mf{sl}_n$ over $k$
and if so, to find an isomorphism.

Because of the algorithm {\tt EnvelopingAlgebra} in 
Chapter~\ref{ch:preliminaries}, Section~\ref{se:LieAssoc} 
this is equivalent to deciding whether a given associative algebra
$A$ is isomorphic to $\M_n(k)$ (see Corollary~\ref{cor:LieAsocIso})
and finding an isomorphism in the affirmative case.
The short algorithm for reducing the problem for Lie algebras
to a corresponding one for associative algebras could read as follows:

\vskip0.4cm
\noindent
\begin{tabular}{ll}
  {\sc Algorithm:} & {\tt TrivializeLieAlgebra}\\
  {\sc Input:}     & $\g$ -- a simple Lie algebra which is a twist of $\mf{sl}_n(k)$.\\
  {\sc Output:}    & an isomorphism $\mf{sl}_n(k)\to\g$, if it exists.
\end{tabular}
\begin{enumerate}
  \item $A, \phi\co\g\hookrightarrow A_{\textrm{Lie}}$ := {\tt EnvelopingAlgebra}($\g$);
  \item {\tt if} $[A:k]>n^2$ {\tt then} \\
    \phantom{xx} {\tt return 'failed'}\\
    {\tt end if};
  \item Let $\varphi\co\M_n(k)\to A$ be an isomorphism of associative algebras;\\
    {\it // now $\varphi$ is also a Lie algebra isomorphism 
    $\mf{gl}_n(k)\to A_\mathrm{Lie}$}
  \item {\tt return} restriction of $\varphi$ to $\mf{sl}_n(k)$.
\end{enumerate}
\vskip0.4cm

In this chapter we describe the step (3) of the algorithm:
for a given central simple algebra of the degree up to $4$
we will construct an isomorphism with $\M_n(k)$.

\section{Cyclic algebras and left ideals}\label{se:csa}

In this section the role of one-sided ideals in a central simple
algebra is explored. We give some hints for finding such ideal
by reducing the problem to solving a norm equation.
Knowing the solution of the norm equation will make it possible
to construct an isomorphism with $M_n(k)$.
There are already efficient algorithms for deciding whether 
a given central simple algebra is a full matrix algebra 
(see for example~\cite{gl}), therefore assuming that there exists
an isomorphism does not impose much restriction. 
Nevertheless, the decision algorithms, like the mentioned one,
do not give an explicit isomorphism in case it exists.

Assume $A\cong\M_n(k)$. 
Then $A$ contains one-sided ideals of dimension $n$.
We will work here with left ideals. 
Every minimal left ideal in this algebra has dimension $n$
and any nontrivial left ideal is a direct sum of minimal left ideals.

Let $\mathcal{L}$ be an $n$-dimensional left ideal in $A$.
For any $a\in A$ we have that $\varphi_a\co\mathcal{L}\to\mathcal{L}$, 
$x\mapsto ax$ is an endomorphism of $\mathcal{L}$ as a vector space.
Let us fix a basis $b_1,\dots,b_n$ of $\mathcal{L}$.
Let $\varphi\co A\to\M_n(k)$ assign to $a\in A$ the matrix of $\varphi_a$ 
with respect to this basis, i.e. the $i$-th column of $\varphi(a)$
contains the coordinates of $ab_i$: $ab_i = \sum_j \varphi(a)_{ji}b_j$.

\begin{prop}\label{pro:leftId}
  $\varphi$ is an isomorphism of algebras.
\end{prop}

\begin{proof}
  Firstly, $\varphi$ is a homomorphism of algebras. 
  Further $\varphi$ is a bijection, since otherwise $\Ker\varphi\neq 0$ would
  be a nontrivial ideal in $A$.
\end{proof}

To make use of the proposition, we have to find an $n$-dimensional left ideal
in the algebra $A$.
It might help, if during our computation we hit a zero divisor $d$.
When for such $d$ we define the vector space endomorphism $\rho_d$ of $A$,
$x\mapsto xd$, 
then both the kernel and the image of $\rho_d$ are nontrivial left ideals of $A$.
If we are lucky enough, or if the algebra $A$ is of degree $2$ or $3$,
then we have found a left ideal of dimension $n$.

In following, a finite extension $k'$ of the field $k$ is said to be {\em cyclic},
if $k' = k(a)$ for some algebraic element $a$ and the minimum polynomial of $a$ factors
completely over $k(a)$. Such field is Galois and $\Gal(k(a)|k)$ is cyclic.

In our algorithm for finding left ideals, the crucial is the following
\begin{df}\label{def:cyclic}
  The algebra $A$ is a {\em cyclic} of degree $n$, if there are
  elements $c,u \in A$ such that
  $\{c^i u^j \mid i,j = 0,\dots,n-1\}$ is a basis of $A$ over its ground field $k$,
  and the multiplication in $A$ is defined by:
  \begin{itemize}
    \item[(i)]   $1,c,\dots,c^{n-1}$ is a basis of a cyclic field extension 
      $k'=k(c)$ of $k$,
    \item[(ii)]  $uc = \sigma(c)u$, where $\sigma$ is a fixed generator of 
      the Galois group $\Gal(k'|k)$,
    \item[(iii)] $u^n = \gamma 1$, where $0\ne\gamma\in k$.
  \end{itemize}
  We will call $c$ a {\em cyclic element} and
  $u$ a {\em principal generator} of $A$ over $k(c)$.
\end{df}

\begin{lem}
  A cyclic algebra is central simple.
  If $c\in A$ is the cyclic element as in Definition~\ref{def:cyclic},
  then the centralizer $\mathbf{C}_A(k(c))=k(c)$, 
  so it is a maximal subfield of $A$.
\end{lem}

\begin{proof}
  cf.~\cite{jacBA2}, Theorem~8.7.
\end{proof}

The main strategy in finding a minimal left ideal is 
writing a given algebra $A$ as a cyclic algebra, which means 
to find a cyclic element $c$ and a principal generator $u$ of $A$ over $k(c)$.
The reason for doing so is given by 

\begin{prop}\label{the:normeq}
  Let $A$ be a cyclic algebra generated by $c$ and $u$ as in
  Definition~\ref{def:cyclic}. Then $A\cong\M_n(k)$ if and only if 
  there exists $s\in k(c)\subset A$ such that
  \begin{eqnarray}\label{eqn:norm}
    N_{k(c)|k}(s) = s\sigma(s)\dots\sigma^{n-1}(s) = \frac{1}{\gamma}.
  \end{eqnarray}
\end{prop}

\begin{proof}
  We prove the proposition for degree $n=2$, the other cases
  can be proven analogously, the proof is just a bit more technical.
  A proof for any central simple algebra can also be found 
  for example in~\cite{pierce}, \S~15.1.

  Since $A$ is central simple, 
  by Wedderburn's structure theorem it can be
  either $\M_2(k)$ or a division algebra over $k$. It is the matrix
  algebra exactly if it has a 2-dimensional left ideal.
  Let $\mathcal{L}$ be such an ideal and let 
  $0\ne c_0 1 + c_1 u \in \mathcal{L}$, $c_i \in k(c)$.
  Then $c_0$ is nonzero, for otherwise $u\cdot c_1 u$
  would be an invertible element in $\mathcal{L}$. 
  Therefore we may suppose $c_0=1$
  and then $\mathcal{L}$ is spanned by $1 + c_1 u$ and $c(1 + c_1 u)$.
  From $u(1 + c_1 u) = \sigma(c_1)\gamma + u \in\mathcal{L}$ 
  we get that $c_1$ is a solution to the norm equation~(\ref{eqn:norm}).
\end{proof}

The proof of the previous proposition delivers also an algorithm 
for finding an isomorphism of a cyclic algebra $A$ and $\M_n(k)$: 
If the norm equation~(\ref{eqn:norm}) is solvable, then
using a solution we can construct an $n$-dimensional left ideal.
The remaining problem is finding a cyclic element in $A$.
In the following sections we discuss this for degrees $2$ and $3$
and explain how to use the developed methods in case of algebras 
of degree $4$.

\section{Algebras of degree 2}\label{subse:M2}

This is a classical method which can be found for example in~\cite{tits}.
We shortly describe it here for the sake of completeness.

When looking for a cyclic element $c$ we pick an arbitrary 
noncentral element in $A$. Its minimum polynomial $\mu_c(\xi)$ 
is quadratic. In case it is reducible over $k$, 
$\mu_c(\xi) = p_1(\xi)p_2(\xi)$, 
we have also found a zero divisor, for example $p_1(c)$, and are therefore done.
Otherwise $c$ generates in $A$ a quadratic field extension of $k$.
By factoring $\mu_c(\xi)$ over $k'=k(c)$ we find $\sigma(c)$.
Then we find a principal generator $u$ by solving
a linear system of equations $uc = \sigma(c)u$. 
There exists such invertible $u$, because the matrices $c$ and $\sigma(c)$
have the same invariant factors and therefore are similar, cf.~\cite{wedd}.
Finally, since $c$ and $u$ are generators of $A$ 
and $u^2c = cu^2$, $u^2\in\mathbf{C}(A) \cong k$.

In case the field $k=\Q$, there are already much more effective algorithms 
for finding an isomorphism $\M_2(\Q)\to A$, see~\cite{cremona, gabor, simon2}.
Since there is an efficient algorithm for finding a rational point on 
a plane conic over $\Q$ implemented in Magma, in our algorithms we reduce 
the degree 2 case to this problem, as described in detail 
in Chapter~\ref{ch:aut}, Section~\ref{se:conics}.

\section{Algebras of degree 3}

In this case finding a cyclic element in the algebra is more technical.
Throughout the whole section we assume that the given algebra $A$
is isomorphic to $\M_3(k)$.

\begin{lem}\label{lem:conj}
  Let $a\in A$ be a noncentral element such that
  the minimal polynomial $\mu_a(\xi)\in k[\xi]$ of $a$
  is irreducible over $k$.
  Then $\deg \mu_a(\xi) = 3$ and
  every $b\in A$, $b\neq 0$ such that $\mu_a(b)=0$ is a conjugate of $a$.
\end{lem}

\begin{proof}
  $a\notin\mathbf{C}(A)$ implies $\deg\mu_a(\xi) > 1$.
  The case $\deg\mu_a(\xi) = 2$ is not possible. 
  For, let $\deg\mu_a(\xi) = 2$, $\mu_a(\xi)$ irreducible. 
  Then the characteristic polynomial $\chi_a(\xi)$ of $a$ is
  $\chi_a(\xi) = \mu_a(\xi)l(\xi)$ with $l(\xi)$ linear.
  The factor $l(\xi)$ of the characteristic polynomial of $a$ 
  is irreducible, therefore it divides the minimal polynomial
  $\mu_a(\xi)$, a contradiction.
  Since $\mu_a(\xi)$ is irreducible over $k$ and $\mu_a(b)=0$, 
  it is also the minimal polynomial of $b$. 
  Then $a$ and $b$ have the same invariant factors 
  and hence are conjugate (cf.~\cite{wedd}).
\end{proof}

In a given algebra of degree $A$ we will first try to construct an 
element $u\in A$ such that $u^3\in\mathbf{C}(A)$, and afterwards
a cyclic element $c\in A$ such that $u$ is a principal generator
of $A$ over $k(c)$.

We start by picking an arbitrary noncentral element $x\in A$. 
If $x$ is not invertible or its minimal polynomial $\mu_x(\xi)$ is 
reducible over $k$, then we are done, since we have found a zero divisor.
So we may assume that $\mu_x(\xi)$ is irreducible. 
Then by Lemma~\ref{lem:conj} we have $\deg\mu_x(\xi) = 3$. 

In the special case when $x$ is already cyclic, 
we set $c=x$. Then $k'=\left<1,c,c^2\right>_k$ 
is a maximal subfield generated by $c$. 
Let $\sigma$ denote a generator of $\Gal(k'|k)$.
By factoring the minimal polynomial $\mu_c(\xi)$ over $k'$ 
we find $\sigma(c)$ and afterwards an element $u\in A$ 
such that $ucu^{-1} = \sigma(c)$ as a nontrivial solution to 
the linear system $uc = \sigma(c)u$.
By Lemma~\ref{lem:conj} such $u$ exists
and since it does not commute with $c$, we have $u\notin k'$.
From $u^3cu^{-3}=c$ follows $u^3\in \mathbf{C}_A(k')=k'$.
We can conclude that $u^3\in k$ since otherwise $u^3$ would generate $k'$
and $u^3u = uu^3$ would imply $u\in\mathbf{C}_A(k')=k'$, a contradiction.
So in the special case the chosen element $x\in A$ is cyclic,
we can easily find a principal generator.
Therefore in our construction we may assume that $x\in A$ is non-cyclic.

By $A[\xi]$ we denote the ring of polynomials in $\xi$ over 
the algebra $A$, where $\xi$ commutes with all elements in $A$. 
As usual we say that $a\in A$ is a root of 
$p(\xi) = c_0 + c_1\xi + \dots + c_n\xi^n$,
if $c_0 + c_1 a + \dots + c_n a^n = 0$.
If $a$ is a root of $p(\xi)\in A[\xi]$, then $p(\xi) = q(\xi)(\xi-a)$
for some $q(\xi)\in A[\xi]$, see~\cite{wedd}.

Recall that by $[a,b]$ we mean $ab-ba$.

\begin{thm}[Wedderburn's factorization theorem]
  Let $a\in A$ have the minimum polynomial $\mu_a(\xi)$ 
  and $m\in A[\xi]$ be such that $\mu_a(\xi) = m(\xi)(\xi-a)$. 
  Then for any $y\in A$ such that $[y,a]$ is invertible, 
  $a^\prime = [y,a]a[y,a]^{-1}$ is a root of $m(\xi)$.
\end{thm}

\begin{proof}
  We partly follow~\cite{jacFDDA}, Theorem~2.9.1, where the claim 
  is proven for division algebras.

  Since $\mu_a$ is a polynomial over the field, we have 
  $\mu_a(\xi) y = y \mu_a(\xi)$ and hence
  $m(\xi)(y\xi - ay) = y m(\xi)(\xi-a)$.
  Here the left hand side can be written as 
  $m(\xi)(y\xi - ay) = m(\xi)y(\xi-a) + m(\xi)[y,a]$, 
  therefore together it gives 
  $m(\xi)[y,a] = (y m(\xi) - m(\xi) y)(\xi-a)$.
  After multiplying by the inverse of $[y,a]$ we get
  $m(\xi) = (y m(\xi) - m(\xi) y)[y,a]^{-1}(\xi - [y,a]a[y,a]^{-1})$.
\end{proof}

Using Wedderburn's factorization theorem we can over $A$ factorize 
the minimal polynomial $\mu_a$ of the chosen element $a\in A$ 
into linear factors
\begin{equation}\label{eq:factors}
  \mu_a(\xi) = (\xi-a_3)(\xi-a_2)(\xi-a_1)
\end{equation}
with $a_i\in A, a_1=a$ and $a_2 = [y,a]a[y,a]^{-1}$ for some $y\in A$,
unless we hit a zero divisor.
Any factorization obtained in this way has useful properties:

\begin{lem}\label{lem:factors}
  For $a\in A$ noncyclic with irreducible minimal polynomial,
  the factorization~(\ref{eq:factors}) fulfills the following:
  \begin{itemize}
    \item[(i)]   $c = [a_1,a_2] = [a_2,a_3] = [a_3,a_1] \neq 0$,
    \item[(ii)]  $c$ is invertible and $ca_ic^{-1} = a_{i+1} \quad$ 
                 (indices reduced mod 3),
    \item[(iii)] $c^3 = \gamma,\quad \gamma\in k$.
  \end{itemize}
\end{lem}

The lemma can be found in~\cite{jacFDDA} as Lemma~2.9.8 for division
algebras. Here we give an alternative proof since it is slightly
different for the matrix algebra $\M_3(k)$.
First we have 2 small technical lemmata

\begin{lem}\label{lem:switch}
  Let $f,g\in A[\xi]$ be monic polynomials such that $fg\in k[\xi]$.
  Then $fg = gf$.
\end{lem}

\begin{proof}
  Let us denote $p=fg$ and $q=gf$.
  Using the facts that $p\in k[\xi]$ and $f$ is monic, 
  after matching coefficients in $pf = fq$ we get that $p$ = $q$.
\end{proof}

\begin{lem}\label{lem:isCyclic}
  Let the minimum polynomial of $a\in A$ be irreducible of degree 3.
  If there is $y\in A$ such that $yay^{-1}\ne a$ and $[yay^{-1},a] = 0$,
  then $a$ is cyclic.
\end{lem}

\begin{proof}
  We follow~\cite{jacFDDA}.
  Since $k(a)$ is a maximal subfield of $A=\M_3(k)$, $[yay^{-1},a] = 0$ 
  implies that $yay^{-1}\in k(a)$. 
  It follows that conjugation by $y$ is a field automorphism of $k(a)$
  and from $yay^{-1}\ne a$ we have it is a nontrivial automorphism.
  Therefore $k(a)$ is Galois over $k$.
\end{proof}

\begin{proof}[Proof of Lemma~\ref{lem:factors}]
  First we show that $[a_1,a_2] \neq 0$.
  For, let $[a_1,a_2]=0$ whereas $a_2=[y,a_1]a_1[y,a_1]^{-1} \neq a_1$.
  But then by Lemma~\ref{lem:isCyclic} is $a_1$ cyclic, a contradiction.

  Next, since $\mu_a(\xi)$ is a polynomial over the field $k$, 
  the factors can be permuted cyclically (see Lemma~\ref{lem:switch}). 
  Then from
  \begin{align*}
    \mu_a(\xi) = \xi^3 + \alpha_2\xi^2 + \alpha_1\xi + \alpha_0
    &= (\xi - a_3)(\xi - a_2)(\xi - a_1) \\
    &= (\xi - a_1)(\xi - a_3)(\xi - a_2) \\
    &= (\xi - a_2)(\xi - a_1)(\xi - a_3) 
  \end{align*}
  we get
  \begin{align*}
    \alpha_1 = a_3a_2 + a_3a_1 + a_2a_1 
             = a_1a_3 + a_1a_2 + a_3a_2
             = a_2a_1 + a_2a_3 + a_1a_3
  \end{align*}
  and (i) follows.

  To prove (ii) and (iii) we first observe, that 
  $a_3a_2a_1 = a_2a_1a_3 = a_1a_3a_2 = -\alpha_0$.
  Using this and the part (i) of the lemma we have
  \begin{displaymath}
    [a_3,a_1]a_1 = [a_2,a_3]a_1 = a_2a_3a_1 - a_3a_2a_1 = 
    a_2a_3a_1 - a_2a_1a_3 = a_2[a_3,a_1].
  \end{displaymath}
  Similarly we deduce $[a_1,a_2]a_2 = a_3[a_1,a_2]$
  and $[a_2,a_3]a_3 = a_1[a_2,a_1]$. Then $c^3a_i = a_ic^3$ for $i=1,2,3$.
  Since $a_i$'s generate the whole $A$, 
  $c^3$ is in the center of $A$ and (iii) follows. 
  To finish the proof, it remains to show that $c$ is invertible. 
  For, let $c^3 = 0$. Let $\mathcal{L} = Ac$ be a left ideal. 
  Since $c\ne 0$, also $\mathcal{L}\ne 0$.
  We observe that $(a_ic)^3 = a_ia_{i+1}a_{i+2}c^3 = 0$ (indices reduced mod 3)
  and conclude that $\mathcal{L}$ is nilpotent (since $a_i$'s generate $A$) 
  and hence contained in a nilpotent two-sided ideal (cf.~\cite{jacToR},
  \S~4.10, Lemma~2), a contradiction, since $A$ is simple and not nilpotent.
\end{proof}

Lemma~\ref{lem:factors} gives an algorithm for finding an element $c\in A$
such that $c^3$ is in the center of $A$. 

\vskip0.4cm
\noindent
\begin{tabular}{lp{11.9cm}}
  {\sc Algorithm:} & {\tt CubicRoot}\\
  {\sc Input:}     & $a$ -- an element in $A$ such that the minimal polynomial 
                       $\mu_a$ is cubic and irreducible over $k$.\\
  {\sc Output:}    & $c$ -- a noncentral element in $A$ such that\\
                   & \hangafter 1 \hangindent 19pt $\quad\bullet$ 
                   $c^3\in k$ and \\
                   & \hangafter 1 \hangindent 19pt $\quad\bullet$ 
                   $(\xi - a')(\xi - c a c^{-1})(\xi - a)$ is a factorization 
                   of $\mu_a$ as in Lemma~\ref{lem:factors}.
\end{tabular}
\begin{enumerate}
  \item $y$ := an element in $A$ such that $[y,a]$ is invertible
    and not commuting with $a$;
  \item $a'$ := $[y,a]a[y,a]^{-1}$;
  \item {\tt return} $[a,a']$.
\end{enumerate}
\vskip0.4cm

\begin{lem}\label{lem:findCyclic}
  Let $b\in A$, $b\notin k$ such that $b^3 = \beta$, $\beta\in k^*$. Let
  $\mu_b(\xi) = \xi^3 - \beta = (\xi - b_3)(\xi - b_2)(\xi - b_1)$, $b_1 = b$
  be the factorization of its minimal polynomial as in Lemma~\ref{lem:factors}.
  Let $c = [b_1,b_2]$. 
  Then $a = b_1c$ is a cyclic element not in $k$.
\end{lem}

\begin{proof}
  The proof follows word by word the construction for division algebras
  which can be found in~\cite{jacFDDA}, in the proof of Theorem~2.9.17.
\end{proof}

\vskip0.4cm
\noindent
\begin{tabular}{lp{11.9cm}}
  {\sc Algorithm:} & {\tt GeneratorsOfCyclicAlgebra}\\
  {\sc Input:}     & $A$ -- a central simple algebra. \\
  {\sc Output:}    & $c$ -- a cyclic element in $A$,\\
                   & $u$ -- a principal generator of $A$ over $k(c)$.
\end{tabular}
\begin{enumerate}
  \item $x$ := a noncentral element in $A$;
  \item {\tt if} $a$ is cyclic {\tt then}\\
    \phantom{xx} $c$ := $x$;\\
    \phantom{xx} $c^\sigma$ := any other root of $\mu_c$ in $k(c)$;\\
    \phantom{xx} $u$ := a nontrivial solution of a linear system
      $c^\sigma u = u c$;\\
    \phantom{xx} {\tt return} $c,u$\\
    {\tt end if};
  \item {\tt if} $x^3\in k$ {\tt then} $u' := x$;\\
    {\tt else} $u'$ := {\tt CubicRoot}($x$);\\
    {\tt end if};
  \item $u$ := {\tt CubicRoot}($u'$);
  \item {\tt return} $u'u$ and $u$.
\end{enumerate}
\vskip0.3cm

To complete the section we add the algorithm which constructs 
a 3-dimensional left ideal in $A$ in case it exists.

\vskip0.4cm
\noindent
\begin{tabular}{lp{11.9cm}}
  {\sc Algorithm:} & {\tt FindMinimalLeftIdeal}\\
  {\sc Input:}     & $A$ -- a central simple algebra. \\
  {\sc Output:}    & $\mathcal{L}$ -- a three dimensional left ideal in $A$.
\end{tabular}
\begin{enumerate}
  \item $c,u$ := {\tt GeneratorsOfCyclicAlgebra}($A$);
  \item $\gamma$ := $u^3$; {\it // $\gamma\in k$}
  \item $k'$ := $k(c)$;
  \item {\it //  solve the norm equation}\\
    Let $s$ be such that $N_{k'|k}(s) = 1/\gamma$;
  \item $b$ := $1 + su + s^\sigma u^2$, where $s^\sigma = u s u^{-1}$;
  \item Let $\mathcal{L}$ be the $3$-dimensional left ideal generated by $b$;
  \item {\tt return} $\mathcal{L}$.
\end{enumerate}
\vskip0.3cm

For each element $a\in A$ which pops up during the computation it is tested 
whether $a$ is a zero divisor or whether the minimal polynomial $\mu_a$ is reducible.
In both cases we would have found a zero divisor $d$. Then either the kernel
or the image of $x\mapsto dx$ is already a 3-dimensional left ideal.
If this happens, we can skip the rest of computation, in particular we can avoid 
time-expensive solving of the norm equation.

\section{Algebras of degree 4}

The main goal in the algorithm for finding an isomorphism of a given algebra 
$A$ and $\M_4(k)$ is finding a zero divisor.
For any $d\in A$ we by $\rho_d$ denote the vector space endomorphism of $A$,
$x\mapsto xd$. 
Both the kernel and the image of $\rho_d$ are left ideals in $A$.
If $d\in A$ is a zero divisor, $\Ker\rho_d$ is clearly nontrivial, 
and the same with $\im\rho_d$, since $1.d\neq 0$.
If $\dim\Ker\rho_d = 4$ or $\dim\im\rho_d = 4$ or 
$\dim(\Ker\rho_d\cap\im\rho_d) = 4$, we are done since we can already
use Proposition~\ref{pro:leftId} and find an isomorphism $A\to M_4(k)$.
Here we first describe how to use other kinds of zero divisors
and afterwards how to find one at all.

\begin{lem}\label{le:zd}
  Let $\varphi\co A\to M_4(k)$ be an isomorphism and let
  $d\in A$ be a zero divisor such that none of 
  $\ \Ker\rho_d$, $\im\rho_d$, $\Ker\rho_d\cap\im\rho_d$
  has dimension~$4$.
  Then $\varphi(d)$ is similar to one of the following block matrices:
  \begin{enumerate}
    \item\label{it:z1}
      $\left(\begin{array}{cc}
        D & 0 \\ 0 & 0
      \end{array}\right),\quad$
      where $D\in M_2(k)$ is an invertible matrix,
    \item\label{it:z2}
      $\left(\begin{array}{cc}
        B & 0 \\ 0 & B
      \end{array}\right),\quad$
      where $B = \left(\begin{array}{cc}
        0 & 1 \\ 0 & 0
      \end{array}\right)$.
  \end{enumerate}
\end{lem}

\begin{proof}
  We consider the Jordan normal form of $\varphi(d)$. 
  Since $\dim\Ker\rho_d=8$, we can conclude that $\varphi(d)$ 
  has at most two nonzero eigenvalues. 
  Therefore the minimum polynomial $\mu_d(\xi)$ is reducible
  and divisible by $\xi^2$.
  If $\mu_d(\xi) = \xi^4$, then we get the case~(\ref{it:z2}) for the 
  Jordan normal form of $\varphi(d)$.
  If $\xi$ appears in $\mu_d(\xi)$ in degree $3$ 
  (i.e.~$\varphi(d)$ has exactly one nonzero eigenvalue),
  then at least one of the three left ideals 
  $\im\rho_d$, $\Ker\rho_d$, $\im\rho_d\cap\Ker\rho_d$ is four-dimensional.
  Lastly, if $\xi$ appears in $\mu_d(\xi)$ exactly in degree $2$,
  we get that the case~(\ref{it:z1}).
\end{proof}

If $\Ker\rho_d \cap \im\rho_d = 0$, then $d$ is a zero divisor of 
type~(\ref{it:z1}) in Lemma~\ref{le:zd}.
We define another vector space endomorphism $\lambda_d$ of $A$,
$x\mapsto dx$. The intersection $A_1 = \im\rho_d\cap\im\lambda_d$
is mapped by $\varphi$ to the subalgebra of all block matrices,
where only the upper left $2\times 2$ block in nonzero, 
so $A_1\cong M_2(k)$. We find a zero divisor $d_1$ in $A_1$ 
as mentioned in Section~\ref{subse:M2}. Then $\im\rho_{d_1}$
is a 4-dimensional left ideal in $A$.

The second case is a bit more tricky. 
Let $d$ be a zero divisor of type~(\ref{it:z2}).
We denote by $A_d$ the centralizer $\C_A(d)$ and by $\mathfrak{R}(A_d)$
the Jacobson radical of $A_d$.
Then there is the natural projection $\pi\co A_d \to A_d/\mathfrak{R}(A_d)$
and for this we have
\begin{lem}
  The algebra $\pi(A_d)$ is isomorphic to $M_2(k)$.
  If $e\in\pi(A_2)$ is a zero divisor, then for a generic element
  $f$ in $\pi^{-1}(e)$ we have $\dim\Ker\rho_f = 4$.
\end{lem}

\begin{rem}
  By saying that something holds for a ``generic element''
  we mean that all elements, for which the assertion is true, 
  form a nonempty Zariski open subset in the set of all considered
  elements. Therefore we can easily find an element 
  satisfying the condition.
\end{rem}

\begin{proof}
  We may suppose that
  $\varphi(d)$ is actually equal to the matrix~(\ref{it:z2}) in Lemma~\ref{le:zd}.
  Then the image of $A_d = C_A(d)$ under $\varphi$ is
  \begin{equation}\label{eq:A_d}
    \varphi(A_d) = \left\{\left(\begin{array}{cccc}
      \alpha_1 & \beta_1  & \alpha_2 & \beta_2  \\
      0        & \alpha_1 & 0        & \alpha_2 \\
      \alpha_3 & \beta_3  & \alpha_4 & \beta_4  \\
      0        & \alpha_3 & 0        & \alpha_4
    \end{array}\right)\ \Bigg|\ \alpha_i,\beta_i\in k\right\}.
  \end{equation}
  The Jacobson radical in $\mathfrak{R}(\varphi(A_d))$ is the set of all 
  $a\in\varphi(A_2)$ such that $\alpha_i=0$, ($i=1,\dots,4$),
  and $\varphi(A_d)/\mathfrak{R}(\varphi(A_d))\cong
  \{a\in\varphi(A_2)\mid \beta_i=0,i=1,\dots 4\}$,
  which is clearly isomorphic to $M_2(k)$.
  Let us denote the natural projection 
  $\varphi(A_d)\to\varphi(A_d)/\mathfrak{R}(\varphi(A_d))$ by $\pi^\prime$.
  If $e$ is a zero divisor in $M_2(k)$, 
  then the preimage $(\pi^\prime)^{-1}(e)$ consists
  of such elements in~(\ref{eq:A_d}), that $\alpha_i$'s are fixed
  and $\alpha_1\alpha_4 = \alpha_2\alpha_3$.

  Now for a matrix $f'$ in $(\pi')^{-1}(e)$ we have 
  that $\dim\Ker\rho_{f'} = 4$ if and only if
  the rank of $f'$ equals $3$.
  The last is equivalent to 
  $\alpha_1\beta_4+\alpha_4\beta_1\ne\alpha_2\beta_3 + \alpha_3\beta_2$.
  Indeed, the rank of $f'$ equals $3$ exactly if there is a $3\times 3$
  nonzero minor of $f'$. Such minor can be found as the determinant of 
  a submatrix containing all $\beta_i$'s.
  Hence for a generic $f'\in(\pi')^{-1}(e)$ we have $\dim\Ker\rho_{f'} = 4$, 
  so $\varphi^{-1}(f')$ gives us a minimal left ideal in $A$.
\end{proof}

From the proof of the previous Lemma it also follows that we can 
compute the the Jacobson radical of $A$ very easily, namely 
$\mathfrak{R}(A_d) = \Ker\rho_d \cap \Ker\lambda_d$.

\newpage
\noindent
\begin{tabular}{lp{11.9cm}}
  {\sc Function:} & {\tt FindMinimalLeftIdeal}\\
  {\sc Input:}    & $A$ -- an associative algebra isomorphic to $M_4(k)$,\\
                  & $d$ -- a zero divisor in $A$.\\
  {\sc Output:}   & $\mathcal{L}$ -- a four-dimensional left ideal in $A$.
\end{tabular}
\begin{enumerate}
  \item {\tt if} $\dim\Ker\rho_d = 4$ {\tt then return} $\Ker\rho_d$; \\
    {\tt elif} $\dim\im\rho_d = 4$ {\tt then return} $\im\rho_d$; \\
    {\tt elif} $\dim\Ker\rho_d\cap\im\rho_d = 4$ {\tt then return}
    $\Ker\rho_d\cap\im\rho_d$; \\ 
    {\tt end if};
  \item {\tt if} $\dim\Ker\rho_d\cap\im\rho_d = 0$ {\tt then} \\
    \phantom{xx} $A_1 := \im\rho_d\cap\im\lambda_d$; \\
    \phantom{xx} $d_1 := $ zero divisor in $A_1$; \\
    \phantom{xx} {\tt return} $\im\rho_{d_1}$ ($d_1$ dealt here as an element in $A$); \\
    {\tt end if};
  \item $\mathfrak{R} := \Ker\rho_d\cap\Ker\lambda_d$; \\
    $A_d := A/\mathfrak{R}$; 
    let $\pi$ be the natural projection $A_d \to A_d/\mathfrak{R}$; \\
    $e := $ zero divisor in $A_d$; \\
    {\it // find $f'\in\pi^{-1}(e)$ such that $\dim\Ker\rho_{f'} = 4$} \\
      fix an element $f'_0\in\pi^{-1}(e)$ and a basis $(b_1,b_2,b_3,b_4)$ of 
      $\pi^{-1}(e) - f'_0$; \\
      $m := 1$; \\
      {\tt repeat} \\
      \phantom{xx} {\tt for} each $(c_1,c_2,c_3,c_4)$ such that 
        $\sum_i |c_i| = m$ {\tt do} \\
      \phantom{xxxx} $f' := f'_0 + \sum_i c_ib_i$; \\
      \phantom{xxxx} {\tt if} $\dim\Ker\rho_{f'} = 4$ 
        {\tt then return} $\Ker\rho_{f'}$; {\tt end if}; \\
      \phantom{xx} {\tt end for}; \\
      \phantom{xx} $m := m+1$; \\
      {\tt until false}.
\end{enumerate}

Note that the repeat-loop in the last step of the algorithm terminates.
Indeed, the preimage $\pi^{-1}(e)$ is a 4-dimensional affine space.
For each point in $\pi^{-1}(e)$ with integral coefficients 
it is tested whether it lies in the Zariski open set given by 
$\alpha_1\beta_4+\alpha_4\beta_1\ne\alpha_2\beta_3 + \alpha_3\beta_2$
($\alpha_i\in k$ and $\beta_i$ are indeterminates)
with respect to some coordinate system in $\pi^{-1}(e)$.
Clearly, the closed complement of this set does not contain a 4-dimensional lattice,
therefore a desired element is hit and the loop is left.

The last point to explain is the first step of the algorithm: 
finding any zero divisor. 
We start by finding a quadratic element, i.e.~an element $a$ such that 
the minimum polynomial $\mu_a(\xi)$ is irreducible quadratic.

\begin{lem}[Rowen]
  Let $c\in A$ have the minimum polynomial 
  $\mu_c(\xi) = \xi^4 + \alpha_2\xi^2 + \alpha_1\xi + \alpha_0$.
  Then for any factorization 
  $\mu_c(\xi) = (\xi^2 + a^\prime\xi+b^\prime)(\xi^2 + a\xi+b)$ 
  in $A[\xi]$ we have $[k(a^2):k] < 4$.
\end{lem}

\begin{proof}
  First denote $\nu(\xi) = (\xi^2 + a\xi+b)(\xi^2 + a^\prime\xi+b^\prime)
  = \xi^4 + \beta_2\xi^2 + \beta_1\xi + \beta_0$ with 
  $\beta_i$'s possibly in the algebra.
  By comparing coefficients in 
  $(\xi^2 + a\xi+b)\mu_c(\xi) = \nu(\xi)(\xi^2 + a\xi+b)$
  we obtain $\nu(\xi) = \mu_c(\xi)$.
  Further we basically follow the proof in ~\cite{rowen}, 
  where a slightly stronger assertion is proven for division algebras.

  By matching coefficients in $\mu_c(\xi)$ we get $a^\prime = -a$ and
  \begin{eqnarray*}
    \alpha_2 &=& b + b^\prime + a^\prime a =  b + b^\prime - a^2, \\
    \alpha_1 &=& a^\prime b + b^\prime a = -ab + b^\prime a, \\
    \alpha_0 &=& b^\prime b = b b^\prime 
      \quad\textrm{(from } \alpha_0 = \beta_0 \textrm{)}.
  \end{eqnarray*}

  {\em Case 1:} $ab = ba$.
  Then also $ab' = b'a$ and $\alpha_1 = (b^\prime - b)a$.
  It follows that
  $\alpha_1^2 = ((b^\prime + b)^2 - 4\alpha_0)a^2 = 
  ((a^2 + \alpha_2)^2 - 4\alpha_0)a^2 =
  (a^2)^3 + 2\alpha_2(a^2)^2 + (\alpha_2^2 - 4\alpha_0)a^2$,
  so $[k(a^2):k]\le 3$.

  {\em Case 2:} $ab\ne ba$.
  By multiplying $\alpha_2 = b + b^\prime - a^2$ by $b$ from left and right 
  we get $a^2 b = b a^2$. Therefore $a\notin k(a^2)$, so 
  $[k(a^2):k] < [k(a):k] \le 4$. 
\end{proof}

The Rowen's lemma gives a recipe for finding a quadratic element as follows.
We start with an arbitrary noncentral $c\in A$. 
If $c$ happens to be a zero divisor or the minimum polynomial $\mu_c$
is reducible then we are done and do not need to continue in finding 
a quadratic element.
So we can assume now that neither $c$ is a zero divisor nor
$\mu_c$ is reducible.
Then the minimum polynomial $\mu_c$ is not cubic.
Indeed, the characteristic polynomial
$\chi_c(\xi) = \mu_c(\xi)\lambda(\xi)$, $\lambda\in k[\xi]$ linear,
and since very irreducible factor of the characteristic polynomials
divides also the minimum polynomial, 
we have $\lambda(\xi)\mid \mu_c(\xi)$.
Hence if $c$ itself is not quadratic,
its minimum polynomial is irreducible of degree $4$. 
After applying a linear substitution eliminating 
the cubic term in $\mu_c(\xi)$ we may use Wedderburn's factorization 
theorem to construct a factorization as in Rowen's lemma.
If $a^2$ is not a zero divisor and the minimum polynomial of $a^2$ 
is not reducible (in which case we would be done), then there are 
3 possibilities left:
\begin{itemize}
  \item[(a)] $[k(a^2):k] = 2$, so $a^2$ is quadratic, or
  \item[(b)] $a^2\in k^*$, then $a$ is quadratic, since otherwise 
    the constructed factorization according to Rowen's lemma 
    would be a factorization over $k$, or lastly
  \item[(c)] $a^2 = a = 0$, then $\mu_c(\xi) = \xi^4 + \alpha_2\xi^2 + \alpha_0$
    and $c^2$ is quadratic.
\end{itemize}

In sequel we will need the well-known
\begin{thm}[Double Centralizer Theorem]
  Let $B$ be a central simple algebra over $k$ and 
  suppose that $C$ is a simple subalgebra of $B$. Then 
  \begin{itemize}
    \item[(i)]   $\C_B(C)$ is simple, 
    \item[(ii)]  $[C:k]\ [\C_B(C):k] = [B:k]$, 
    \item[(iii)] $\C_B(\C_B(C)) = C$. 
  \end{itemize}
\end{thm}

\begin{proof}
  cf.~\cite{pierce}, p. 232.
\end{proof}

Let $a\in A$ be a quadratic element, 
so $k(a)$ is a subfield of $A$ and $[k(a):k] = 2$.
Then by the Double Centralizer Theorem, 
the centralizer of $a$ in $A$ is a simple algebra
of dimension $8$ over $k$. 
The center of $\C_A(a)$ is the field $k(a)$, therefore $\C_A(a)$
can be understood as a central simple algebra over $k(a)$.
We denote this algebra by $A_2$.

\begin{lem}\label{le:A2}
  $A_2$ is isomorphic to $M_2(k(a))$.
\end{lem}

\begin{proof}
  Let $\mathcal{L}$ be a $4$-dimensional left ideal in $A$.
  If we understand $A_2 = \C_A(a)$ as an algebra over $k$, 
  then $\mathcal{L}$ is also a $4$-dimensional $A_2$-module.
  Let $0\ne b_1\in\mathcal{L}$ and let $b_2\in\mathcal{L}$ be
  such that $b_2\notin \C(A_2)b_1$. Then $b_1, b_2$ is a basis
  of $\mathcal{L}$ over $k(a)$. 
  So we have a $2$-dimensional $A_2$-module over $k(a)$, 
  where $A_2$ is now taken to be an algebra over $k(a)$.
  Since $A_2$ acts faithfully on $\mathcal{L}$,
  this gives an embedding of $A_2$ into $M_2(k(a))$.
  Now the assertion of the Lemma follows from $[A_2:k(a)] = 4$.
\end{proof}

After finding the algebra $A_2$ as the centralizer of a quadratic 
element $a$, we write $A_2$ as a cyclic algebra over $k'=k(a)$,
so we find a cyclic element $c\in A_2$ and $u^\prime$ such that 
$(u^\prime)^2 = \gamma\in (k')^*$.
By Proposition~\ref{the:normeq} there is $s\in k'(c)$ such that
$s\sigma(s) = 1/\gamma$.
Then also $u=su^\prime$ is a principal generator of $A_2$ over $k'(a)$
and moreover $u^2 = 1$.
We found an element $u\in A$ with the reducible minimum polynomial 
$\mu_u(\xi) = \xi^2 - 1$, therefore $u+1$ is a zero divisor.

\vskip0.4cm
\noindent
\begin{tabular}{ll}
  {\sc Function:} & {\tt FindZeroDivisor}\\
  {\sc Input:}    & $A$ -- an associative algebra isomorphic to $M_4(k)$.\\
  {\sc Output:}   & $d$ -- a zero divisor in $A$.
\end{tabular}
\begin{enumerate}
  \item {\it // find a quadratic element $a\in A$} \\
    $c$ := a noncentral element in $A$; \\
    $c := c + c_3/4*1$, where $c_3$ is the cubic coefficient in $\mu_c$; \\
    {\tt if} $\deg\mu_c = 2$ {\tt then} $a := c$; \\
    {\tt elif} $\mu_c(\xi) = \xi^4 + \alpha_2\xi^2 + \alpha_0$ {\tt then} $a := c^2$;
      {\it // case (c)} \\
    {\tt else} \\
      \phantom{xx} find a factorization 
        $\mu_c(\xi) = (\xi^2 + a^\prime\xi+b^\prime)(\xi^2 + a\xi+b)$ over $A$;\\
        \phantom{xx} {\tt if} $a^2\notin k$ {\tt then} $a := a^2$; {\tt end if};
          {\it \ // cases (a) and (b)}\\
    {\tt end if};
  \item $A_2$ := the centralizer of $a$ in $A$ 
    regarded as a four-dimensional central algebra over $k(a)$.
  \item {\it \ // write $A_2$ as a cyclic algebra over $k(a)$} \\
    $c$ := a noncentral element in $A_2$; \\ 
    find $\sigma(c)\in k(a,c)$ such that $\mu_c(\xi) = (\xi-c)(\xi-\sigma(c))$; \\
    $u^\prime$ := a nonzero solution of the linear system 
    $u^\prime c = \sigma(c)u^\prime$; \\
    $\gamma := (u^\prime)^2$; {\it // ($\gamma\in k(a)$)} 
  \item {\tt if} $\exists r\in k(a)$ such that $r^2 = \gamma$ {\tt then} 
    $u := u^\prime/r$; \\
    {\tt else} \\
      \phantom{xx} $s$ := a solution of the norm equation 
        $s\sigma(s) = 1/\gamma$; {\it // ($s\in k(a)$)}; \\
      \phantom{xx} $u := su^\prime$; \\
    {\tt end if};
  \item {\tt return} $u+1$.
\end{enumerate}
\vskip0.3cm

As in the case of algebras of degree 3, for each element popping up during
the computation it is tested whether it already accidentally gives a zero divisor.
Any kind of the zero divisor would lead to avoiding solving the norm relative equation,
which is the most expensive step in the algorithm.

\nocite{hartshorne}  
\bibliographystyle{alpha}

\end{document}